\documentclass[reqno,12pt,pdfsync]{amsart}

\makeindex
\usepackage{amsmath}
\usepackage{latexsym} 
\usepackage{amscd}
\usepackage{amsfonts}
\usepackage{pb-diagram}
\usepackage[ps,dvips,all]{xypic}
\usepackage[mathcal,mathscr]{eucal}
\usepackage{amsthm}
\usepackage{epsfig}
\usepackage{pb-diagram}
\usepackage{lamsarrow}  
\usepackage{pb-lams}
\usepackage{epsfig}
\usepackage{bbm}
\usepackage[colorlinks]{hyperref}
\usepackage{url}
\usepackage{pifont}
\usepackage{amsbsy}
\usepackage{graphicx}
\usepackage{epstopdf}
\usepackage{times}
\newcounter{theoremcount}
\addtocounter{theoremcount}{1}
\newtheorem{theorem} {Theorem}[section]
\newtheorem{prop}[theorem]{Proposition}
\newtheorem{lemma}[theorem]{Lemma}
\newtheorem{corollary}[theorem]{Corollary}
\numberwithin{equation}{section}

\theoremstyle{definition}
\newtheorem{remark}[theorem]{Remark}

\theoremstyle{definition}
\newtheorem{definition}[theorem]{Definition}

\theoremstyle{definition}
\newtheorem{example}[theorem]{Example}

\DeclareMathOperator{\Sym}{\mathscr{S}ym}
\DeclareMathOperator{\Range}{Range}
\DeclareMathOperator{\Ran}{Ran}
\DeclareMathOperator{\Hom}{Hom}
\DeclareMathOperator{\Isom}{Isom}

\DeclareMathOperator{\Cay}{Cay}

\DeclareMathOperator{\End}{End}
\DeclareMathOperator{\Ker}{Ker}
\DeclareMathOperator{\Ad}{Ad}
\DeclareMathOperator{\Ri}{Ri}

\DeclareMathOperator{\Tr}{Tr}
\DeclareMathOperator{\SFred}{\mathscr{S}Fred}
\DeclareMathOperator{\BFred}{BFred}

\DeclareMathOperator{\Sch}{Sch}
\DeclareMathOperator{\PD}{PD}

\DeclareMathOperator{\ind}{ind}
\DeclareMathOperator{\even}{even}
\DeclareMathOperator{\odd}{odd}

\DeclareMathOperator{\ch}{ch}

\DeclareMathOperator{\id}{id}

\DeclareMathOperator{\codim}{codim}

\DeclareMathOperator{\Mas}{\mathscr{M}as}

\DeclareMathOperator{\Lag}{\mathscr{L}ag}
\DeclareMathOperator{\Gr}{Gr}
\DeclareMathOperator{\Real}{Re}

\DeclareMathOperator{\Imag}{Im}
\DeclareMathOperator{\R}{\mathscr{R}}

\newcommand\bR{{\mathbb R}}

\newcommand\bZ{{\mathbb Z}}

\newcommand{\eC}{\EuScript{C}}

\newcommand{\eE}{\EuScript{E}}
\newcommand{\eF}{\EuScript{F}}
\newcommand{\eG}{\EuScript{G}}
\newcommand{\eH}{\EuScript H}

\newcommand{\eN}{\EuScript{N}}

\newcommand{\eS}{\EuScript{S}}
\newcommand{\eT}{\EuScript{T}}
\newcommand{\eU}{\EuScript{U}}

\newcommand{\eX}{\EuScript{X}}
\newcommand{\eY}{\EuScript{Y}}

\newcommand{\ra}{\rightarrow}
\newcommand{\hra}{\hookrightarrow}
\newcommand{\Lra}{{\longrightarrow}}

\newcommand{\eSy}{\EuScript{S}ymp}

\newcommand{\Bom}{\boldsymbol{\Omega}}
\newcommand{\bom}{\boldsymbol{\omega}}

\begin{document}

\title{\bf Localization formulae in odd $K$-theory}
\author{ Daniel Cibotaru}

\begin{abstract} 
We describe a class of real Banach manifolds, which classify $K^{-1}$. These manifolds are Grassmannians of (hermitian) lagrangian subspaces in a complex Hilbert space. Certain finite codimensional real subvarieties described by incidence relations define  geometric representatives for the generators of the cohomology rings of these classifying spaces. Any family of self-adjoint, Fredholm operators parametrized by a closed manifold comes with a map to one of these spaces.  We use these Schubert varieties to describe the Poincare duals of the   pull-backs to the parameter space of the  cohomology ring generators.  The class corresponding to the first generator is the spectral flow.
\end{abstract}
\maketitle

\noindent
\tableofcontents

\section{Introduction}

In this paper we answer a question posed by Isadore Singer in the mid 80's concerning the complex odd 
K-theoretic functor $K^{-1}$. 

The classifying spaces for $K^{-1}$ have various homotopically equivalent realizations. There is the classical description as the unitary group $U(\infty)$ and there is the Atiyah-Singer realization as the space $\BFred_*(H)$, a certain component of the space of \textit{bounded}, self-adjoint, Fredholm operators on a separable complex Hilbert space $H$.

The cohomology ring $H^*(U(\infty),\mathbb{Z})$  is an exterior algebra with canonical generators $(x_k)_{k\geq 1}$, deg $x_k=2k-1$. The degree one generator $x_1$ has a very useful geometric interpretation. If $A:S^1\rightarrow \BFred_*(H)$ is a continuous loop of Fredholm operators then the Poincare dual of the class $A^*x_1\in H^1(S^1,\mathbb{Z})$ can be represented by a cycle supported on the degeneracy locus $\{\theta\in S^1~|~\Ker{A_\theta}\neq 0\}$. The integer $\int_{S^1}A^*x_1$ is called the spectral flow of the family and, under generic conditions, can be described as a count with sign of the zero eigenvalues of the family. 

Singer asked for a similar description of the classes $A^*x_k$, where $A:M\rightarrow\BFred_*(H)$ is a continuous map of self-adjoint Fredholm operators.  Moreover it is desirable to design an approach that deals with families of \textit{unbounded} Fredholm operators directly without passing to the associated bounded operators via functional calculus.

Our approach to this problem is based on symplectic techniques and has the added bonus that it provides an elegant way of dealing with unbounded operators as well. The dissertation is roughly divided in three parts. In the first part we describe several smooth models for the classifying space of $K^{-1}$ while in the second part we describe various finite codimensional, cooriented stratified spaces that determine all the cohomology classes corresponding to products of the canonical generators $x_i$.

In order to describe the main results we introduce a bit of terminology.

 On the direct sum of a complex Hilbert space with itself $\hat{H}:=H\oplus H$, there is a natural extra complex structure: $J=\left(\begin{array}{cc}  0& 1 \\ -1 &0 \end{array}\right)$.  A lagrangian subspace $L\subset \hat{H}$ is a subspace which is taken by $J$ isomorphically to its orthogonal complement, $JL=L^\perp$.  For example the graph of every (closed) self-adjoint operator, bounded or unbounded is a lagrangian. 
   
 The Grassmannian of such spaces, denoted $\Lag(\hat{H})$, can be turned into a \textit{real} Banach manifold. In fact the space of bounded, self-adjoint operators on $H$ which we denoted by $\Sym (H)$ embedds as an open dense subset of $\Lag$ via the graph map:
 \[ \Sym(H)\ni A\rightarrow \Gamma_A\in\Lag(\hat{H}), ~~~~\Gamma_A:=\{(x,Ax)~|~x\in H\}
 \]
 
 On the other hand $\Sym(H)$ embedds into the space of unitary operators, $\mathscr{U}(H)$ via the Cayley transform. Our first result, which is the generalization to infinite dimensions of a result by Arnold \cite{A} says that $\Lag$ is diffeomorphic with $\mathscr{U}(H)$ and the diffeomorphism is  constructed by extending the Cayley transform to the whole $\Lag$. (see Theorem \ref {Arnold} for details)
 
   By a famous result of Kuiper,  $\mathscr{U}(H)$ and hence $\Lag$ is contractible. Nevertheless an open subset  of this space $\Lag^-$, has the homotopy type of the inductive limit of unitary groups, thus  classifying  $K^{-1}$. By definition, $\Lag^-$ is the set of all lagrangians which are Fredholm pairs with a fixed one, namely the vertical space, $H^-:=0\oplus H$. It turns out that, inside $\Lag^-$, there is a whole zoo of classifying spaces, $\Lag_{\mathscr {I}}$ which correspond via the mentioned diffeomorphism to the Palais unitary groups $U_{\mathscr{I}}(H)$, modelled on  two-sided, symmetrically normed ideals $\mathscr{I}$.  Examples of such ideals are the Schatten class operators. The reader can find an intrinsic characterization of $\Lag_{\mathscr{I}}$ in section \ref {Lag Gr}.

 The main technique used in proving that $\Lag^-$ is a classifying space for $K^{-1}$ is (linear) \textit {symplectic reduction}. More precisely,  for every closed subspace $W\in H^-$ let $H_W$ be the orthogonal complement of $W\oplus JW$ in $\hat{H}$. The symplectic reduction is a map:
 \[ \mathscr{R}:\Lag^-\rightarrow \Lag{H_W}
 \]
 which is continuous on a certain open subset denoted $\Lag^W$  of $\Lag^-$. It is, in fact, diffeomorphic with a vector bundle over $\Lag{H_W}$.
   To investigate the homotopy type of $\Lag^-$ we  take a complete, decreasing flag of $H^-$ 
 \[ H^-=:W_0\supset W_1\supset W_2\supset \ldots
 \]
 and this flag determines a filtration of $\Lag^-$ by open subsets $\Lag^{W_i}\subset\Lag^{W_{i+1}}$. Each $\Lag^{W_i}$ is homotopy equivalent with $\Lag{H_{W_i}}$, hence with $U(i)$ by the finite version of Arnold's theorem. The limit when $i->\infty$ of $\Lag^{W_i}$, which is $\Lag^-$, has the homotopy type of $U(\infty)$. (see  Section \ref {Symp red} and Theorem \ref{weak homotopy equivalence})

  The Atiyah-Singer space, $\BFred_*$ embedds in $\Lag^-$, simply by associating to an operator its \textit{switched} graph. 
   \[   \{(v,Av)~|~v\in H\}~~\mbox{ graph }\quad \leftarrow \quad A \mbox { operator } \quad\rightarrow \quad\{(Av,v)~|~v\in H\}\quad \mbox { switched graph}
 \]
 We prove that this map is a weak homotopy equivalence.  In this way, the classical index, as defined by Atiyah and Singer of a family of self-adjoint, Fredholm operators is simply the homotopy class of the family of switched graphs. This definition extends to the unbounded case.
 
 In the second part of this paper we build geometric representatives for the generators of the cohomology ring of $\Lag^-$ which we identify with some "canonically" defined cohomology classes of $\Lag^-$. These  canonical classes are  described as follows.  The natural inclusions $i_n:U(n)\rightarrow\Lag^-$ induce isomorphisms in cohomology:
 \[ H^k(\Lag^-)\rightarrow H^k(U(n)) , \forall k< 2n-1
 \]
 The cohomology ring of $U(n)$ is an exterior algebra over $\mathbb{Z}$ with $n$-generators $x_i(n)$. These generators are obtained by transgressing the Chern classes $c_i(n)$ of the universal rank $n$ bundle over $S^1\wedge U(n)$. Therefore there exist \textit{unique} cohomology classes $x_i\in H^{2i-1}(\Lag^-)$ that pull-back via the inclusion maps $i_n$ to these generators $x_i(n)$.
 
 The geometrical representatives of the classes $x_i$ are build from some finite codimensional stratified spaces of $\Lag^-$. These spaces are the analogues of the Schubert varieties on the usual finite Grassmannian with the crucial difference that the strata are  \textit{real} manifolds.
 
   Here is briefly the construction. The set
  \[ Z_k:=\{L~|~\dim{L\cap W_{k-1}}=1\}
  \] 
  is a cooriented, codimension $2k-1$ submanifold of $\Lag^-$. Its set of singularities, $\partial {Z_k}:=\overline{Z_k}\setminus Z_k$ has codimension  bigger than $2k+1$  in the ambient space. This implies that 
  \[H^k(\Lag^-)\simeq H^k(\Lag^-\setminus\partial Z_k)
  \]
  is an isomorphism. The submanifold  $Z_k$ is closed in $\Lag^-\setminus\partial Z_k$. The cohomology class, denoted $[Z_k,\omega_k]$, defined  by $Z_k$,  with the natural coorientation $\omega_k$ is the image of $1$ via the composition of maps:
  \[ H^0(Z_k)\simeq H^{2k-1}(\Lag^-\setminus \partial {Z_k},\Lag^-\setminus \overline {Z_k})\rightarrow H^{2k-1}(\Lag^-\setminus\partial {Z_k})\simeq H^{2k-1}({\Lag^-})
  \]
   Above,  the first map is  the Thom isomorphism   combined with excision and the second and the third map are the natural pull-backs. We have the following result (for details see  Theorem \ref{intersection}) 
   
   \medskip
   
   \noindent  
\textbf{Theorem A}  \textit{   The geometric cohomology class $[Z_k,\omega_k]$ coincides with the canonical class $x_k$.}\\

 Suppose now that $M$ is a closed manifold and suppose one has a smooth family of self-adjoint, Fredholm operators parametrized by $M$. Then taking their (switched) graphs one gets a smooth map $f:M\rightarrow\Lag^-$.  In transversal conditions, made precise in Section \ref {IntFor},  the preimage set $f^{-1}(\overline{Z}_k)$ is a stratified space with no singularities in codimension one with both a coorientation and an orientation on the top stratum. We denote by $[f^{-1}(\overline{Z}_k), f^*\omega_k]^*$ the cohomology class  and  respectively $[f^{-1}(\overline{Z}_k), f^*\omega_k]_*$ the Borel-Moore homology class it determines. We have the following result.
 
 \medskip
 \noindent
 \textbf {Localization Theorem}  \textit {Let $f:M\rightarrow \Lag^-$ be the map determined by a  smooth family of self-adjoint operators. Then}
 \[ \tag{$a$}
 f^*(x_k)= [f^{-1}(\overline{Z}_k), f^*\omega_k]^*,
 \]
 \[\tag{$a^*$} \PD  f^*(x_k)=[f^{-1}(\overline{Z}_k), f^*\omega_k]_*,
 \]
 \textit{where $\PD$ denotes  Poincar\'e duality.}
\qed \\

    In the case when $M$ has complementary dimension to $Z_k$, i.e. $\dim{M}=2k-1$, then the preimage consists of a bunch of points with signs.  The relevance of these $0$-cycles for index theory is that they represent up to a fixed constant the Poincar\'e dual of the $2k-1$ th component of the cohomological index (see Section \ref{Chern} and Theorem \ref{chern and poincare dual}).
    In the case when $M=S^{2k-1}$ then the degree of the $0$-cycle,i.e. the sum of the local intersection numbers determines completely the homotopy type of the map $f$. This number is always divisible by $(k-1)!$ by Bott divisibility theorem. (see Theorem \ref{homotopy type and cycles}).
    
      We devote the last part of this paper to find formulae for the local intersection numbers in terms of the differentials of the family of operators.  One of the difficulties in doing intersection theory with the Schubert varieties we described is finding a useful characterization of their normal bundle. In order to achieve this we used a general form of reduction. (see section \ref {Gen Red}). Arnold's  theorem comes to rescue in some key technical points in achieving this description. (see Proposition \ref{algebraic complement}) 
      
      In the case $k=1$ we recover the classical description of the spectral flow. For $k\geq 2$ as one might expect the local intersection numbers, unlike the spectral flow depend not only on the variation of the eigenvalues but also on the variation of the eigenspaces. An example in that direction is Proposition \ref{strong transversality}.
            
      In the last Chapter of the paper we also discuss what it means for a family of operators to be smooth. We use our criterion of differentiability on the  universal family of self-adjoint elliptic operators, parameterized by the group $U(N)$. This family associates to a unitary matrix, the Dirac operator on the $\mathbb{C}^N$ vector bundle on a circle obtained from the unit interval by gluing the end vector spaces via the corresponding unitary operator. Doing symplectic reduction on this example we get that the associated family of graphs is homotopy equivalent with the inclusion $U(N)\hookrightarrow \Lag^-$ reproving a result that first appeared in  \cite {KM}.  This universal family is key to the proof that Atiyah-Singer classifying space is homotopy equivalent with $\Lag^-$.
      
      In an Appendix we include a collection of known facts  that we deemed important for understanding how one builds cohomology classes out of stratified spaces endowed with a coorientation and having no singularities in codimension $1$.

   It is hard to overestimate the influence that my advisor, Prof. Liviu Nicolaescu had on this work.  I benefitted greatly from our conversations, questioning and constant encouragement and I am more than grateful to him. Over the years I have enjoyed mathematical discussions with many people. I learned a lot from Prof. Stolz, Prof. Williams, Prof. Evans, Prof. Hind, Prof. Xavier, Prof. Hall. My colleagues and friends Florin Dumitrescu, Iulian Toader, John Harper, Allegra Berliner Reiber, Katie Grayshan, Mark Collarusso, Ryan Grady, Fernando Galaz-Garcia, Inanc Baikur, Stuart Ambler have shared with me ideas and enthusiasm. Their friendship is priceless. 
   
   The good things in this paper are dedicated to my parents.

   \section{Classifying spaces for  odd $K$-theory} 
\setcounter{equation}{0} 
In this chapter we will introduce our main objects of study. These  are various real Banach manifolds classifying for $K^{-1}$ that can be described either as subsets of the unitary group of a complex Hilbert space or as subsets of the infinite Lagrangian Grassmannian determined by the same Hilbert space. The connection between the two different types of classifying spaces is provided by a generalization to  infinite dimensions of a theorem of Arnold. Symplectic reduction turns out to be a useful technique that reduces many of our questions to their finite dimensional counterpart.

\subsection{The infinite unitary group}

This section is designed to recall some well-known facts about spaces of unitary operators and to introduce notation and terminology that will be used throughout.

The space of $\textit{\textbf{bounded}}$ operators on a fixed, complex, separable Hilbert space $H$ is denoted by $\mathscr{B}$ or $\mathscr{B}(H)$.

The group of \textit{\textbf{ unitary operators}} $\mathscr{U}:=\{U\in \mathscr{B}~|~UU^*=U^*U=I\}\subset\mathscr{B}$  has the structure of a Banach manifold modelled on the space of \textit {\textbf {self-adjoint, bounded}} operators $\Sym:=\{A\in\mathscr{B}~|~A=A^*\}$.\index{$\Sym$} 

The charts are given by the Cayley transforms.  For a fixed $U_0\in \mathscr{U}$ the set
 \[\mathcal{A}_{U_0}:=\{U\in \mathscr{U},~1+UU_0^{-1}~  \mbox{ invertible }\} \] is open and the map:
\begin{center}
$U\rightarrow i(1-UU_0^{-1})(1+UU_0^{-1})^{-1}$
\end{center}
is a homeomorphism $\mathcal{A}_{U_0}\simeq \Sym$. In fact, $\mathscr{U}$ is a Banach-Lie group since multiplication and taking inverse (actually adjoint) are obviously differentiable maps.
\begin{definition} The open set $\mathcal{A}_{U_0}$ together with the Cayley transform  is called the $\textit{\textbf{Arnold chart}}$ \index{Arnold chart} around $U_0$.
\end{definition}
A famous result of Kuiper, \cite {Ku}, says that the group of unitary operators on a Hilbert space has trivial topology. Nevertheless, certain subspaces are more interesting.

We recall the following
\begin{definition} A bounded operator $T:H\rightarrow H$ is called \textit{\textbf {Fredholm}}\index{Fredholm operator} if the dimension of both its kernel and cokernel are finite dimensional spaces and if its image is a closed subspace of $H$.

An unbounded operator $T:D(T)\subset H\rightarrow H$ is Fredholm if it is closed,i.e. the graph is closed in $H\oplus H$, its image is closed and its kernel and cokernel are finite dimensional spaces.
\end{definition}

Let us now introduce the following subset of the group of unitary operators on $H$ \index{$\mathscr{U}_{-1}$}
\begin{center}
 $\mathscr{U}_{-1}:=\{U:H\rightarrow H~|~UU^*=U^*U=I,~ 1+U$ is Fredholm $\}.$
\end{center}
Another way of saying that $1+U$ is Fredholm is $-1\notin\sigma_{ess}(U)$, that is $-1$ is not an essential spectral value of $U$.  We will see in the next sections that $\mathscr{U}_{-1}$ is an open subset of the space of unitary operators  $\mathscr{U}$ that has the homotopy type of $U(\infty)$. 
 This space also appears in \cite{BLP}.

\begin{remark} Notice  that $\mathscr{U}_{-1}$ is not a group since if $U$ is a unitary operator such that $\pm 1$ are not spectral values then so is $-U^{*}$ and then $1+U(-U^*)=0$.  
\end{remark}

\begin{definition} Let $\mathscr{K}$ be the ideal of compact operators in $\mathscr{B}$ and let $\mathscr{I}\subset\mathscr{K}$ be a non-trivial, two-sided, subideal with a topology at least as strong as the norm topology. The \textit {\textbf {Palais unitary group}} \index{Palais unitary group} $\mathscr{U}_{\mathscr{I} }$ of type $\mathscr{I}$ is the subgroup of $\mathscr{U}$ consisting of operators of type $U=I+T$ where $I$ is the identity operator and $T\in \mathscr{I}$.
\end{definition}
 Let us note that any two-sided ideal is$~^*$-closed. \\

 The topology on $\mathscr{U}_{\mathscr{K}}$ is the norm topology, the topology on $\mathscr{U}_{\mathscr{I} }$ is the topology induced by $\mathscr{I}$ on its subset $\mathscr{U}_{\mathscr{I} }-I$. 
 
 Palais \cite {Pa} has shown that the Palais unitary groups are classifying for odd $K$-theory, that is to say, they are homotopy equivalent with:
 \[ U(\infty):=\lim_{\rightarrow} {U(n)}
 \]
 where $U(n)$ is the unitary group on $\mathbb{C}^n$. We will reprove his result, using different methods in section $\ref {Symp red}$.
 
The two-sided ideals have been classified by Calkin in \cite {Ca}. For a quick description of the relevant aspects of the theory a good reference is \cite {S}. We will be content to describe the ideals of Schatten class operators.(see also \cite {Q})

\begin {definition} Fix a number $p, ~p\in[1,\infty]$. A compact operator $K$ is of the \textit{\textbf{Schatten class}} \index{Schatten class} $p$  if $\Tr (K^*K)^p<\infty$. An operator of Schatten class $p=\infty$ is just a compact operator. 
\end{definition}

 We will use the notation $\Sch^p$ for the set of all operators of Schatten class $p\in[1,\infty]$. These are naturally Banach spaces with the norm: 
\begin{center}
$\|K\|_p=(\Tr (K^*K)^p)^{1/p} $ for $p\in[1,\infty)$
\end{center}
or the operator norm in the case $p=\infty$.  Let us notice that for a compact operator the following relation about their spectra is true $\sigma(K^*K)=\sigma(KK^*)$. This is because both have $0$ in their spectrum and it is easy to see that an eigenvalue for $T^*T$ is an eigenvalue for $TT^*$ as well and vice-versa.  It follows that the spaces $\Sch^p(H)$ are $*$-invariant. Also we have the following important inequality .
\begin{center}
$\|TK\|_p \leq\|T\|_{p_1}\|K\|_{p_2}$ whenever $p^{-1}\leq p_1^{-1}+p_2^{-1}$ 
\end{center}
which holds for all $p_1,p_2\in[1,\infty]$, $K\in \Sch^{p_2}$ and $T\in \Sch^{p_1}$ and moreover it holds for $T\in \mathscr{B}$ in which case the norm $\|\cdot\|_{p_1}$ is the operator norm. This turns $\Sch^p$ into  a closed $~^*$ subideal of $\mathscr{K}$. It is also saying that the Schatten ideals increase with $p$.

 The cases $p=1$ of trace class operators and $p=2$ of Hilbert-Schmidt operators are most likely the more familiar examples. The spaces $\Sch^p(H)$ should be taken as abstract analogues of the usual $L^p$ spaces. Indeed one other similar feature is the following.  Let $\hat{p}:=\frac{p}{p-1}$ for $1<p<\infty$, $\hat{1}=\infty$ and $\hat{\infty}=1$. Then the following duality relations hold:
\begin{center}
 $(\Sch^p)^*$ = $\Sch^{\hat{p}}$ for $1<p<\infty$ \\  
$(\Sch^\infty)^*=(\mathscr{K})^*=\Sch^1$
\\ $(\Sch^1)^*=\mathscr{B}$\\
\end{center}

The Palais groups are all Banach-Lie groups modelled on the space of self-adjoint operators. 
\begin{center}
  $\Sym_{\Sch^p}:=\{A\in \Sch^p~|~A=A^*\}$.
\end{center} 

\begin{remark} We clearly $\mathscr{U}_{\mathscr{I}}\subset \mathscr{U}_{-1}$ for every ideal $\mathscr{I}$ since the only essential spectral value of $U=I+K$ is $1$. 
\end{remark}

\bigskip

\subsection{The complex Lagrangian Grassmannian} \label {Lag Gr}

Our main object of study is the space of all lagrangians on a Hilbert space, endowed with an extra complex structure. We will give first the main definitions. 
\begin{definition}
 Let $H$ be a separable, complex Hilbert space and let  $\hat{H}=H\oplus H$. We denote by $H^+$ \index{$\hat{H}$}the space $H\oplus 0$ and call it the \textit{\textbf{horizontal subspace}}\index{horizontal subspace} and by $H^-$  the space $0 \oplus H$ and call it the \textit{\textbf{vertical subspace}}\index{vertical space}.
 \end{definition}
 Let $J\hat{H}\rightarrow \hat{H}$ be the unitary operator which has the block decomposition relative $\hat{H}=H\oplus H$ 
 \[ J=\left(\begin{array}{cc}
 0 &1 \\
 -1 & 0
 \end{array}\right)
 \] 
 The essential properties of $J$ are $J=-J^*=-J^{-1}$ so we can think of $J$ as a complex structure on $H$.
 
  \begin{definition}  A subspace $L\subset \hat{H}$ is called \textit{\textbf{lagrangian}}\index{lagrangian} if $L$ 
 $JL=L^\perp$.
 
    The (hermitian) \textit{\textbf{Lagrangian Grassmannian}}, \index{Lagrangian Grassmanian} \index{ $\Lag$} $\Lag(\hat{H},J)$ or simply $\Lag$ is the set of all lagrangian subspaces of $\hat{H}$.
\end{definition}
\noindent
\begin{remark} Notice that $JL=L^\perp$ implies that $L$ is closed since $L^\perp$ is always closed.
\end{remark}

\noindent
\begin {remark} The notion of Lagrangian Grassmannian of a complex space already appears in the literature for example in \cite{FP} . We caution the reader that the  notions which we use here is not the same as the one used there.  The symplectic structure that underlines the definition of a lagrangian in our case, $\langle J(\cdot), \cdot \rangle$ is skew symmetric in the hermitian sense, that is $\langle J(x), y\rangle=-\overline{\langle J(y), x\rangle}$. As we will see, our $\Lag$ is only a real manifold. 
\end{remark}
\begin{example}
\begin{enumerate}
\item[i)] Each of the spaces $H^\pm$ is a lagrangian subspace and $JH^\pm=H^\mp$.
\item[ii)] Given a self-adjoint, Fredholm operator, bounded or unbounded $T:D(T)\subset H\rightarrow H$, its  \textit {\textbf {graph}}:
\[ \Gamma_{T}:=\{(v,Tv)~|~v\in D(T)\}
 \]
 and its \textit{\textbf {switched graph}}\index{switched graph}\index{$\tilde{\Gamma}$}:
  \[ \tilde{ \Gamma}_{T}:=\{(Tw,w)~|~ w\in D(T)\}
 \]
 are both lagrangian subspaces in $\hat{H}$.
 \end{enumerate}

\end{example}
 The  Lagrangian Grassmannian is naturally endowed with a topology as follows. To each lagrangian $L$  we associate the orthogonal projection $P_L\in \mathscr{B}(\hat{H})$ such that $\Ran{P_L}=L$. The condition that $L$ is a lagrangian translates into the  obvious  relation
 \[J P_L=P_{JL}J=P_L^{\perp}J=(1-P_L)J
 \]
 which if we let $R_L:=2P_L-1$ be the reflection in $L$ becomes
 \[ JR_L=-R_LJ
 \]
 It is easy to see that if $R$ is an orthogonal reflection that anticommutes with $J$ then $\Ker{(I-R_l)}$ is a lagrangian subspace. In other words we get a bijection
 \[\Lag \leftrightarrow \{R\in\mathscr{B}(\hat{H})~|~ R^2=1, \mbox{ } R=R^*, \mbox{ } RJ=-JR\}
 \]
  and so  $\Lag$ inherits a topology as a subset of $\mathscr{B}{(\hat{H})}$.
 
 The following lemma is well- known.

\begin{lemma} \label{projection comp} \begin {enumerate}
\item [(a)] If $L$ is a lagrangian and $S\in\Sym{(L)}$ is a self-adjoint operator then the graph of $JS:L\rightarrow L^\perp$ is a lagrangian as well.
\item[(b)] For a fixed lagrangian $L$, if $L_1$ is both lagrangian and the graph of an operator $T:L\rightarrow L^\perp$ then $T$ has to be of the type $JS$ with $S\in\Sym(L)$ self-adjoint. 
\end{enumerate}
\end{lemma}
\noindent
\textbf{Proof:} $a)$ The graph of $JS$ is closed since $JS$ is bounded. Then for every $v\in L$ one has 
\[ J(v,JSv)=J(v+JSv)=Jv+SJ (Jv)=SJw+ w=(SJw,w)
\] 
where $w=Jv\in L^\perp$. It is easy to see that
\[ \langle (v,JSv), ( SJ w, w)\rangle = 0 ,~~~~\forall v\in L, ~ w\in L^\perp
\]
 Hence $J\Gamma_{JS}\perp \Gamma_{JS}$. In order to finish the proof one has to show that 
 \[ J\Gamma_{JS}+\Gamma_{JS}=\hat{H}
 \]
 which comes down to showing that for every $(a,b)\in L\oplus L^\perp=\hat{H}$ the system
 \[\left\{\begin{array}{ccc}
 v+SJw&=&a \\
 JSv+w&=&b
 \end{array}
 \right .
 \]
 has a solution $(v,w)\in L\oplus L^\perp$. There is nothing easier than that:
 \[\begin{array}{ccc}
  v&=&(1+S^2)^{-1} (a-SJb)\\
  w&=&(1+JS(JS)^*)^{-1}(b-JSa)
  \end{array}
 \]
 $b)$ The orthogonal complement of the graph of $T$ is the switched graph of $-T^*$. Hence
 \[ J\Gamma_{T}=\tilde{\Gamma}_{-T^*}
 \]
 implies that $JTJ=T^*$ which is another way of saying that $JT$ is self-adjoint.
 \qed\\
 
It is a known fact that, in the finite dimensional case the  sets $\Sym{(L)}$ are mapped to open subsets of $\Lag$ around $L$, turning the Lagrangian Grassmannian into a manifold. The situation in the infinite dimensional case is identical. However we need the following important proposition.

\begin{prop}\label {Ar. chart} Let $L$ be a lagrangian space, $L\in \Lag $. The following are equivalent: 
\begin{enumerate}
\item [$(a)$] $L$ is the graph of an operator $JA:L_0\rightarrow L_0^\perp$ where $A\in \Sym(L_0)$.
\item [$(b)$] $L\cap L_0^\perp=\{0\}$ and $L+L_0^\perp$ is closed.
\item[$(b')$]  $\hat{H}= L\oplus L_0^\perp$.
\item[$(b'')$] $\hat{H}=L+L_0^\perp$.
\item [$(c)$] $R_L+R_{L_0}$ is invertible.
\end{enumerate}
\end{prop}

\noindent
\textbf{Proof:} 
$(a)\Rightarrow (b)$ Clearly if $L$ is the graph of an operator $L_0\rightarrow L_0^\perp$ then $L$ is a linear complement of $L_0^\perp$.

\noindent

$(b)\Rightarrow (b')\Rightarrow (b'')\Rightarrow b)$ We have the following equality:
\[ (L+L_0^\perp)^\perp=L^\perp\cap L_0=J(L\cap L_0^\perp)=\{0\}
 \] 
 and this  proves that $(b)\Rightarrow (b')$. Clearly $(b')$ implies $(b'')$ and $(b'')$ implies $(b)$ because if $z\in L\cap L_0^\perp$ then $Jz\perp L$ and $Jz\perp L_0$ and so $z=0$.\\

\noindent
$(b)\Rightarrow (c)$ It is easy to check that
\[ \Ker {(P_L-P_{L_0^\perp})} =L\cap L_0^\perp\oplus L^\perp\cap L_0 =\{0\}
\] 
which combined  with $R_L+R_{L_0}=R_L-R_{L_0^\perp}=2(P_L-P_{L_0^\perp})$ proves the injectivity.

 Part $(b')$ gives also $\hat{H}=L^\perp\oplus L_0$. This implies that:
  \[  \Range{~(R_L-R_{L_0^\perp})}=\Range{~(P_L-P_{L_0^\perp})} =L+L_0^\perp=\hat{H}\]
 which proves the surjectivity.
 
  To see the second equality pick first $z\in L$. Then $\hat{H}=L^\perp\oplus L_0$ implies that $z$ can be written uniquely as $z=-z^\perp+y$ with $z^\perp\in L^\perp$ and $y\in L_0$.  Therefore $(P_L-P_{L_0^\perp})(z+z^\perp)=z$. 
   
    Similarly given $y^\perp\in L_0^\perp$ there exists a unique $y\in L_0$ such that $y-y^\perp\in L^\perp$ and so $(P_L-P_{L_0^\perp})(y-y^\perp)=y^\perp$. \\
    
 \noindent
$(c)\Rightarrow (a)$ We show that  the restriction to $L$ of the projection onto $L_0$, $P_{L_0}|_L$ is an isomorphism. 
First  $\Ker{P_{L_0}|_L}=L\cap L_0^\perp$ and since $\Ker{(R_L+R_{L_0})}=L\cap L_0^\perp\oplus L^\perp\cap L_0$ one concludes that $L\cap L_0^\perp=\{0\}$. 

Surjectivity comes down to showing that the adjoint $(P_{L_0}|_L)^*$ is bounded below \cite{B}.   But $(P_{L_0}|_L)^*$ is nothing else but $P_{L}|_{L_0}$. 

For $x\in L_0$ one has the following string of equalities:
\[ 
\| P_L|_{L_0}(x)\|=\| P_LP_{L_0}(x)\|=1/4\| (R_L+1)(R_{L_0}+1)(x)\|=
\]
\[
1/4 \|(R_L+R_{L_0}-R_{L_0}+1)(R_{L_0}+1)(x)\|=1/4\|(R_L+R_{L_0})(R_{L_0}+1)(x)\|=
\]
\[
=1/2\|(R_L+R_{L_0})P_{L_0}(x)\|=1/2\|(R_L+R_{L_0})(x)\|
\]
and the bound below follows from the invertibility of $R_L+R_{L_0}$. 

It is clear  that $L$ is the graph of an operator $T:L_0\rightarrow L_0^\perp$, $T=P_{L_0^\perp}\bigr |_L\circ (P_{L_0}|_L)^{-1}$. This operator has to be of the type $JA$ with $A\in \Sym(L_0)$ by part $b)$ in the previous lemma. 
\qed\\

\begin{corollary}
 The set $\{L\in \Lag~|~ L$ is the graph of an operator $L_0\rightarrow L_0^\perp\}$  is an \textit {open} neighbourhood around $L_0$
\end{corollary}
\noindent
\textbf{Proof:} The invertibility of $R_L+R_{L_0}$  is an open condition. \qed

\begin{definition} For a fixed lagrangian $L$, the map  $\mathcal{A}_L:\Sym(L)\rightarrow \Lag$ which associates to an operator $S$ the graph of $JS$ is called the \textit{\textbf{Arnold chart}}\index{Arnold chart} around $L$. We will sometimes use the same notation, $\mathcal{A}_L$ to denote the image of this map in $\Lag$.
\end{definition}

The only ingredient missing from turning $\Lag$ into a Banach manifold modelled on $\Sym(H)$ is to make sure that the transition maps are differentiable. To see that this is indeed the case we pick  a  unitary isomorphism $U:L_0\rightarrow L_1$. Then
\[ U^\sharp=\left ( \begin{array}{cc}
U &0 \\
0 & JUJ^{-1}
\end{array}
\right )
\] 
is a unitary isomorphism of $\hat{H}$  written in block decomposition as a map \linebreak $L_0\oplus L_0^\perp\rightarrow L_1\oplus L_1^\perp$. 

Let $L\in\mathcal{A}_{L_0}\cap\mathcal{A}_{L_1}$, that is $L=\Gamma_{JT}=\Gamma_{JS}$ where $S\in\Sym{(L_1)}$ and $T\in\Sym{(L_0)}$.  Let $\tilde{S}=U^{-1}SU\in\Sym{(L_0)}$. We pick $v\in L_1$ and let $w=U^{-1}v\in L_0$. Then
\[ v+JSv=Uw+JSUw=Uw+JUJ^{-1}JU^{-1}SUw=U^\sharp(w+J\tilde{S}w)
\] 
So $\Gamma_{JS}=U^\sharp\Gamma_{J\tilde{S}}=\Gamma_{JT}$.  Therefore 

\[ T=-JP_{L_0^\perp}\circ U^\sharp\circ(I,J\tilde{S})\circ(P_{L_0}\circ U^\sharp\circ(I,J\tilde{S}))^{-1}\
\]
where $(I,J\tilde{S}):L_0\rightarrow \hat{H}$ is the obvious operator whose range is the graph of $J\tilde{S}$. The differentiability is now clear so we have just proved:
\begin{prop} The Arnold charts turn the Lagrangian Grassmannian $\Lag$ into a Banach manifold modelled on the space of self-adjoint operators $\Sym{(H)}$.
\end{prop}

Let us a give an application of what we did so far. We will need this computation later.

\begin{lemma}\label{projection derivative} Let $P:\Lag\rightarrow \Sym{(\hat{H})}$ be the map that associates to the lagrangian $L$ the orthogonal projection onto $L$, i.e., $P(L):=P_L$. Then the differential $d_LP:\Sym{(L)}\rightarrow \Sym{(\hat{H})}$ is given by the following expression relative $\hat{H}=L\oplus L^\perp$
\[
d_LP(\dot{S})=\left(\begin{array}{cc}
0 & \dot{S}J_L^{-1} \\
J_L\dot{S} &  0
\end{array}
\right)=\left(\begin{array}{cc}
0 & -\dot{S}J \\
J\dot{S} &  0
\end{array}
\right)
\]
where $J_L:L\rightarrow L^\perp$ is the restriction of $J$ to $L$. \end{lemma}
\noindent
\textbf{Proof:}  We need an expression for the projection $P_{\Gamma_{JS}}$ onto the graph $JS:L\rightarrow L^\perp$. That comes down to finding $v$ in the equations:
\[ \left\{ \begin{array}{ccc}
a& = & v-SJ^{-1}w \\
b & = & JS v +w
\end{array}
\right.
\]
where $a,v\in L$ and $b,w\in L^\perp$. We get 
\[ v= (1+S^2)^{-1}(a +SJ^{-1} b)
\]
and the projection has the block decomposition relative $L\oplus L^\perp$. 
\begin{equation} \label{projection expression}
 P_{\Gamma_{JS}}=
\left(\begin{array}{cc}
(1+S^2)^{-1} & (1+S^2)^{-1}SJ^{-1} \\
J(1+S^2)^{-1}S &  J(1+S^2)^{-1}S^2J^{-1}
\end{array}
\right)
\end{equation}
Differentiating this at $S=0$ we notice that the diagonal blocks vanish since we deal with even functions of $S$ and so the product rule delivers the result.
\qed\\

We will see later that the tangent space of $\Lag$ can be naturally identified with the "tautological" bundle 
\[ \mathscr{T}: =\{(L,S)\in \Lag(H)\times \mathscr{B}(H\oplus H)~|~ S\in \Sym(L)\}.\]  

The space $\Lag$ is not very interesting from a homotopy point of view and  in the next  section we will prove that it is  diffeomorphic with the unitary group $\mathscr{U}(H)$ and so it is contractible. To get something non-trivial  we restrict our attention to the  subspace of vertical, Fredholm lagrangians.  
\begin{definition} A pair of lagrangians $(L_1,L_2)$ is called a \textit{\textbf{Fredholm pair}} if the following two conditions hold
\[ \dim(L_1\cap L_2)<\infty  \quad \mbox{ and } \quad  L_1+L_2 \mbox{ is closed}.
\]

The Grassmannian of \textit{\textbf{vertical, Fredholm lagrangians}}\index{vertical, Fredholm lagrangian}   \index{$\Lag^-$} is 
\[ \Lag^{-}:=\{L\in \Lag ~|~(L,H^-) \mbox{ is a Fredholm pair}\}.
\]
\end{definition}

    Fredholm pairs have been studied before both from the point of view of the projections (\cite{ASS})  and from the point of view of closed subspaces of a linear space (\cite{N2}). We summarize the main definitions and properties from \cite{ASS}.
    
    \begin{definition} \begin{itemize}
    \item[(a)] A pair of orthogonal projections $P$ and $Q$ in a separable Hilbert space $H$ is said to be a \textit{\textbf{Fredholm pair}} \index{Fredholm pair} if the linear operator
    \[ QP:\Ran{P}\rightarrow \Ran{Q}
    \] 
    is Fredholm.  
    \item[(b)] A pair of closed subspaces $U$ and $V$ of $H$ is said to be a Fredholm pair if 
    \[ \dim{U\cap V}<\infty, \quad \dim{U^\perp \cap V^\perp}<\infty \quad \mbox{ and }  \quad U+V \mbox { closed }.
        \]
       \item[(c)] Two subspaces $U$ and $V$ are said to be \textit {\textbf{commensurable}} if $P_U-P_V$ is a compact operator.
       \end{itemize}
    \end{definition}
    
    When $U$ and $V$ are lagrangian subspaces the middle condition in the definition of a Fredholm pair is superfluous.
        
    \begin{prop} Let $(P,Q)$ be a pair of projections. Then the following statements are equivalent.
    \begin{itemize}
    \item[(a)] The pair $(P,Q)$ is a Fredholm pair.
    \item[(b)] The pair $(Q,P)$ is a Fredholm pair.
    \item[(c)] The operators $P-Q\pm 1$ are Fredholm.
    \item[(d)] The pairs of subspaces $(\Ran{P},\Ker{Q})=(\Ran{P}, (\Ran{Q})^\perp)$  and $(\Ran{Q}, \Ker{P})=(\Ran{Q}, (\Ran{P})^\perp)$ are Fredholm pairs. 
    \end{itemize} 
    \end{prop}
    \noindent
    \textbf{Proof:}  For the equivalence of the first three claims see Proposition 3.1  and Theorem 3.4 (a) in \cite{ASS}. 

\noindent    
$(c) \Rightarrow (d)$ Let $U=\Ran{P}$ and $V=\Ker{Q}=\Ran{(1-Q)}$. Then  
\[  U\cap V=\Ker{(1-P+Q)} \quad \quad \mbox{ and }  \quad \quad U^\perp\cap V^\perp=\Ker{(1-Q+P)}.\] 
    
 If $\Ran{(1-Q+P)}$ is closed then the following sequence of inclusions proves that $U+V$ is closed
\[ (U^\perp\cap V^\perp)^\perp=\Ker{(1-Q+P)}^\perp=\Ran{(1-Q+P)}\subset U+V\subset (U^\perp\cap V^\perp)^\perp
 \]
 \noindent
 $(d) \Rightarrow (a)$ We need to prove that the map:
 \[ P_{V^\perp}: U\rightarrow V^\perp
 \]
 is Fredholm. Its kernel can be identified with $U\cap V$ so it is finite dimensional. Denote by $W$ the subspace
 \[W= U+(U^\perp\cap V^\perp)= U+(U+V)^\perp
 \]
 Clearly $W$ is a closed subspace since $(U^\perp\cap V^\perp)$ is a closed subspace of $U^\perp$. Moreover $W+V=H$ so that the operator $P_{V^\perp}:W\rightarrow V^\perp$ is Fredholm because it is surjective and has finite dimensional kernel. The inclusion $U\hookrightarrow W$ also Fredholm and so $P_{V^\perp}$ is Fredholm as a composition of two Fredholm operators
 \[ U\hookrightarrow W \rightarrow V^\perp.
 \]\qed 
 
 \begin{prop}\label{commensurate} Suppose $(U,V)$ is a Fredholm pair of closed subspaces and that $W$ is another subspace commensurable with $V$. Then the pairs $(V^\perp, W)$ and $(U,W)$ are  Fredholm pairs.
 \end{prop}
 \noindent
 \textbf{Proof:} This follows from the previous proposition and Theorem 3.4 (c) in \cite{ASS}.\qed
 \\
 
   Let $P^{\pm}|_L$\index{$P^{\pm}$} be the orthogonal projections on $H^\pm$ restricted to the lagrangian $L$ and $P_L|_{ H^\pm }$ be the projection on $L$ restricted to $H^\pm$. The following is just a corollary of the definitions and the first of the previous propositions.
 
 \begin{lemma} \label {char proj} The set of vertical, Fredholm lagrangians,  $\Lag^{-}$, coincides with the set 
 \[\{L\in \Lag ~|~P^+|_L \mbox{ is Fredholm of index } 0\}.
 \]
\end{lemma}
 
 \bigskip

\subsection{Arnold's Theorem}

 In this section we generalize a result by Arnold (\cite {A}, see also \cite {N1}) to infinite dimensions. In his article, Arnold showed that the finite Lagrangian Grassmannian $\Lag(N)\subset \Gr(N,2N)$ is diffeomorphic to the unitary group $U(N)$.  
 
 We introduce now the main suspects.  Consider the $\mp i$ eigenspaces of $J$, i.e. $\Ker{(J\pm i)}$ and let
 \[\Isom{(\Ker{(J+i)},\Ker{(J-i)})}
 \] 
 be the group of Hilbert space isomorphisms between the two eigenspaces. To each lagrangian $L$ we associate the restriction of the  reflection $R_L$ to $\Ker{(J+i)}$. Since $R_L$ anticommutes with $J$ we get a well-defined \textit{\textbf {reflection map}}
 \[ R_{-}:\Lag\rightarrow \Isom{(\Ker{(J+i)},\Ker{(J-i)})},\quad \quad R_{-}(L)= R_L\bigr|_{\Ker{(J+i)}}\]

 Notice that
  \[\label{reflection} R_L=\left (\begin{array}{cc} 
 0 &  R_{-}(L)^* \\
 R_{-}(L) & 0
 \end{array}  \right)
 \]
 relative to the decomposition $\hat{H}=\Ker{(J+i)}\oplus \Ker{(J-i)}$. On the other hand, to each isomorphism $T\in\Isom{(\Ker{(J+i)},\Ker{(J-i)})}$ we can associate its graph $\Gamma_{T}$ which is a subspace of $\hat{H}$. It is, in fact, a lagrangian. Indeed 
 \[ w\in J\Gamma_T \Leftrightarrow w=J(Tv+v)=-iTv+iv=z-T^*z\in \tilde{\Gamma}_{-T^*}\]
  for some $v\in \Ker{(J+i)} $, with  $z=-iTv\in \Ker{(J-i)}$. It is standard that the switched graph $\tilde{\Gamma}_{-T^*}$ is the orthogonal complement of $\Gamma_T$. Hence, we get  a second well-defined  \textit{\textbf{graph map}}
  \[ \Gamma_{-}: \Isom{(\Ker{(J+i)},\Ker{(J-i)})}\rightarrow \Lag,\quad \quad  \Gamma_{-}(T)=  \Gamma_T
  \]
  
  Notice that given a lagrangian $L$ there are canonical Hilbert space isomorphisms 
    \[  \phi_{\mp}(L):L\rightarrow \Ker{(J\pm i)}, \quad \quad \quad \phi_{\mp}(L): \quad  v\mapsto 1/\sqrt{2} ( v\pm iJv)
  \]
\noindent
\textbf{Notation:} When there is no possibility for confusion we will use $\phi_{\mp}:=\phi_{\mp}(L)$.
 
  Every isomorphism $T\in\Isom{(\Ker{(J+i)},\Ker{(J-i)})}$ comes from a unitary operator $U_T\in \mathscr{U}(L)$. 
  \[ \xymatrix{  L\ar_{\phi_{- }}[d] \ar^{U_T}[r] & L \ar^{\phi_{ +}}[d] \\
  \Ker{(J+i)} \ar^{T}[r]  & \Ker{(J-i)}
  }
  \]
  
It is straightforward to see that the graph,  $\Gamma_T$ is expressed in terms of the unitary operator $U_T$ as the set
\[ \Gamma_T=\{(1+U_T)v+iJ(1-U_T)v~|~v\in L\}
\]

For each lagrangian $L$, we will call the \textit{\textbf{Cayley graph map}} the following application
 \[ \mathscr{C}_L : \mathscr{U}(L) \rightarrow \Lag\]
 \[ \mathscr{C}_L(U): =\Gamma_{-}(\phi_{ i}\circ U\circ (\phi_{-})^{-1})=\Ran\{L\ni v\mapsto (1+U)v+iJ(1-U)v\in\hat{H}\}
 \]

We have the following 
  
 \begin{theorem}[\textbf{Arnold}]\index{Arnold Theorem}  \label{Arnold}
  \begin{itemize}
 \item[(a)] The reflection and the graph maps, $R_{-}$ and $\Gamma_{-}$ are inverse to each other. 
 \item[(b)] For every lagrangian $L$, the Cayley graph map, \index{Cayley graph map} $\mathscr{C}_L:\mathscr{U}(L)\rightarrow \Lag$ is a diffeomorphism of real Banach manifolds.
 \item[(c)] The restriction of the Cayley graph map induces a diffeomorphism of the following open sets  
$\mathscr{U}_{-1}(L):=\{U\in \mathscr{U}(L)~|~1+U \mbox { is Fredholm }\} $   and  \[ \Lag^-{(L)}:=\{L_1\in\Lag~|~(L_1,L) \mbox {Fredholm pair}  \}\]
 \end{itemize}
  \end{theorem}
 \noindent
 \textbf{Proof:} $(a)$ The identity $R_{-}\circ\Gamma_{-}=\id$ boils down to computing the reflection in the graph an the isomorphism $T:\Ker{(J+i)}\rightarrow\Ker{(J-i)}$ in terms of $T$. Notice that the next operator on $\hat{H}$ written relative the decomposition $\hat{H}=\Ker{(J+i)}\oplus\Ker{(J-i)}$,
 \[ \left(\begin{array}{cc}
 0& T^* \\
 T & 0
 \end{array}
 \right)
 \] 
 is an orthogonal reflection whose eigenspace corresponding to the eigenvalue $1$ is $\Gamma_T$. Hence $R_{-}\Gamma_{-}(T)=T$.
 
 In order to see that $\Gamma_{-}\circ R_{-}=\id$ notice that  to show that two lagrangians are equal is enough to prove that one contains the other. Take  $v\in\Ker{(J+i)}$. Then
 \[ \Gamma(R_{-}(L))\ni v+R_L(v)=2P_L(v)\in L
 \]
 
 $(b)$ The only issue one needs to be concerned with is differentiability. Fix a lagrangian $L$ and let  $\Sym(L)\ni S\rightarrow\Gamma_{JS}\in\Lag$ be the Arnold chart centered at $L$. Then using a previous computation, (\ref{projection expression})  we get the  following expression for the reflection $R_{\Gamma_{JS}}$ relative $\hat{H}=L\oplus L^\perp$
  \[ R_{\Gamma_{JS}}=2 P_{\Gamma_{JS}}-1=
\left(\begin{array}{cc}
 (1-S^2)(1+S^2)^{-1} & 2S(1+S^2)^{-1}J^{-1} \\
2JS(1+S^2)^{-1} &  -J(1-S^2)(1+S^2)^{-1}J^{-1}
\end{array}
\right)
\]
  which is differentiable function of $S$. Since $R_{-}=P^{i}\circ R\bigr |_{\Ker{(J+i)}}$, where $P^{i}$ is the projection on the $i$ eigenspace of $J$, we conclude that $R_{-}$ is smooth and therefore its inverse $\Gamma_{-}$ is smooth and so is $\mathscr{C}_L$. 
 
 $(c)$   Notice that it is enough to prove the claim for a single lagrangian $L$, which we will take to be $H^+$. Let $\mathscr{C}=\mathscr{C}_{H^+}$  Since the Fredholm property is an open condition it follows that $\mathscr{U}_{-1}$ is open in $\mathscr{U}(H)$. We will, in fact, see a proof in the next section that $\Lag^-$ is an open set of $\Lag$.
 
  By standard spectral theory, the Fredholm property of $1+U$ implies that \linebreak $\Ker{(1+U)}=\mathscr{C}(U)\cap H^{-}$ is finite dimensional and also $-1\notin \sigma (U|_{\Ker{(1+U)^\perp }})$. 
  
   We can now factor out $\Ker(1+U)$. To that end, let $\breve{H}=\Ker {(1+U)^\perp}$, and ${\breve{H}}^\pm$ be the horizontal/vertical copy of $\breve{H}$ in $\breve{H}\oplus \breve{H}$ and $U'=U|_{\breve{H}}$.  Since $1+U'$ is invertible, the Cayley graph of $U'$ is in the Arnold chart of ${\breve{H}}^+$ and so $\mathscr{C}{(U')}+{\breve{H}}^-$ is closed by proposition $\ref {Ar. chart}$. On the other hand \[ \mathscr{C}(U)+H^{-}=\mathscr{C}{(U')}+{\breve{H}}^-+\Ker{(1+U)}
   \] where $\Ker{(1+U)}\subset H^-$ is finite dimensional and this proves that $\mathscr{C}(U)+H^{-}$ is closed. \\\par
   
 Conversely, let $L\cap H^-$ be finite dimensional and $L+H^-$ be closed. If we let $L'$ and $\breve{H}^-$ be the orthogonal complements of $L\cap H^-$ in $L$ and $H^-$ respectively, these two spaces are lagrangians in $J\breve{H}^-\oplus \breve{H}^-$ whose intersection is empty, and whose sum is closed. Their sum is closed because of the relation
  \[  L'+\breve{H}=(L\cap H^-)^\perp
  \]
   where the orthogonal complement is taken in $L+H^-$. The $\subseteq$ inclusion in this relation is obvious.
 
   For the sake of completeness, we do the $\supseteq$ inclusion.    Every 
  \[ z=z_1+z_2\in L+H^- \]
  can be written as
  \[ \begin{array} {cccccc}
  z=z_1+z_2=&x_1+y_1&+&x_2&+&y_2\\
  &  L\cap H^- & &L' && \breve{H}
   \end{array}
  \] 
  
   So $z=(x_1+y_1)+(x_2+y_2)$ with $x_1+y_1\in L\cap H^-$ and $x_2+y_2\in (L\cap H^-)^\perp$ and therefore if $z \in(L\cap H^-)^\perp$ then $z=x_2+y_2\in L+H^-$.  
   
  We have just proved that $L'$ is in the Arnold chart of $J\breve{H}^-$, hence $L'=\Gamma_S$ where $S:J\breve{H}^-\rightarrow J\breve{H}^-$ is a self-adjoint operator. Finally, we have $L=\mathscr{C}(U)$ where $U$ is the extension by $-1$ on $J(L\cap H^-)$ of the Cayley transform of $S$. It is clear that $1+U$ is Fredholm.\qed\\
  
  Since our main interest is in $\Lag^-$ we will deal with this case separately. \\
  \noindent
  \textbf{Notation:}  To simplify notation we will use $\mathscr{C}$ to denote the Cayley graph map at $H^+$, i.e. $\mathscr{C}:=\mathscr{C}_{H^+}$.
  \begin{corollary} The Cayley graph map $\mathscr{C}:\mathscr {U}(H)\rightarrow \Lag$
  \[ \mathscr{C}(U):=\Ran\{H\ni v\rightarrow ((1+U)v,-i(1-U)v)\in H\oplus H \}
  \]
   induces a diffeomorphism between $\mathscr{U}_{-1}$ and $\Lag^-$. 
 \end{corollary}
 \noindent
 \textbf{Proof:} This is just the case $L=H^+$ in Arnold's theorem. The reason for $-i$ in the second component is  that under the  canonical identifications $H=H^+$ and $H=H^-$, $J\bigr |_{H^+}$ acts as minus the identity.\qed
 
 \begin{remark} \label{sign convention} Our choice of  the reflection map, $R_{-}$  to go from $\Ker{(J+i)}$ to  $\Ker{(J-i)}$ rather then the other way around was not accidental.  The Arnold chart at the vertical space $H^-=0\oplus H$ associates to a self-adjoint operator $A\in \Sym(H^-)$ the switched graph $\tilde{\Gamma}_A:=\{(Av,v)~|~v\in H \}\subset H\oplus H$. In the case when $H=\mathbb{C}$ we want the composition 
\[ \xymatrix {
  S^1 \ar[r]^{\mathscr{C}\quad} & \Lag(1)  \ar[r]^{(\tilde{\Gamma}_{A})^{-1}} & \Sym(\mathbb{C})
   }
\]
defined where it makes sense (i.e. for $\lambda \neq 1$) to be orientation preserving. This is related to the definition of the spectral flow; see \ref{spectral flow}. In our case, the composition is
\[ \lambda\rightarrow i\frac{1+\lambda}{1-\lambda}
\]
 This is indeed orientation preserving as a map from the unit circle minus a point to the real axis. 
\end{remark}
 
For each lagrangian $L$ we introduce the change of basis isomorphism 
\[  \tag{*} \label{*}
\Phi_L:L\oplus L^\perp\rightarrow\Ker{(J+i)}\oplus \Ker{(J-i)},\quad \quad \Phi_L=\left(\begin{array}{cc}
 \phi_{-}& 0 \\
 0 & \phi_{+}\circ J^{-1}
\end{array}
\right)
\]
where $J^{-1}:L^\perp\rightarrow L$ is the inverse of the restriction $J:L\rightarrow L^\perp$. As a map $L\oplus L^\perp\rightarrow L\oplus L^\perp$, $\Phi_L$ has the expression
\[\Phi_L=\frac{1}{\sqrt{2}}\left(\begin{array}{cc}
 1 &  J^{-1}\\
iJ & -i
\end{array}
\right)
\]
 Notice that $\Phi_L$ diagonalizes $J$  relative to the decomposition $\hat{H}=L\oplus L^\perp$, i.e.
\[ \Phi^{-1}_LJ\Phi_L=\left(\begin{array}{cc}
-i &0 \\
0 & i
\end{array}
\right)
\]

\begin{lemma}  Let $U\in\mathscr{U}(L)$ be a unitary map. Then the reflection  in the lagrangian \linebreak $L_U:=\mathscr{C}_L(U)$ has the following expression relative $\hat{H}=L\oplus L^\perp$
\[ R_{L_U}=\Phi_L\left (\begin{array}{cc}
0 &U^*J^{-1} \\
JU & 0
\end{array}
\right)\Phi^{-1}_L
\]
\end{lemma}

\noindent
\textbf{Proof:} Let $T_U:\Ker{(J+i)}\rightarrow\Ker{(J-i)}$ be the isometry that corresponds to $\mathscr{C}_L(U)$, in other words $T_U=(\phi_{+})U(\phi_{-})^{-1}$. Then
\[ R_{L_U}=\left(\begin{array}{cc}
0 &T_U^* \\
T_U & 0
\end{array}
\right)=\left(\begin{array}{cc}
\phi_{-} & 0 \\
 0& \phi _i
\end{array}
\right)\left(\begin{array}{cc}
0 &U^* \\
U & 0
\end{array}
\right)\left(\begin{array}{cc}
\phi_{-}^{-1} & 0 \\
0 & \phi_{+}^{-1}
\end{array}
\right)
\]
The claim follows from the expression for $\Phi_L$.\qed\\

Using the previous lemma we get a different characterization of $\Lag^-$ which we record:

\begin{corollary} The space of vertical, Fredholm lagrangians has the following characterization $~\Lag^{-}=\{L\in \Lag~|~R_L+R_{H^+} \mbox{ is Fredholm}\}$.
\end{corollary}
\noindent
\textbf{Proof:}  If $U$ is the operator $L$ is coming from via the Cayley graph map then
\[R_L+R_{H^+}=
\Phi_{H^+}
\left(\begin{array}{cc}
0&(1+U)^*J^{-1}\\
J(1+U)&0
\end{array}
\right)\Phi_{H^+}^{-1}
\]
which is Fredholm if and only if $1+U$ is Fredholm. \qed\\

\begin{corollary}\label{Cayley and Arnold charts} The Cayley graph map, $\mathscr{C}$ takes the Arnold chart around $U_0$ bijectively onto the Arnold chart around $L_0:=\mathscr{C}(U_0)$.
\end{corollary}
\noindent
\textbf{Proof:} Let $T:=U+U_0$, $L_U:=\mathscr{C}(U)$. Then 
\[ R_{L_U}+R_{L_0 }=\Phi_{H^+}\left(\begin{array}{cc}
0&-T^*\\
-T&0
\end{array}
\right)\Phi_{H^+}^{-1}
\]  is invertible if and only if $T$ is invertible and hence by proposition $\ref {Ar. chart}$ we get $\mathscr{C}_{H^+}(U) \in\mathcal{A}_{L_0}$ if and only if $T$ is invertible.\qed\\

\begin{corollary} \label{reflection diff} Let $M$ be a differentiable manifold and $F:M\rightarrow \Lag$  be a map. Then $F$ is differentiable if and only if $R\circ F$ is differentiable where $R:\Lag\rightarrow\mathscr{B}(\hat{H})$ associates to every lagrangian $L$ the reflection in it. 
\end{corollary}
\noindent 
\textbf{Proof:} Clearly $F$ is differentiable if and only if $F_1:=(\mathscr{C})^{-1}\circ F :M\rightarrow \mathscr{U}(H^+)$ is differentiable. Now $R\circ \mathscr{C}:\mathscr{U}(H)\rightarrow\mathscr{B}(\hat{H})$ has the following expression
\[ R\circ \mathscr{C}(U)=\Phi_{H^+}
\left(\begin{array}{cc}
0&- U^*\\
-U&0
\end{array}
\right)\Phi_{H^+}^{-1}
\]
So $F_1$ is differentiable if and only if $R\circ \mathscr{C}\circ F_1$ is differentiable.\qed\\

A closer look at the Cayley graph map suggests a useful reformulation of Arnold's theorem. Notice that for every unitary operator $U\in\mathscr{U}(L)$,  with $L_U:=\mathscr{C}_L(U)$ the map $\tilde{U}:L\rightarrow L_U$
\begin{center}
$\tilde{U} (v)= 1/2((1+U)v+iJ(1-U)v))$
\end{center}
is a \textit {canonical} Hilbert space isomorphism that carries $L$ into $L_U$. We can build out of $\tilde{U}$ an automorphism  $U^\sharp$ of $\hat{H}$ by taking the direct sum  of $\tilde{U}$ with $J\tilde{U}J^{-1}:L^\perp\rightarrow L_U^\perp$.  Written in the decomposition $\hat{H}=L\oplus L^\perp$ this automorphism has the expression:
\[ U^\sharp= \frac{1}{2} \left(\begin{array}{cc}
1+U & -i(1-U)J^{-1}\\
iJ(1-U)& J(1+U)J^{-1} 
\end{array}
\right) =\Phi_L\left(\begin{array}{cc}
1 & 0\\
0 & JUJ^{-1}
\end{array}
\right)\Phi_L^{-1}
\]

We are now ready to give the reformulation of Arnold's theorem.
\begin{theorem} \begin{itemize}
\item [(a)]
 For a fixed lagrangian $L$, the following map is a Banach-Lie group embedding 
\[ \mathscr{O}_{L}: \mathscr{U}(L)\mapsto \mathscr{U}(\hat{H}), \quad \quad \label{unitaryinhat}
\mathscr{O}_{L}(U)=\Phi_L\left(\begin{array}{cc}
1 & 0\\
0& JUJ^{-1} 
\end{array}
\right)\Phi_L^{-1}
\]
where the decomposition is relative $\hat{H}=L\oplus L^\perp$. The orbit of the lagrangian $L$ under the action of this subgroup, i.e., $\{\mathscr{O}_{L}(U)L~|~U\in \mathscr{U}(L)\}$ is the entire space $\Lag$ and the stabilizer is trivial. The bijection
\[ \mathscr{U}{(L)}\rightarrow \Lag, \quad \quad U\rightarrow \mathscr{O}_{L}(U)L
\]
is a diffeomorphism of real Banach manifolds.

\item [(b)] The subgroup of $\mathscr{U}(\hat{H})$ determined by $\mathscr{O}_{L}$ does not depend on $L$, but only on $J$. Relative to the decomposition $\hat{H}=\Ker{(J+i)}\oplus\Ker{(J-i)}$, the map $C_L$ has the expression:
\[ \mathscr{O}_{L}(U)=\left(\begin{array}{cc}
1 & 0\\
0& \phi_{+}U \phi_{+}^{-1}
\end{array}
\right)
\] 

\item [(c)] The bundle $\mathscr{U}_{\tau }\subset\Lag\times\mathscr{B}(\hat{H})$ over  $\Lag$ whose fiber over a lagrangian $L$  consists of unitary operators $U\in\mathscr{U}(L)$ is canonically  trivializable and the map:
\[\mathscr{O}: \mathscr{U}_{\tau }\rightarrow \mathscr{U}(\hat{H}), \quad \quad \mathscr{O}(L,U):=\mathscr{O}_L(U)
\]
is differentiable.

\end{itemize}
\end{theorem}
\noindent
\textbf{Proof:} $(a)$ Note that  $\mathscr{O}_{L}(U)L=\mathscr{C}_L(U)$ and so the job is done by Theorem \ref{Arnold}.

$(b)$ Self-explanatory.

$(c)$  Let us notice that the tautological bundle $\tau\subset \Lag\times \hat{H}$, $\tau:=\{(L,v)~|~v\in L\}$ over $\Lag$ is naturally trivializable. A natural trivialization is given as follows.  For every lagrangian $L$, let $U_L:=(\mathscr{C})^{-1}(L)$ be the unitary operator on $H^+$ from which $L$ is coming. Then the  following map is a trivialization of the tautological bundle
\[ \alpha: \tau\rightarrow \Lag\times H^+, \quad \quad \alpha(L,v)=(L, \left(\mathscr{O}_{H^+}(U_L)\right)^{-1}(v))
\]
since both $(\mathscr{C})^{-1}$ and $\mathscr{O}_{H^+}$ are differentiable. It is straightforward now that $\mathscr{U}_{\tau}$ is naturally trivializable. 

In order to show that $\mathscr{O}$ is differentiable it is enough to show that the map
\[\Phi: \Lag\rightarrow \mathscr{U}(\hat{H}),\quad \quad L\rightarrow \Phi_L
\]
is differentiable. Since $\mathscr{O}_{H^+}$ and  $ \mathscr{C}_{H^+}$ are differentiable, the following identity proves this claim.
\[ \Phi ( \mathscr{C}(U))=\mathscr{O}_{H^+}(U)\Phi_{H^+}=\frac{1}{\sqrt{2}}\left(\begin{array}{cc}
1 & -U \\
-i & -iU
\end{array}\right), \quad \quad \forall U\in \mathscr{U}(H^+)
\]
The decomposition is relative $\hat{H}=H\oplus H$.
To see why the identity is true, let $L_U:=\mathscr{C}_{H^+}(U)$ and $v\in H^+$.  Then 
\[ x:=\mathscr{O}_{H^+}(U)(v)=1/2((1+U)v,-i(1-U)v)\in L_U 
\] and 
\[ \Phi_{L_U}(x)=1/\sqrt{2}(x+iJx)=1/\sqrt{2}(v,-iv)\]
Let $w\in H^-$. Then
 \[ y:=\mathscr{O}_{H^+}(U)(w)=1/2(i(1-U)w,(1+U)w)\in L_U^\perp\] and 
\[ \Phi_{L_U}(y)=1/\sqrt{2}(J^{-1}y-iy)=1/\sqrt{2}(-Uw,-iUw).
\]
\qed\\

\noindent
\textbf{Notation:} For a lagrangian $L$ and a unitary map $U\in\mathscr{U}(L)$ we will use the notation $\tilde{U}$ for the Hilbert space isomorphism:
\[ \tilde{U}:L\rightarrow \mathscr{C}_L(U), \quad \quad \tilde{U}(v)=\mathscr{O}_{L}(U)(v)=1/2((1+U)v+Ji(1-U)v)
\]\qed

We see that Arnold's theorem gives more than just a diffeomorphism between two real Banach manifolds. It shows that  given two lagrangians $L_1$ and $L_2$, there exists \textit{canonical} unitary operators 
 \[
 U(L_1,L_2)\in \mathscr{U}(L_1), \quad \quad \quad \quad U(L_1,L_2):=(\mathscr{C}_{L_1})^{-1}(L_2)\]
  \[ U(L_2,L_1)\in \mathscr{U}(L_2),\quad \quad \quad \quad U(L_2,L_1):=(\mathscr{C}_{L_2})^{-1}(L_1) \]   that induce isomorphisms between the two lagrangians
\[ \tilde{U}(L_1,L_2):L_1\rightarrow L_2 \quad \quad \mbox{ and } \quad \quad \tilde{U}(L_2,L_1):L_2\rightarrow L_1
\]
Notice that for every lagrangian $L$ we have:
\[ \tilde{U}(L,L)=\id
\]
In fact we have the following result
\begin{lemma} For every two lagrangians $L_1$  and $L_2$ and  for every unitary map $X\in\mathscr{U}(L_1)$ the following identities hold: 
\begin{itemize}
\item[(a)] $U(L_1,L_2)=(\phi_{+}{(L_1)})^{-1}\circ\phi_{+}{(L_2)}\circ(\phi_{-}{(L_2)})^{-1}\circ\phi_{-}{(L_1)}$;
\item[(b)] $\tilde{U}(L_1,L_2)=\phi_{-}^{-1}(L_2)\circ \phi_{-}(L_1) $;
\item[(c)] $\mathscr{O}(L_1, U(L_1,L_2)XU(L_1,L_2)^{-1})=\mathscr{O}(L_2, \tilde{U}(L_1,L_2)X\tilde{U}(L_1,L_2)^{-1})$.
\end{itemize}
\end{lemma}
\noindent
\textbf{Proof:} $(a)$ This identity follows from the commutativity of the diagram
\[
\xymatrix{ L_1 \ar^{U(L_1,L_2)}[r]  \ar_{\phi_{-}{(L_1)}}[d]&  L_1 \ar^{\phi_{+}{(L_1)}}[d] \\
\Ker{(J+i)}  &  \Ker{(J-i)} \\
L_2 \ar^{\id}[r]  \ar^{\phi_{-}{(L_2)}}[u]&  L_2 \ar_{\phi_{+}{(L_2)}}[u] 
}\]
To see why this diagram is commutative think that both  $\id_{L_2}$ and $U(L_1,L_2)$ induce an isomorphism $T\in \Isom{(\Ker{(J+i)},\Ker{(J-i)})}$ whose graph is exactly $L_2$. \\

\noindent $(b)$ We have
\[\tilde{U}(L_1,L_2)=\left(\begin{array}{cc}
1 & 0\\
0& \phi_{+} (L_1)U(L_1,L_2) (\phi_{+}(L_1)) ^{-1}
\end{array}
\right) \Bigr |_{L_1}=\left(\begin{array}{cc}
1 & 0\\
0& \phi_{+} (L_2) (\phi_{-}(L_2)) ^{-1}
\end{array}
\right) \Bigr |_{L_1}
\]
where the decomposition is relative $\hat{H}=\Ker{(J+i)}\oplus\Ker{(J-i)}$. On the other hand for $v\in L_1$ we have:
\[ v=1/\sqrt{2}(\phi_{-}(L_1)v,\phi_{+}(L_1)v)\in\Ker{(J+i)}\oplus\Ker{(J-i)}
\]
and so 
\[\tilde{U}(L_1,L_2)v=1/\sqrt{2}(\phi_{-}(L_1)v,\phi_{+} (L_2) (\phi_{-}(L_2) ^{-1}\phi_{+}(L_1)v)=1/\sqrt{2}(\phi_{-}(L_2)w,\phi_{+}(L_2)w)
\]
 for some $w\in L_2$. The identity is now  obvious. \\
 
 \noindent $(c)$  Let $U_{12}:=U(L_1,L_2)$.  We consider $(\tilde{U}_{12} v,J\tilde{U}_{12}J^{-1}w)\in L_2\oplus L_2^\perp$ where  $(v,w)\in L_1\oplus L_1^\perp$. Notice that for all pairs $(v,w)$ we have

\[ \mathscr{O}(L_2, \tilde{U}_{12}X\tilde{U}_{12}^{-1}) (\tilde{U}_{12} v,J\tilde{U}_{12}J^{-1}w)=(\tilde{U}_{12} a,J\tilde{U}_{12}J^{-1}b)
\]
where
\[ \left(\begin{array}{c}
a\\
b
\end{array}
\right)= \frac{1}{2} \left(\begin{array}{cc}
1+X & -i(1-X)J^{-1}\\
iJ(1-X)& J(1+X)J^{-1} 
\end{array}
\right) \left(\begin{array}{c}
v\\
w
\end{array}
\right)=\mathscr{O}(L_1,X)(v,w)
\]
The same relation can be written  as
\[ \mathscr{O}(L_2, \tilde{U}_{12}X\tilde{U}_{12}^{-1}) \mathscr{O}(L_1, U_{12})=\mathscr{O}(L_1, U_{12})\mathscr{O}(L_1,X)
\]
from which the identity follows.\qed

\begin{corollary}
 Let $L_i$, $i\in\{1,2,3\}$ be three lagrangians and let  $\tilde{U}(L_1,L_2)$  be the canonical isomorphisms between them.  Then
\[ \tilde{U}(L_2,L_3)\circ \tilde{U}(L_1,L_2)=\tilde{U}(L_1,L_3)
\]
\end{corollary}

\begin{prop}\label {tangent space}
The tangent space of $\Lag$ is isomorphic  with  the tautological bundle \linebreak $\mathscr{T}: =\{(L,S)\in \Lag \times \Sym{(\hat{H})}~|~ S\in \Sym(L)\}$ and they are both naturally trivializable.
\end{prop}
\noindent
\textbf{Proof:} The isomorphism between the tangent space of $\Lag $ and $\mathscr{T}$ is given by:
\begin{center}
$(L,[\alpha])\rightarrow(L, \alpha'(0))$
\end{center}
where $[\alpha]$ is an equivalence class of curves and $\alpha'(0)$ is the derivative of a representative of $\alpha$ in the Arnold chart of $L$. 

For every lagrangian $L$ let $U_L:=\mathscr{C}^{-1}_{H^+}(L)$. Consider now the  following automorphism of the trivial bundle 
\[  \begin{array} {c}
\Lag \times \Sym{(\hat{H})} \mapsto \Lag \times \Sym{(\hat{H})}, \quad  \quad
(L,S)\mapsto (L,\mathscr{O}(U_L)^{-1}S\mathscr{O}(U_{L}))
\end{array}
\]
It restricts to an isomorphism
\[ 
\mathscr{T}\rightarrow \Lag \times \Sym(H^+), \quad \quad
(L,S)\rightarrow (L,\tilde{U}_L^{-1}S\tilde{U}_L)
\]
\qed\\

We saw in Corollary \ref{Cayley and Arnold charts} that the Arnold chart around a unitary operator $U_0$ is taken by the Cayley graph map, $\mathscr{C}_{H^+}$,  to the Arnold chart around the corresponding lagrangian $L_0:=\mathscr{C}(U_0)$. The Arnold chart of $U_0$ consists of those operators $U$ such that $U_0+U$ is invertible. These operators live in the image of the composition
\[ \xymatrix {  \Sym(H) \ar[r] \ar @/_1pc/ [rr]_{\Cay_0} & \mathscr{U}(H) \ar^{U_0\cdot}[r] & \mathscr{U}(H),}
 \quad \quad A\mapsto \frac{i-A}{i+A}=:U \mapsto U_0 U
\]
where the first map is the Cayley transform at identity. When $H=\mathbb{C}^n$ this Cayley transform preserves the natural orientation on $U(n)$. We denote this composition by $\Cay_0$.  On the other hand,  $\Sym(L_0)$ is naturally isomorphic with $\Sym (H)$ by conjugation with the Hilbert space isomorphism $\tilde{U}_0: H\rightarrow L_0 $ given by 
\[ v\rightarrow 1/2((1+U_0)v, -i(1-U_0)v)
\]
Let $\tilde{\mathscr{A}}_0$ be the application that maps $\Sym(H)$ bijectively to the Arnold chart of $L_0$.
\[ \xymatrix { \Sym(H) \ar[r] \ar @/_1pc/ [rr]_{\tilde{\mathscr{A}}_0} &  \Sym(L_0)  \ar[r] & \Lag,} \quad \quad A\mapsto \tilde{U}_0A\tilde{U}_0^{-1}=:X \mapsto \Gamma_{JX}
\]
\begin{prop}\label{useful for orientation} The composition $\tilde{\mathscr{A}}^{-1}_0\circ \mathscr{O} \circ \Cay_0:\Sym(H)\rightarrow\Sym(H)$ is the identity.
\end{prop}
\noindent
\textbf{Proof:} We know that the Arnold chart at $U_0$ maps to the Arnold chart at $L_0$.  Every unitary operator $U\in\mathscr{U}(H)$ such that $1+U$ is invertible maps via left multiplication by $U_0$ to an operator $U'$ such that $U_0+U'$ is invertible.  We fix such an operator $U$ and let $\tilde{A}\in\Sym(L_0)$ be defined by
\[ \tilde{A}:= \tilde{U}_0\frac{i(1-U)}{1+U}\tilde{U}_0^{-1}
\] 

We show that  $\tilde{A}$ corresponds via the map $\tilde{A}\mapsto \Gamma_{J\tilde{A}}$ to the lagrangian determined by $U_0U$.  Indeed, let $ x\in H$ and $v\in L_0$ be two vectors related by $v=\tilde{U}_0(1+U)x$. Then
\[ v+J\tilde{A}v=\left(\tilde{U}_0(1+U)x+J\tilde{A}\tilde{U}_0(1+U)x\right)=\left(\tilde{U}_0(1+U)x +iJ\tilde{U}_0 (1-U)x\right)
\]
\[ =1/2\bigr((1+U_0)(1+U)x,-i(1-U_0)(1+U)x\bigr)  +\bigr(1-U_0)(1-U)x, -i(1+U_0)(1-U)x \bigr)
\]
\[ =\bigr( (1+U_0U)x, -i(1-U_0U)x\bigr)\in\mathscr{O}(U_0U)
\]
Since the correspondence $v\leftrightarrow x$ is bijective the claim follows.

Plugging in $U=(i-A)(i+A)^{-1}$ with $A\in \Sym(H)$ finishes the proof.
\qed

\bigskip

\subsection{More examples of Lagrangians}\label{moreexamples}

In this section we will discuss what happens with the Palais groups $\mathscr{U}_{\mathscr{I}}$ through the Cayley graph bijection.  Let
 \[
 \Lag^-_{\mathscr {I}} :=\mathscr{C}(\mathscr{U}_\mathscr{I})\]
 be the space of lagrangians that correspond to the Palais groups.   The computation of the reflection $R_{L_U}$ in Arnold's theorem proves again useful for an intrinsic characterization of $\Lag^-_{\mathscr {I}}$. Here $L_U:=\mathscr{C}(U)$.

\begin {center}
$R_{L_U}+R_{H^-}=\Phi_{H^+}\left(\begin{array}{cc}
0&1-U^*\\
1-U&0
\end{array}
\right)\Phi_{H^+}^{-1}$ 
\end{center}
Hence, if we let $\mathscr{I}_\Phi$ to denote the  Banach subspace of $\mathscr{B}(\hat{H})$
\[
\mathscr{I}_\Phi:=\{X\in \mathscr{B}(\hat{H})~|~ X=\Phi_{H^+}\left(\begin{array}{cc}
0&T^*\\
T&0
\end{array}
\right)\Phi_{H^+}^{-1}, ~T\in \mathscr{I}\}=\] 
\[=\left\{X \in \mathscr{B}(\hat{H})~|~X=\left(\begin{array}{cc}
T+T^*& i(T-T^*)\\
i(T-T^*)& -(T+T^*) 
\end{array}
\right), ~T\in \mathscr{I}\right\}
\]
 we have the following description/definition for $\Lag^-_{\mathscr {I}}$.
 
\begin{lemma} \label{lagrangian and ideals}  The space of  $\mathscr{I}$-commensurable lagrangians with $H^+$, $\Lag^-_{\mathscr {I}}$ (or just commensurable when $\mathscr{I}=\mathscr{K}$) has the following description as a subset of $\Lag^-$
\[
\Lag^-_{\mathscr{I}}=\{L\in \Lag^- ~|~ R_L+R_{H^-}\in \mathscr{I}_\Phi\}=
\{L\in \Lag^- ~|~ P_L-P_{H^+}\in \mathscr{I}_\Phi\}.
\]
\end{lemma}

The space $\mathscr{I}_\Phi$ is not an ideal of $\mathscr{B}(H\oplus H)$, (it is not even an algebra) but the obvious topology that it inherits from $\mathscr{I}$, topology which we will denote by superscript $\mathbf{s}$ (from $\mathbf{s}$trong) makes the identity map continuous in one direction:
\begin{center}
$\id: \mathscr{I}_\Phi^{\mathbf{s} }\rightarrow\mathscr{I}_\Phi^{\|\cdot\| }$
\end{center}
With this topology on $\Lag^-_{\mathscr{I}}$ the Cayley graph map becomes a homeomorphism when restricted to $\mathscr{U}_{\mathscr{I}}$.\\
\noindent
\begin{remark} It is tempting to have a description of $\Lag^-_{\mathscr {I}}$ in terms of the projections in the spirit of the Lemma $\ref {char proj}$. Such a description suggests itself and we  would like to say that if the lagrangian $L$ is in $\Lag^-_{\mathscr {I}}$ then $P^-|_L \in \mathscr{I}$. Nevertheless such a statement must be taken with a grain of salt.
This is because the set of projections $P^-|_L $ is not a subspace of $\mathscr{B}(H)$ and in general it does not have an algebraic structure.

On the other hand it is true that for  $L\in\Lag^{-}_{\mathscr{I}}$, one has $P^-|_L\circ \tilde{U}_L\in\mathscr{I}$, where $\tilde{U}_L$ is the unitary isomorphism $H\rightarrow L$ provided by the Cayley map and conversely if $P^-|_L\circ U_L\in\mathscr{I}$ then $L\in\Lag^{-}_{\mathscr{I}}$. Therefore, if we take $P^-|_L \in \mathscr{I}$ to mean \textit{only} a certain boundedness condition on the singular values of $P^-|_L$, the same one that describes $\mathscr{I}$ (\cite {S} ), e.g., trace class or Hilbert-Schmidt condition, then the previous description makes perfect sense:
\[\Lag^{-}_{\mathscr{I}}=\{L\in \Lag(H)~|~P^+|_L \mbox{  is Fredholm of index }0 \mbox{ and }P^-|_L \in \mathscr{I}\}\]\qed
\end{remark}

  We  will put charts on $\Lag^-_{\mathscr{I}}$ and turn it into a manifold modeled on the Banach space $\Sym_\mathscr{I}$. With this manifold structure the Cayley graph map becomes a diffeomorphism 
  \[ \xymatrix{  \mathscr{U}_\mathscr{I}\ar^{\mathscr {C}}[r] &\Lag^-_{\mathscr{I}}}\]

 We know for example that $\mathcal{A}_{H^+}\subset \Lag^-$ and it is easy to see that in order to have  $ \Gamma_{T}\in \Lag^-_{\mathscr{I}}$ for a self-adjoint operator $T:H\rightarrow H$ one needs to have $T\in\mathscr{I}$. Hence 
$\mathcal{A}_{H^+}\cap \Lag^-_{\mathscr{I} }\simeq \Sym_\mathscr{I}(H)$. The first set is open in the norm topology, so it is open in the stronger $\mathbf{s}$ topology. This way we obtained a first chart on $\Lag^-_{\mathscr{I} }$.  In order to put other charts on this space notice that we can unambiguously talk about the ideal $\mathscr{I}\subset\mathscr{B}(H)$ as an ideal of $\mathscr{B}(L)$ for any separable, Hilbert space $L$. Indeed since every ideal is conjugation invariant one can just take a unitary isomorphism between $H$ and $L$ and "transfer" $\mathscr{I}$, via conjugation,  to a subset of $\mathscr{B}(L)$.\footnote {In other words we do not have to define what we understand by a compact or Hilbert-Schmidt operator on each separable, Hilbert space. } \par
  Let us fix a finite codimensional space $V\subset H^+$. Then 
  
\begin{lemma} The lagrangian $V\oplus JV^\perp$ belongs to all subsets $\Lag_{\mathscr{I}}^-$. The intersection of the Arnold chart in $\Lag^-$ around $V\oplus JV^\perp$ with $\Lag^-_{\mathscr{I} }$ is an Arnold chart of $\Lag^-_{\mathscr{I} }$ around $V\oplus JV^\perp$. In other words the graph map induces a bijection
\[\Sym_\mathscr{I}(V\oplus JV^\perp) \simeq \mathcal{A}_{V\oplus JV^\perp}\cap  \Lag^-_{\mathscr{I} }. 
\]
\end{lemma}
\noindent
\textbf{Proof:} The first claim is true because $P^+|_{V\oplus JV^\perp}$ is clearly Fredholm and $P^-|_L$ is a finite rank operator.

One easy observation is that the unitary map $U\in \mathscr{U}(H^+)$ in Arnold's theorem which satisfies $\mathscr{O}(U)H^+=V\oplus JV^\perp$ is:
\begin{center}
$U=\left\{\begin{array}{cc}
I    & \mbox { on } V \\
-I  & ~\mbox{ on } V^\perp
\end{array}
\right.
$
\end{center}
Let $T\in \Sym(V\oplus JV^\perp)$. It is not hard to see that $\Gamma_{JT}=\mathscr{O}(U)\Gamma_{J\tilde{T} }$ where $\tilde{T}=U^{-1}TU\in\Sym_\mathscr{I}(H^+)$. (here $U$ is the unitary isomorphism $H^+\rightarrow V\oplus JV^\perp$ induced by the Cayley map of $U\in \mathscr{U}(H^+)$ defined before). Indeed we have :
\[ \mathscr{O}(U)\Gamma_{J\tilde{T} }=
\left\{
\left(\begin{array}{cc}
U & 0 \\
0 &JUJ^{-1}
\end{array}
\right) 
\left(\begin{array}{c}
v\\
J\tilde{T}v
\end{array}
\right) 
~|~v\in H^+\right\}
\]
and the latter set is just $\Gamma_{JT}$.

If we let $U_{\tilde{T}}:=\frac{1+i\tilde{T} }{1-i\tilde{T}}$ be the Cayley transform of $\tilde{T}$ then $\Gamma_{JT}=\mathscr{O}(U)\mathscr{O}(U_{\tilde{T}})H^+=\mathscr{O}(UU_{\tilde{T}})(H^+)$. Hence $UU_{\tilde{T}}$ is the unitary transform in Arnold's theorem which takes $H^+$ to $\Gamma_{JT}$.  It is easy to see that $UU_{\tilde{T}}\in 1+\mathscr{I}$ if and only if $\tilde{T}\in \mathscr{I}$ iff $T\in\mathscr{I}$.
\qed\\

 \begin{remark} The fact that the Arnold charts $\mathcal {A}^{\mathscr{I}}_{V\oplus JV^\perp}:=\mathcal{A}_{V\oplus JV^\perp}\cap  \Lag^-_{\mathscr{I} }$  with $V$ cofinite dimensional in $H^+$, cover $\Lag^-_{\mathscr{I}}$  is just an observation in the next section.
 \end{remark}

 Simple examples of vertical, Fredholm lagrangians come from the  graphs of  bounded, self-adjoint operators $S:H\rightarrow H$.   More interesting examples arise when one looks at switched graphs of operators. In this case in order for   the \textit{switched} graph to be Fredholm pair with the vertical space, $T$ itself has to be Fredholm. Most importantly, $T$ need not even be bounded.

 The usual framework is as follows.  Let  $T$ be a closed, densely defined, self-adjoint operator $T:D(T)\subset H\rightarrow H$ with a compact resolvent. In particular $T$ has only discrete spectrum and it is Fredholm.  Not only is its switched graph a vertical, Fredholm lagrangian, it is in fact commensurable with $H^+$, which is another way of saying that it sits in  $\Lag_{K}^-=\mathscr{C}(U_{\mathscr{K}}(H))$. To see this, let $K_1:(~ \Ker{T})^\perp\rightarrow(~ \Ker{T})^\perp$ be the  self-adjoint, compact operator which is the inverse of $T|_{(\Ker T)^\perp }$ (by $T|_{(\Ker T)^\perp}$ we mean the operator defined on the projection of $D(T)$ along $\Ker{T}$). If 
\begin{center} 
$K:=\left(\begin{array}{cc}
\frac{2K_1}{i-K_1}&0\\
0&-2I
 \end{array} \right)$
\end{center}
 is an extension of $K_1$ by a finite-rank operator to a compact operator on
\begin{center}
$H=(~ \Ker {T})^\perp\oplus~ \Ker{T}$
\end{center} 
  one can check that 
\begin{center}
$U:=1+K=\left(\begin{array}{cc}
\frac{i+K_1}{i-K_1}&0\\
0&-I\end{array}
\right)$
\end{center} is unitary and that $\mathscr{C}(U)=\tilde{\Gamma}_T$.  

If  the compact parametrix $K_1\oplus P|_{\Ker{T}}$ is in the ideal $\mathscr{I}$ which contains the finite rank operators then $\tilde{\Gamma}_T\in \mathscr{C}(U_\mathscr{I}(H))$.

Let $T$ be a self-adjoint, Fredholm operator and let $U=\mathscr{C}^{-1}(\tilde\Gamma_{T})$ be the unitary operator it corresponds to via the Arnold isomorphism. Then the operator $X=1-U$ is bounded and induces a bijection $X:H\rightarrow D(T)$. If we let 
\[\langle v,w\rangle_{D(T)}=\langle X^{-1}v,X^{-1}w\rangle_{H}, \quad \forall v,w\in D(T)
\]
be an inner product on $D(T)$ induced by $X$ then with this inner product two things are true, the inclusion $D(T)\rightarrow T$ is continuous and $X:H\rightarrow D(T)$ is a Hilbert space isomorphism. The following result is now straightforward:

 \begin{lemma}\label{new norm} For every closed, self-adjoint operator $T:D(T)\subset H\rightarrow H$ there exists a  bounded operator $X\in \mathscr{B}(H)$ which is  a bijection onto $D(T)$ and  which induces an inner product on $D(T)$ such that $T: D(T)\rightarrow H$ becomes a bounded operator and 
 \[ T=iX^{-1}(2-X)
 \] 
 \end{lemma}
 
 The new norm on $D(T)$ is nothing else but the graph norm as the next result shows.
 \begin{definition}\label{graph norm}  For every closed, densely defined operator $T:D(T)\subset H\rightarrow H $  \textit{\textbf{the graph norm}} \index{graph norm} on $D(T)$ is 
    \[ \|v\|_g=\|Tv\|+\|v\|, \quad \forall v\in D(T)\]
\end{definition}

\begin{lemma}\label{elliptic estimates} Let $H_0$ be a dense subset of $H$ and let $T:H_0\rightarrow H$ be a self-adjoint operator. Consider the norm $\langle\cdot,\cdot\rangle_0$ on $H_0$ induced by the operator $X$ from the previous lemma. Then this is equivalent with the graph norm.
\end{lemma}
\noindent
\textbf{Proof:} Let $v\in D(T)$ and let $w=X^{-1}(v)\in H$. Notice that 
\[ \|v\|_0\leq 1/2\|v\|_g
\]
Indeed the inequality is equivalent with
  \[  \|w\|\leq1/2(\|TXw\|+\|Xw\|)   \] 
which becomes by (\ref{new norm})
  \[ \|w\|\leq 1/2(\|(2-X)w\|+\|Xw\|)
  \]
  This is just the triangle inequality. 
  
  Conversely one has 
 \[ \|v\|_g=\|(2-X)w\|+\|Xw\|\leq (\|X\|+\|2-X\|)\|w\|=M\|v\|_0
 \]\qed
 
We have the following result
\begin{lemma} If $U=\mathscr{C}^{-1}(\tilde\Gamma_{T})$ is the unitary operator that corresponds to the switched graph of $T$ and $R_{T}(-i)=(T+i)^{-1}$ is the resolvent of $T$ at $-i$ then
\[ U=1-2iR_T(-i)
\]
\end{lemma}
\noindent
\textbf{Proof:} Let $X:=1-U$. By Lemma \ref{new norm} we have
\[ X=2i(T+i)^{-1}=1-U
\]\qed
\begin{corollary} The switched graph of a self-adjoint, Fredholm operator $T$ is in $\Lag_{\mathscr{I}}$ if and only if the resolvent $R_{T}(\lambda)\in \mathscr{I}$ for some $\lambda\notin\sigma(T)$.
\end{corollary}
\noindent
\textbf{Proof:} By the previous lemma $1-U\in\mathscr{I}$ if and only if $R_T(-i)\in\mathscr{I}$.\qed

\bigskip

\subsection{Symplectic reduction} \label{Symp red}

Our big goal for this section is to prove that $\Lag^{-}$, together with its little brothers, $\Lag_{\mathscr{I}}^{-}$ are all classifying spaces for odd $K$-theory, i.e. they all have the same homotopy type as $U(\infty)$. To achieve this we use the technique called (linear) symplectic reduction which we now describe. 

\begin{definition} An \textit{\textbf{isotropic}}\index{isotropic subspace} subspace of the complex Hilbert space $\hat{H}=H\oplus H$, endowed with a complex structure as in the previous section is a \textit{closed} subspace $W\subseteq \hat{H}$ such that $JW\subseteq W^\perp$.

The space $W^\omega:=(JW)^\perp$ is called the \textit{\textbf{annihilator}}\index{annihilator} of $W$. 

For an isotropic space $W$, the orthogonal complement of $W$ in $W^\omega$, denoted $H_W$ \index{$H_W$} is called the \textit{\textbf{(symplectically) reduced}} space of $\hat{H}$. 
\end{definition} 
    One trivial observation is that $H_W$ is the orthogonal complement of $W\oplus JW$ in $\hat{H}$. Notice that $H_W$ is $J$-invariant since its orthogonal complement is.
 \begin{definition} 
 The isotropic space $W$ is called \textit{\textbf{cofinite}}\index{isotropic ! cofinite} if  $H_W$ is finite dimensional and \linebreak $\dim{\Ker{(i\pm J|_{H_W})}}=\frac{1}{2}$dim$~H_W$.
 \end{definition} 
 \noindent
 \begin{remark} The reason for considering the signature zero condition for $iJ|_{H_W}$ in the previous definition is because we want the Lagrangian Grassmannian $\Lag(H_W)$ to be nontrivial whenever $H_W$ is not trivial.   We want it nontrivial  because this guarantees that every maximal, cofinite, isotropic space is actually a lagrangian. Indeed with our definition in place, if a maximal, cofinite, isotropic space were not a lagrangian, it would mean that $H_W$ is non-zero and hence $\Lag(H_W)$ is nontrivial and so a choice of a lagrangian in $H_W$ added to the initial isotropic space would deliver a bigger isotropic space contradicting the maximality.
 \end{remark}
 \noindent
 \begin{example} Every finite codimensional subspace $W\subset H^\pm$  is in fact cofinite, isotropic and the same is true about every finite codimensional subspace of any lagrangian, not only $H^\pm$ . 
 \end{example}
 \begin{definition}
 For a fixed, cofinite, isotropic space $W$ we say that the lagrangian $L$ is $\textit{\textbf {clean}}$\index{clean lagrangian} with $W$ if it belongs to the set:\index{$Lag^W$}
\begin{center}
$\Lag^W:=\{L\in \Lag, ~L\cap W=\{0\}, ~L+W$ closed$\}$
\end{center}
\end{definition}
 If $W=L_0$ is a lagrangian itself, we have $\Lag^W =\mathcal{A}_{L_0^\perp }$.
 
  The next proposition/definition is fundamental.
\begin{prop}\index{symplectic reduction}\label{symp red} For every cofinite, isotropic space $W$  the following are true:
\begin{itemize}
\item[a)] The  map 
\[\mathscr{R}:\Lag^W\mapsto \Lag(H_W), \quad \quad~~~~~ L\mapsto~\Range{ P_{H_W}|_{L\cap W^\omega}}\]
is well-defined. Here $P_{H_W}|_{L\cap W^\omega}$ is the orthogonal projection onto $H_W$ restricted to $L\cap W^\omega$. The map $\mathscr{R}$ is called \textit{\textbf{symplectic reduction}}.

\item[b)]  Given any lagrangian $L\in\Lag^W$ there exists another lagrangian $L_0\supset W$ such that $L\in\mathcal{A}_{L_0^\perp}\subset\Lag^W$. Hence $\Lag^W$ is an open subset of $\Lag$.

\item[c)] If $L_0\supset W$ is a lagrangian and $W^\perp$ is the orthogonal complement  of $W$ in $L$  then $\mathscr{R}(\mathcal{A}_{L_0^\perp})=\mathcal{A}_{JW^\perp }$. In these Arnold coordinates $\mathscr{R}(T)=P_{JW^\perp }TP_{JW^\perp }$, for every operator $T:L_0^\perp\rightarrow L_0$. Hence $\mathscr{R}$ is differentiable.

\end{itemize}
\end{prop}

\noindent
\textbf{Proof:}  $a)$ The case when $W$ is a lagrangian itself, is trivial since then $H_W={0}$. We will therefore suppose that $W$ is not lagrangian in what follows.

Let us  fix $L\in \Lag^W$. The first thing to notice is that $L\cap W^\omega\neq 0$. Indeed suppose $L\cap W^\omega=\{0\}$. Then $J(L\cap W^\omega)=\{0\}$ and so $(L+W)^\perp=\{0\}$ which implies that $L+W=H\oplus H$. Since  $L\cap W^\omega=\{0\}$  we get that $W$ is a linear complement of $L$ and because $W^\omega\supset W$ the only way $L\cap W^\omega=\{0\}$ is if $W^\omega =W$, that is $W$ is lagrangian, which is not our assumption. \par
 
  Let now $\ell:=P_{H_W}(L\cap W^\omega)$. We must have $\ell\neq \{0\}$ as well, since  ${P_{H_W}}|_{L\cap W^\omega}$ is injective because of the equality
   \[ 
   \Ker {P_{H_W}}|_{L\cap W^\omega}=L\cap W^\omega \cap (W\oplus JW)=L\cap W=\{0\}.
   \] 
   Moreover $J\ell=P_{H_W}(L^\perp\cap (JW)^\omega)$. This is true because $P_{H_W}J=JP_{H_W}$ which is another way of saying that $H_W$ is $J$ invariant. 
  
  We will check that $\ell\perp J\ell$. Let $x=y_1+z_1$ with $x\in L\cap W^\omega$, $y\in H_W$ and $z\in W$. Let also $x^\perp\in L^\perp\cap (JW)^\omega$ be decomposed as $x^\perp=y_2+Jz_2$ with $y_2\in H_W$ and $z_2\in W$. We notice that
   \[  W^\omega \ni x\perp Jz_2\in JW \mbox{ and }  (JW)^\omega\ni x^\perp\perp z_1\in W
   \]
   The next relation is now straightforward, thus proving the claim.
  \[ \langle y_1,y_2\rangle= \langle x-z_1,x^\perp-Jz_2\rangle=0
  \] 
  
  So $\ell$ is an isotropic subspace of $H_W$  and $L_0=J\ell+W$ is an isotropic space of $\hat{H}$  such that $L_0+L$ is closed and $L_0\cap L=\{0\}$. 
We will prove this claim.  Take $z=x+y\in (J\ell+W) \cap L$, with $x\in J\ell\subset H_W$ and $y\in W$. Then $z\in L\cap W^\omega=L\cap ( H_W\oplus W)$ which also means that $x=P_{H_W}(z) \in \ell$ and since $x\in J\ell$ we conclude that $x=0$ and therefore $z=y\in L\cap W=\{0\}$.\par

We notice now that $L_0$ has to be a maximal, isotropic space, hence lagrangian. If it were not maximal, we could repeat the whole process, with $L_0$ instead of $W$. This means that letting $H_{L_0}\subset H_W$ be the orthogonal complement of $L_0+JL_0$ and $\ell_0$ be the projection of $L\cap L_0^\omega\subset L\cap W^\omega$ onto $H_{L_0}$ then $\ell_0\subset H_{L_0}$ would be non-zero and also a subset of $\ell$. Since $\ell\perp H_{L_0}$,  we get a contradiction. \\

 $b)$  We only need to show that $\mathcal{A}_{L_0^\perp }\subset  \Lag^W $ where $L_0$ is as before and we are done. The finite dimensionality follows immediately from $L\cap W\subset L\cap L_0$. 

Notice that we have the set equality $L+W=L_0^\perp+W$ because if $L=\Gamma_T$ for some operator $T:L_0^\perp\rightarrow L_0$, then for every pair $(a,b)\in L_0^\perp+W$ the equation  $(x,Tx)+(0,y)=(a,b)$ has a unique solution.  The space $L_0^\perp+W$ is the orthogonal complement of $W^\perp$ in $\hat{H}$ and so it is closed. \\

$c)$ We want to see what happens when we restrict $\mathscr{R}$ to $\mathcal{A}_{L_0^\perp }$. First of all  we can see $W^\perp$ as a distinct lagrangian in $\Lag(H_W)$. 

We claim that $\mathscr{R}(\mathcal{A}_{L_0^\perp })=\Lag^{W^\perp}(H_W)$ and the last set is just the Arnold chart $\mathcal{A}_{JW^\perp}$ in $H_W$. To prove the claim we will compute the symplectic reduction of a graph $\Gamma_T\subset L_0^\perp\oplus L_0$ of an operator $T:L_0^\perp=JW\oplus JW^\perp\rightarrow L_0=W\oplus W^\perp$ with components:
\begin{center}
$T=\left(
\begin{array}{cc}
T_1& T_2\\
T_3 & T_4
\end{array}
\right)$
\end{center}
This is actually easy to do and the answer is $\mathscr{R}(\Gamma_T)=\Gamma_{T_4}\subset  JW^\perp\oplus W^\perp$.  This computation also proves the last claim.
\qed\\

\begin{remark} Given a  cofinite, isotropic space $W$ and a lagrangian $L$ such that $\dim {L\cap W}<\infty$ and $L+W$ closed, the symplectic reduction of $L$ is still well-defined. However the symplectic reduction as a map is not continuous on this set. We will see that it is continuous and in fact differentiable when restricted to each subspace of lagrangians with a fixed dimension of the intersection $L\cap W$. See section \ref {Gen Red}.
\end{remark}

The next set-theoretic equality is a useful by-product of the previous lemma.
\begin{lemma} \label{byproduct} If $L_0$ is a lagrangian which contains $W$, then $\mathcal{A}_{L_0^\perp}=\mathscr{R}^{-1}(A_{JW^\perp})\subset \Lag^W$, where $W^\perp$ is the orthogonal complement of $W$ in $L_0$.
\end{lemma}
\noindent
\textbf{Proof:} The inclusion $\subset$ is clear from the proof of the previous lemma. For the other one, let $L$ be a lagrangian such that $\mathscr{R}(L)\cap W^\perp=\{0\}$. 

 Consider now $x\in L\cap L_0=L\cap (W\oplus W^\perp)\subset L\cap W^\omega$. This means $x=a+b$ with $a\in W$ and $b\in W^\perp$ and so $b=P_{H_W}(x)\in\mathscr{R}(L)$. Therefore $b=0$ which implies $x=a\in L\cap W=\{0\}$.
 
 The sum $L+L_0$ is closed because it is the sum of a closed space $L+W$ and a finite dimensional one $W^\perp$. 
\qed\\

 We want to show next that the symplectic reduction is actually a linear fibration. First let us notice that we have a canonical section of $\mathscr{R}$, namely
\[\mathscr{S}:\Lag(H_W)\mapsto \Lag^W, \quad \quad~~~~~~ \ell\mapsto \ell\oplus JW.\]

Remember that by Lemma $\ref {tangent space} $ the tangent bundle of $\Lag^W$ is  identified with the bundle whose fiber at $L$ are just the self-adjoint operators $\Sym(L)$.  For every $\ell\in \Lag(H_W)$,  the lagrangian $\ell\oplus JW$ is clean with $W$. We saw in lemma $\ref {symp red} $, part $c)$ that in the Arnold chart of $L=\ell\oplus JW$, the symplectic reduction with $W$ has the simple expression of projection onto the $\ell\times \ell$ block. We conclude that $\Ker {d_{L}\mathscr{R}}$ consists of self-adjoint operators $S:L\rightarrow L$ with block decomposition:
\begin{center}
$S=\left(\begin{array}{cc}
0& S_2^*\\
S_2 & S_3
\end{array}
\right)$
\end{center}

We have the following equivalent of the tubular neighborhood theorem:

\begin{prop} \begin{itemize}

\item[a)] $(\Ker {d\mathscr{R}})\bigr |_{\Lag(H_W)}$ is a trivializable bundle over $\Lag(H_W)$.
\item [b)] Let ${\Lag(H_W)}$ be embedded in  $\Lag^W$ via $\mathscr{S}$.  Then the map

 \[ \mathscr{N}:(\Ker {d\mathscr{R}})\bigr |_{\Lag(H_W)}\mapsto \Lag^W, \quad \quad~~~~~~
(L,S)\mapsto\Gamma_{JS}\subset L\oplus L^\perp \]
 is a diffeomorphism which makes the diagram commutative:

\[ \xymatrix{ \ar @{} [dr] |{}
(\Ker{d\mathscr{R}})\bigr |_{\Lag(H_W)} \ar[rr]^{\quad\mathscr{N}} \ar[dr] && \Lag^W \ar[ld]_{\mathscr{R}} \\
 & \Lag(H_W)  &                    }
\]
\end{itemize}
\end{prop}
\noindent 
\textbf{Proof:} $a)$ In this infinite dimensional context we need to make sure first that \linebreak $(\Ker {d\mathscr{R}})\bigr |_{\Lag(H_W)}$ is a manifold. 

We fix  $\ell_0\subset\Lag{(H_W)}$ and let $L_0:=\ell_0\oplus JW=\mathscr{S}(\ell_0)$.

We will show that there is a natural homeomorphism  of   $(\Ker {d\mathscr{R}})\bigr |_{\Lag(H_W)}$ with a product \linebreak$ \Lag(H_W)\times\Sym^0(L^0)$, which commutes with the projections to $\Lag(H_W)$ and is linear on each fiber. Here $\Sym^0(L^0):=\Ker d_{L_0}\mathscr{R}$ is a Banach subspace of $\Sym(L_0)$.

The space  $\Sym^0(L_0)$ has a concrete description. It is the set of all self-adjoint operators on $L_0$ which are zero on the $\ell_0\times\ell_0$ block.

 Recall that the finite version of Arnold's theorem  tells us that for every lagrangian $\ell_1\in\Lag{(H_W)}$ there exists a canonical unitary isomorphism $\tilde{U}_{\ell_1}:\ell_0\rightarrow \ell_1$. Let $L_1:=\ell_1\oplus JW$.  The canonical unitary isomorphism $\tilde{U}_{L_1}:L_0\rightarrow L_1$ has a block decomposition 
 \begin{center}
$\tilde {U}_{L_1}:=\left(\begin{array}{cc}
{U}_{\ell_1} & 0\\
0& 1
\end{array}
\right)$
\end{center}
  as a map $\ell_0\oplus JW\rightarrow\ell_1\oplus JW$. 
  
  In these conditions, if the $\ell_1\times\ell_1$ block of a self-adjoint operator $S\in\Sym{(L_1)}$ is zero,   then the $\ell_0\times\ell_0$ block of $\tilde{U}_{L_1}^{-1} S\tilde{U}_{L_1}$ is zero as well.
  
  The map:
  \[ 
(\Ker {d\mathscr{R}})|_{\Lag(H_W)}\mapsto \Lag(H_W)\times \Sym^0(L^0), \quad \quad~~~~
(L_1,S)\mapsto(\ell_1, \tilde{U}_{L_1}^{-1}S\tilde{U}_{L_1}) 
  \]
  is the homeomorphism we were after. The continuity of the second component follows by noticing 
  as in corollary $\ref {tangent space} $ that $\tilde{U}_{L_1}^{-1}S\tilde{U}_{L_1}$ is the $L_0\times L_0$ block of $\mathscr{C}(U_{L_1})^{-1}S\mathscr{C}(U_{L_1})$. \\
 
$b)$  We will build an inverse for $\mathscr{N}$. Let $L\in \Lag^W$ and  let $\ell=\mathscr{R}(L)$ be its symplectic reduction. 

We consider the lagrangian $L_1=\ell\oplus JW$. It has the property that  $L+L_1^\perp=L+W+J\mathscr{R}(L)$ is closed and also $L\cap L_1^\perp=\{0\}$. 

 We prove the last claim. Let us take  $x\in L\cap W\oplus J\ell $ decomposed as $x=a+b+c$ where $a\in W,~b\in JW$  and $c\in H_W$. We have  $b=0$ because $x,a$ and $c$ are in $W\oplus H_W$. It follows that  $c$ is the projection of $x$ onto $H_W$ and by definition this is just $\mathscr{R}_L(x)\in \ell$. 
On the other hand $c$ is in $J\ell$ because $a+c\in W\oplus J\ell$ and $c\perp W$ and so $c=0$. Therefore $x=a\in L\cap W=\{0\}$ and we are done with the claim.
 
  By proposition $\ref {Ar. chart}$, $L$ is the graph of $JS$ with $S\in\Sym{(L_1)}$. We claim that the $\ell\times \ell$ block of $S$ is zero. 
  
  Notice first that $\mathscr{R}(L)=\mathscr{R}(L_1)=\ell$ and $L=\Gamma_{JS}$.  But the symplectic reduction of the graph of $JS$ with $W$ is the graph of $JS_1$ where $S_1$ is the $\ell\times\ell$ block  of $S$. Since the graph of $JS_1$ is $\ell$ we conclude that $S_1=0$.
   
    Therefore the inverse is:
\begin{center}
$L\rightarrow (\mathscr{R}(L), JP_{L_1^\perp }\circ(P_{L_1}|_L)^{-1})$
\end{center}

The differentiability of the maps is immediate when one works in the Arnold charts of $\ell$ on the base space of symplectic reduction and uses Lemma $\ref {byproduct}$.
\qed\\

 We will  go back now to our space $\Lag^-$. We fix a complete, decreasing flag of cofinite, closed subspaces  in $H^-$. 
\begin{center}
$H^-=W_0\supset W_1\supset W_2 \supset \ldots$
\end{center}

We will also use the associated increasing flag:
\[ \{0\}=W_0^\perp\subset W_1^\perp\subset W_2^\perp \ldots \subset H^-
\]

Let us briefly recall that every closed subspace of a lagrangian is isotropic.
\begin{lemma} \label {open exhaust}
For any flag of cofinite, isotropic subspaces of the vertical space $H^-$ we have: 
$$\bigcup_{i} \Lag^{W_i}=\Lag^- ,\mbox{ and } \Lag^{W_i}\subset \Lag^{W_{i+1}} \mbox{ for all } i$$
Hence $\Lag^-$ is an open subset of $\Lag$.
\end{lemma}
\noindent
\textbf{Proof:} The inclusion $\bigcup \Lag^{W_i}\subset\Lag^-$ is straightforward  since if $L\in \Lag^{W_i}$ then dim$~L\cap H^-\leq$ codim$~W_i=i$ and $L+H^-=L+W_i+W_i^\perp$ is closed since the orthogonal complement $W_i^\perp$ finite dimensional. \par
  Conversely the decreasing sequence of finite dimensional spaces $L\cap W_i$ has trivial intersection so it must be that there is an $i$ such that $L\cap W_i={0}$.  The fact that $L+W_i$ is closed for every $i$ follows by noticing that $W_i$ is commensurate with $H^-$ and using \ref{commensurate}.
  
  The inclusion $ \Lag^{W_i}\subset \Lag^{W_{i+1}}$ is straightforward.
 \qed\\
 
We now describe how a choice of a complete flag of cofinite subspaces for $H^-$ defines an atlas on $\Lag^-$.
 
 Let $W_{i}^\perp$ be the orthogonal complement of $W_i^\perp$ in $H^-$ and let $F_{i}$ be the orthogonal complement of $W_{i-1}^\perp$ in $W_{i}^\perp$ for $i\geq1$ (or the orthogonal complement of $W_i$ in $W_{i-1}$ for that matter). Each $F_i\subset H^-$ is one dimensional. For each $k$-tuple $I=\{i_1,i_2,\ldots ,i_k\}$ let 
$F_I:=\oplus_{i\in I } F_i$\index{$F_I$} and $F_{I^c}$ be the orthogonal complement of $F_I$ in $H^-$. Furthemore we will define 
\begin{center}
$H_I^+=F_I\oplus JF_{I^c}$\\
$H_I^-=JF_I\oplus F_{I^c}$
\end{center}
In other words each $H^\pm_I$ is a lagrangian consisting of a direct sum between a subspace of $H^+$ and one of $H^-$. We chose $\pm$ in order to  suggest that $H_I^\pm$ has a finite codimensional space in common with $H^\pm$.  In particular $H_I^\pm$ and $H^\pm$ are commensurate. 

\begin{lemma} \label{flag charts}
\begin{itemize}
\item[a)] For $n\geq\max\{i~|~i\in I\} $ we have $\mathcal{A}_{H_I^+}\subset \Lag^{W_n}$. 
\item[b)] The  Arnold charts $\mathcal{A}_{H_I^+}$ cover $\Lag^-$, that is :
$$\bigcup_{I } \mathcal{A}_{H_I^+}= \Lag^-$$
\end{itemize}
\end{lemma}
\noindent
\textbf{Proof:} At $a)$ notice that for $n\geq\max\{i~|~i\in I\} $ we have $F_{I^c}\supset W_n$ and so $H_I^-\supset W_n$ and therefore $\mathcal{A}_{H_I^+}\subset \Lag^{W_n}$.

For  $b)$ we write $T:H_I^+\rightarrow H_I^-$ as

\[T= \left(\begin{array}{cc}
A&B\\
C&D
\end{array}
\right)
\]
relative to the decompositions $H_I^+=F_I\oplus JF_{I^c}$ and $H_I^-=JF_I\oplus F_{I^c}$.
 It is easy to see that $\Gamma_T\cap H^-=\Ker C$. Now the sum $\Gamma_T+H^-$ is closed being the sum of a finite dimensional space and the set $\{v+Bv+Dv~|~v\in JF_{I^c}\}+F_{I^c}$.  The later set is closed because it is just the graph of $JF_{I^c}\oplus F_{I^c}\ni (a,b)\rightarrow Ba \in JF_I$. So $\mathcal{A}_{H_I^+}\subset Lag^-$. \par
 Conversely, we know by the previous lemma that $L\in \Lag^-$ is in some set $\Lag^{W_n}$. We let $H_{W_n}=W_n^\perp\oplus JW_n^\perp$. It is known (see \cite {N1}  ) that in the finite dimensional case the Arnold charts of $F_I\oplus JF_{I^{c}_n}$ cover $\Lag(H_{W_n})$. Here $I=\{i_1,i_2,\ldots ,i_k\}$ is a $k$-tuple in $\{1,2,\ldots ,n\}$ and $I^{c}_n$ is its complement. Now notice that $JW_n\oplus JF_{I^{c}_n}=JF_{I^c}$ and hence by Lemma $\ref {byproduct}$ we must have that $L\in \Lag^-$ is in some Arnold chart as above.  
\qed \\

In a very similar manner one can prove the next result:
\begin{lemma}  Let $V$ denote a finite codimensional subspace of $H^+$ and $V^\perp$ its orthogonal complement in $H^+$. An atlas of $\Lag^-$ is given by the collection of Arnold charts around $V\oplus JV^\perp$:
$$\bigcup_{V\subset H^+}\mathcal{A}_{V\oplus JV^\perp}=\Lag^-$$
\end{lemma}
\noindent
\begin{remark} This lemma, together with $\Lag^-_{\mathscr{I} }\subset \Lag^-$ finishes the proof that $\Lag^-_{\mathscr{I} }$ is a manifold modelled on $\Sym_\mathscr{I}$ with the charts described at the end of last section. In fact all the previous lemmas of this section are true, with $\Lag^-_{\mathscr{I} }$ replacing $\Lag^-$.
\end{remark}

 The following important theorem  is an application of symplectic reduction.
\begin{theorem} \label{weak homotopy equivalence}
The space $\Lag^-$ is weak homotopy equivalent of $U(\infty)$. The same is true about $\Lag^-_{\mathscr{I} }$.
\end{theorem}
\noindent
\textbf{Proof:} If we can prove that for a fixed $k$ there exists an $n$ big enough such that the pair $(\Lag^-,\Lag^{W_n})$ is $k$-connected than we are done because $\Lag^{W_n}$ is homotopy equivalent with $\Lag(H_{W_n})$ and this will imply that the induced map 
\[\Lag(\infty):=\lim_{n}\Lag(H_{W_n})~~\longmapsto~~\lim_{n}\Lag^{W_n}=\Lag^-
\]
is a weak homotopy equivalence. 

We have of course that $\Lag^{W_n}=\Lag^-\setminus \{L~|~\dim {L\cap W_n}\geq 1\} $ and we will see  in the next section that the set $\overline{Z}_{n+1}:=\{L~|~\dim {L\cap W_n}\geq 1\} $ is a finite codimensional stratified subset of $\Lag^-$ whose "highest" stratum has codimension $2n+1$.  We therefore fix $n>1/2(k-1)$ and show the induced map on homotopy groups
\[\pi_k(\Lag^{W_n})\mapsto \pi_k(\Lag^-)
\]
is an isomorphism.

Every continuous map $\sigma: S^k\rightarrow\Lag^-$  is contained  in an open set $\Lag^{W_N}$ for some $N>n$ big enough so one can deform it to map $S^k\rightarrow\Lag(H_{W_N})\hookrightarrow\Lag^-$  simply by composing with the symplectic reduction which is a deformation retract. The new map $\sigma_1:S^k\rightarrow \Lag(H_{W_N})$ can be deformed into a smooth map and can also be put into transversal position with $\overline{Z}_{n+1}\cap \Lag(H_{W_N})$ which is a Whitney stratified set of codimension $2n+1$ in the finite dimensional manifold $\Lag(H_{W_N})$. But for $k<2n+1$ this means that there is no intersection and hence the resulting map $\sigma_2$ has its image in $\Lag^{W_n}$. This proves the surjectiviy of the map on homotopy groups.

The injectivity follows by noticing that every map $I\times S^k\rightarrow \Lag^-$ can be deformed to a map $I\times S^k\rightarrow \Lag^{W_n}$ by the same type of argument as before for $2n>k$.

The same proof  works for $\Lag^-_{\mathscr{I} }$.\qed\\

We want to give an example of how symplectic reduction works in a concrete case:

\begin{example}[\textit{The universal family}] \label{ex univ family}
\index{universal family} \par
For each unitary map $U\in U(N)$, the Hilbert space $ L^{1,2}_U([0,1])$ is the completion of the space of smooth maps 
\[ C^{\infty }_U([0,1]):=\{f:[0,1]\rightarrow \mathbb{C}^N\}~|~f\in C^{\infty }([0,1]; \mathbb{C}^N), ~f(1)=Uf(0)\}\] in the $L^{1,2}$ norm:
\[ \|f\|_{1,2}^2=\int _0^1 |f'(t)|^2+|f(t)|^2~dt
\]
 The  differential operator:
\[
T_U: L^{1,2}_U([0,1])\rightarrow L^2([0,1]),  \quad \quad T_U=-i\frac{d}{dt}
\]
is a closed, densely defined, self-adjoint, elliptic (hence Fredholm) with compact resolvent. If we let $U\in U(N)$ vary we get a family of differential operators and by taking the switched graphs, a map:
\[ T: U(N)\rightarrow \Lag^-
\]
We will prove later that this family is  in fact differentiable, Corollary \ref{univ family}.

The isotropic space we choose to do symplectic reduction with, $W_N\subset H^-:=L^2([0,1])$ will be the orthogonal complement of the space $\mathbb{C}^N$ of constant functions. Notice that in general this is not a subspace of the domain of $T_U$. The annihilator, $W_N^\omega$ is the space $\mathbb{C}^N\oplus L^2([0,1])\subset L^2([0,1])\oplus L^2([0,1])$. 

First, let us check that $T_U$ is clean with $W_N$ for all $U$. This comes down to proving that the system
\[  \left\{\begin {array} {c}
T_U(\phi)=0 \\

  \int_0^1 \phi(t)~dt=0
  \end{array}\right.
\]
admits only the trivial solution. And that is easy to do.

We look now at the intersection of the switched graph of $T_U$ with the annihilator. This means solving the equations:
\[ \begin {array} {ccc}
T_U\phi = c  \\
\phi(1) = U\phi(0)
\end{array}
\]
with $c\in\mathbb{C}^N$.
We get $\phi(t)=ict+b$ and $ic+b=Ub$. So $c=i(1-U)b$ and one has solutions only for those constants $c$ that lie in the image of $i(1-U)$ in which case $\phi(t)=(U-1)bt+b$. We have to project this to the subspace of constant functions, which means computing the integral:
\[\int_0^1\phi(t)~dt= 1/2(U+1)b
\]

 In the end, the symplectic reduction of the switched graph of $T_U$ with $W_N$ is the subspace $L_U$ of $\mathbb{C^N}\oplus \mathbb{C}^N$ described by the following:
\[ x\in L_U \Leftrightarrow x=(i(1-U)b, 1/2(1+U)b) \mbox { for some } b\in \mathbb{C}^N
\]

This looks almost like the map that gives the finite dimensional Arnold isomorphism. The following is true. The map:

\[ 
U(N)\rightarrow U(N), \quad \quad \quad U\mapsto (1-3U)(3-U)^{-1}
\] 
is a diffeomorphism and in fact an involution.  If we let
  \[ U_1:=(1-3U)(3-U)^{-1} \quad \mbox{ and }\quad a:=i/4(3-U)b\]
   then 
\[ x=(i(1-U)b, 1/2(1+U)b)= ((1+U_1)a,-i(1-U_1)a)
\]
Recall that Arnold's  theorem says that the Cayley graph map
\[ 
 U(N)\rightarrow \Lag(\mathbb{C}^N\oplus\mathbb{C}^N), \quad \quad~~~~  U\rightarrow \Ran [v\rightarrow ((1+U)v,-i(1-U)v)]
\]
is a diffeomorphism. 

We conclude that the composition   of the inverse of the  Cayley graph map with the symplectic reduction of the universal family,   $(\Gamma^{\mathscr{C}})^{-1}\circ \mathscr{R}_{W_N}\circ T:U(N)\rightarrow U(N)$ is the involution of $U(N)$ given by
\[ U\rightarrow  (1-3U)(3-U)^{-1}
\] 

The following is a family of diffeomorphisms of the unitary group
\[ 
\mathbb{C}\times \mathbb{C}\setminus \{(\lambda,\mu)~|~|\lambda|=|\mu|\}\times U(N)\longmapsto U(N), \quad \quad       (\lambda,\mu, U)\longmapsto (\lambda-\mu U)(\bar{\mu}-\bar{\lambda}U)^{-1}
\]

The real hypersurface $|\lambda|=|\mu|$ splits $\mathbb{C}\times\mathbb{C}$ into two connected components \linebreak $(\mathbb{C}\times\mathbb{C})_\pm:=\{(\lambda,\mu)~|~sgn(|\mu|-|\lambda|)=\pm\}$. So $(1,3)$ can be connected with $(0,1)$ and therefore the previous involution is homotopy equivalent with the identity map.

We have therefore proved:
\begin {prop} The universal family is  homotopy equivalent with the inclusion $U(N)\hookrightarrow\Lag^-$ that the Arnold isomorphism provides.
\end{prop}
\end{example}

\noindent
\begin{remark} In the discussion above the expression for Arnold's isomorphism plays a very important role in finding the homotopy class of the universal family. This is because the group of diffeomorphisms of $U(N)$ is not connected. Hence if instead of the Cayley transform $(1+U,-i(1-U)) $ we would have chosen the minus Cayley transform 
$(1+U,i(1-U))$ for the Arnold isomorphism, then what we called the universal family above would have been equivalent with the inclusion $U(N)\hookrightarrow\Lag^-$ composed with the conjugation map $U(N)\rightarrow U(N)$. Nevertheless, it is true that in such a situation, the family $T_U=i\frac{d}{dt}$ with the same boundary condition as above would have been equivalent with the inclusion. Notice also that changing the boundary condition, for example to $U\phi(1)=\phi(0)$, changes the homotopy type of the universal family to its conjugate.
\end{remark}

When one looks back at Arnold's theorem, a legitimate question is certainly what does symplectic reduction mean for the unitary group. We conclude this section with the answer.

\begin{prop}
Let $U\in \mathscr{U}(H)$ be a unitary operator such that $1+U$ is Fredholm and let $W\subset H$ be a finite codimensional, closed subspace. Suppose  that $\Ker{(1+U)}\cap W=\{0\}$. Let $W^\perp$ be the orthogonal complement of $W$ and let 
\[ U=\left (\begin {array}{cc}
 X & Y \\
 Z & T
\end{array}
\right )
\]
be the block decomposition of $U$ relative $H=W\oplus W^\perp$. Then $1+X:W\rightarrow W$ is invertible  and the operator  $\mathscr{R}(U):W^\perp\rightarrow W^\perp$:
\[ \mathscr{R}(U):= T-Z(1+X)^{-1}Y
\]
is unitary. Moreover $\Ker{(1+\mathscr{R}(U))}=P_{W^\perp}(\Ker{(1+U)})$ and in particular \linebreak $\dim{\Ker {(1+\mathscr{R}(U))}}=\dim{\Ker{(1+U)}}$.
\end{prop}
\noindent
\textbf{Proof:}  We write $1+X=P_W \circ (1+U) |_W$. The operator $(1+U)|_W$ is Fredholm as a composition of two Fredholm operators, $1+U$ and $i_W:W\hookrightarrow H$. So the image of $(1+U)|_W$ is closed and it has the same codimension as $W$ in $H$. This follows from 
\[ \ind {(1+U) i_W}=\ind {(1+U)}+ \ind{ i_W}
\]
and from noticing that $\ind{(1+U)}= \{0\}$ since $\Ker{(1+U)}=\Ker{U(1+U^*)}=\Ker{(1+U^*)}$. 

If we can show that $\Range{(1+U)|_W}\cap W^\perp=\{0\}$ then $W^\perp$ would be an algebraic complement of $\Range{(1+U)|_W}$ and so this would imply that $P_W\circ (1+U) |_W$ is invertible. Let now $w\in W$ be such that 
\[ \langle (1+U)w,v\rangle=0 \quad\forall v\in W
\]
In particular this is true for $v=w$ and so:

\begin{equation} \label{equ}
 \|w\|^2+\langle Uw,w\rangle =0
\end{equation}
or in other words $|\langle Uw,w\rangle|=\|w\|^2$ which  is the equality case in the Cauchy inequality since $U$ is unitary and so $Uw=\lambda w$. Going back to $(\ref {equ})$ one sees that $\lambda =-1$ and so $w\in W\cap\Ker{(1+U)}=\{0\} $.

Let us check now that the resulting operator is unitary. First let us write down the relations one gets because $U$ is unitary.
\[ \left \{ \begin{array}{ccc}
X^*X+Z^*Z &=&1\\
X^*Y+Z^*T&=&0\\
Y^*Y +T^*T &=&1\\
Y^*X+T^*Z&=&0
\end{array}
\right .
\]

So:
\[ (T^*-Y^*(1+X^*)^{-1}Z^*)(T-Z(1+X)^{-1}Y)= 1-Y^*Y+ Y^*(1+X^*)^{-1}X^*Y+
\]
\[+ Y^*X(1+X)^{-1}Y+Y^*(1+X^*)^{-1}(1-X^*X)(1+X)^{-1}Y=
\]
\[ 1- Y^*[(-1+(1+X^*)^{-1}X^*)+X(1+X)^{-1}+(1+X^*)^{-1}(1-X^*X)(1+X)^{-1}]Y
\]
The sum in the square brackets is equal to:
\[ -(1+X^*)^{-1}+X(1+X)^{-1}+(1+X^*)^{-1}(1+X)^{-1}-(1+X^*)^{-1}X^*X(1+X)^{-1}=
\]
\[=(1+X^*)^{-1}(-1+(1+X)^{-1})+(1-(1+X^*)^{-1}X^*)X(1+X)^{-1}=
\]
\[ =-(1+X^*)^{-1}X(1+X)^{-1}+(1+X^*)^{-1}X(1+X)^{-1}=0
\]
and we are done proving that $\mathscr{R}(U)$ is unitary. 

Let us take $(w,w^\perp)\in\Ker{(1+U)}$. This means that
\[\left \{ \begin{array} {ccc} (1+X)w+Yw^\perp&=&0\\
Zw+(1+T)w^\perp&=&0
\end{array}
\right.
\]
Since $1+X$ is invertible one gets that $w=-(1+X)^{-1}Yw^\perp$ and therefore 
\[ (1+T)w^\perp-Z(1+X)^{-1}Yw^\perp=0
\]
This means, of course, that $w^\perp\in\Ker{(1+\mathscr{R}(U))}$. 

Conversely if $w^\perp\in\Ker{(1+\mathscr{R}(U))},$ then it is straightforward to see that \linebreak $(-(1+X)^{-1}Yw^\perp,w^\perp)\in\Ker{(1+U)}$.

The fact that the dimensions are equal follows from  the injectivity of $P_{W^\perp}|_{\Ker{(1+U)}}$. \qed

\begin{remark} If $\lambda+U$ is Fredholm, where $|\lambda|=1$ is a unit complex number and $W$ is closed, cofinite such that $\Ker{(\lambda +U)}\cap W=\{0\}$ then one should replace $\mathscr{R}(U)$ with $\mathscr{R}_\lambda(U)=T-Z(\lambda +X)^{-1}Y$. The conclusion is that $\mathscr{R}_\lambda(U)$ is unitary and the relation $\Ker {(\lambda+\mathscr{R}_\lambda (U))}=P_{W^\perp}(\Ker{(\lambda +U)})$ holds, as one can easily check. 
\end{remark}

\begin{definition} Let $U\in U(H)$ be a unitary operator such that $1+U$ is Fredholm and let $W$ be a finite codimensional subspace. Then $U$ is said to be clean with $W$ if $\Ker{(1+U)}\cap W=\{0\}$.
\end{definition}
\noindent
\textbf{Example:} Every unitary matrix $U\in U(2)$ can be written in a unique way as
\[ U=\left ( \begin{array} {cc}
z & -\lambda\bar{w}\\
w & \lambda\bar{z}
\end{array}
\right )
\]
where $(\lambda,z,w)\in\mathbb{C}^3$ such that $|\lambda|=1$ and $|z|^2+|w|^2=1$. Here $\lambda=\det{U}$. A unitary matrix $U$ is clean   with $W:=\mathbb{C}\oplus 0$ if and only if $z\neq -1$. 

The reduction map associates to every $U\in U(2)\setminus\{U~|~z=-1\}$ the unit complex number:
\[ \lambda\bar{z}-w\frac{1}{1+z}(-\lambda \bar{w})=\lambda\bar{z}+\lambda\frac{1-|z|^2}{1+z}=\lambda \frac{1+\bar{z}}{1+z}
\]

When $\lambda =1$ this descends to a map
\[ \mathscr {R}: SU(2)\setminus\{{-1}\}\mapsto S^1\setminus \{{-1}\}, \quad \quad~~~~~ \left ( \begin{array} {cc}
z & -\bar{w}\\
w & \bar{z}
\end{array}
\right ) \mapsto \frac{1+\bar {z}}{1+z}
\]

\section{The topology of the vertical Lagrangian Grassmannian}
\setcounter{equation}{0} 

In this chapter we describe a Schubert cell stratification of $\Lag^-$ that parallels the one given in the finite dimensional case by Nicolaescu in \cite {N1}. These cells are finite codimensional submanifolds of the Lagrangian Grassmannian whose closures naturally determine cohomology classes that correspond to the generators of the cohomology group of $\Lag^-$. For the "canonical" generators of the cohomology \textit{ring} of $\Lag^-$ there is a different stratification which is more suitable for doing intersection theory.

\bigskip

\subsection{Schubert cells and  varieties} 

The topological structure of $\Lag^-$ is intimately connected with the structure of the finite Lagrangian Grassmanians which are nothing else but the classical unitary groups. A detailed topological study of these spaces has been undertaken by Nicolaescu  in \cite {N1}.  In that paper, the author shows that the Poincar\'e duals of the generators of the cohomology group of  $U(n)$ can be represented by integral currents supported by semialgebraic varieties. That approach is not available in our infinite-dimensional context. However we have on our side symplectic reduction that reduces most of the problems to their finite dimensional counterpart.

In section \ref{Symp red} we introduced a complete, decreasing flag:

\[  H^-=W_0\supset W_1\supset W_2\supset\ldots 
 \]
 We now  fix a Hilbert basis $\{f_1,f_2,\ldots \}$ of $H^-$ such that $W_n^\perp=\langle f_1,f_2,\ldots ,f_n\rangle$ and we set $e_i:=Jf_i$ such that $\{e_1,e_2,\ldots \}$ is an  orthonormal basis of $H^+$.
 
 To every $k$-tuple of positive integers $I=\{i_1<i_2<\ldots <i_k\}$ we associated the following vector spaces 
  
 \begin{center}
$F_I=\langle f_i~|~i\in I\rangle$, $F_{I^c}=\overline{\langle f_i~|~i\in I^c\rangle}$ and $ H_I^+=\langle f_i~|~i\in I\rangle\oplus \overline {\langle e_j~|~j\in I^c\rangle } $
\end{center}

\begin{definition}Let  $I=\{i_1<i_2<\ldots <i_k\}$ be a $k$-tuple of positive integers. Set $i_0:=0$ and $i_{k+1}:=\infty$. The weight of the $k$-tuple is the integer:
\[ N_I:=\sum_{i\in I}(2i-1)
\]

The \textit{\textbf{Schubert cell}}\index{Schubert cell} of type $I$ denoted $Z_I$ is a subset of $\Lag^-$ defined by the following incidence relations with respect to a fixed flag 
\begin{center}
 $Z_I=\{L\in \Lag^-~|~$dim$~L\cap W_j=k-p,~\forall~0\leq p\leq k, ~\forall~ j$ such that $i_p\leq j< i_{p+1} \}$
\end{center}
\end{definition}
\begin{remark} One way to look at the incidence relations is by thinking that the $k$-tuple $(i_1,i_2,\ldots ,i_k)$ records the "nodes" in the flag where the dimension of the intersection with the lagrangian $L$ drops by one. \qed
\end{remark}

\noindent
\begin{remark}
Notice that the orthogonal complement $W_n^\perp$ of $W_n$ is naturally a lagrangian in $H_{W_n}:= W_n^\perp\oplus JW_n^\perp$ and $W_n\subset H_{W_n}$ will play the role of $H^-$. The flag $W_0=H^-\supset W_1\supset W_2$ induces a complete, decreasing flag of $W_n^\perp$:
\begin{center} $\tilde{W_0}:=W_n^\perp\supset \tilde{W_1}:=W_1/W_n\supset\ldots \supset\tilde{W_n}:= W_n/W_n=\{0\}$
\end{center}  
We let  $Z_{I}(n)$  be the Schubert cell in  $\Lag(H_{W_n})$ described by the same incidence relations as the sets $Z_I$ above with $\tilde{W_i}$ replacing $W_i$.
\qed
\end{remark}

The following description of Schubert cells proves that they are actually Banach spaces when regarded in the right charts.
 \begin{prop} \label{Schubert cell description}
The {Schubert cell}  $Z_I$ is a closed vector subspace of the Arnold chart $\mathcal{A}_{H^+_I}$ of codimension $N_I$. More precisely $\Gamma_{JA}\in\mathscr{A}_{H_I^+}\cap Z_I$ if and only if the bounded self-adjoint operator satisfies the linear equations
\[ 
 \langle Af_i, f_j \rangle=0, ~\forall ~j\leq i, ~i,j\in I \]
 \[ \langle Af_i, e_j \rangle=0, ~\forall ~j\leq i, ~i\in I, ~j\in I^c
 \]
\end{prop}
\noindent
\textbf{Proof:} We will show first that $ Z_I\subset\mathcal{A}_{H^+_I}$.  Let $L\in Z_I$. Notice that $(L,H_I^-)$ is a Fredholm pair  by Proposition \ref{commensurate} since $H_I^-$ is commensurate with $H^-$.

 We will show that $L\cap H_I^-=\{0\}$ thus proving that $L=\Gamma_{JA}\in\mathcal{A}_{H_I^+}$ with $A\in\Sym{H_I^+}$.
 
  Let us remark that $L\cap F_{I^c}=\{0\}$ because otherwise the dimension of  $L\cap W_j$ would drop at "nodes" other than $i_1,i_2, \ldots i_k$, (take $v=\sum_{j\in I^c} a_jf_j\in L\cap F_{I^c}$ with $p=\min{\{j\in I^c ~|~a_j\neq 0\}}$ then $v\in L\cap W_{p-1}\setminus L\cap W_p$). This is saying that $L\cap H^-$ is the graph of an operator $T:F_{I}\rightarrow F_{I^c}$.
 
   To see that $L\cap H_I^-=\{0\}$, let $x=v_1+v_2\in L\cap JF_{I}\oplus F_{I^c}$. Then $Jx\in L^\perp$ and so $\langle Jx, w+Tw\rangle=0$, for all $w\in F_{I}$. This implies:
   \[
   \langle Jv_1,w\rangle=0 ~\forall ~w\in F_{I}
   \] 
   We get $v_1=0$ and so $x=v_2\in L\cap F_{I^c}={0}$, thus finishing the proof that $L=\Gamma_{JA}\in\mathcal{A}_{H_I^+}$.
   
 Let  now $A\in \Sym{(H_I^+)}$ such that $\Gamma_{JA}\in Z_I$ and let 
\begin{center}
$A=\left( \begin{array}{cc}
A_1 & A_2\\
A_3 & A_4
\end{array}
\right)$
\end{center}
be the block decomposition of $A$ relative to $H_I^+=F_I\oplus JF_{I^c}$.

One checks immediately that the intersection $\Gamma_{JA}\cap H^-$ is just the graph of the restriction $JA_3|_{\Ker{A_1}}$ which has the same dimension as $\Ker{A_1}\subset F_{I}$. Since $F_I$ has dimension $k$, one concludes that 
\[  \dim\Gamma_{JA}\cap H^-=k \Longleftrightarrow A_1=0 \Longleftrightarrow \langle Af_i, f_j \rangle=0, ~\forall ~j\leq i, ~i,j\in I\] 
 
To prove the rest of the relations, i.e., $\langle Af_i, e_j \rangle=0, ~\forall i\in I, ~j\in I^c, ~j\leq i$ we observe first that 
\[   \Gamma_{JA}\cap H^-=\Gamma_{JA_3} \mbox{ and } \langle Af_i, e_j \rangle=-\langle JA_3f_i,f_j\rangle\].
   The graph of $JA_3=T:F_{I}\rightarrow F_{I^c}$ satisfies the incidence relations if and only if the required coefficients vanish, otherwise we would have dimension drops at the wrong places again.\qed \\\

\noindent
\begin{remark} For every two-sided symmetrically normed ideal $\mathscr{I}$ we can define $Z_I(\mathscr{I})=Z_I\cap \Lag^-_{\mathscr{I} }$. Since the next results are true for $Z_I$, as well as for $Z_I(\mathscr{I})$ making only the minimal changes, we choose to work with $Z_I$ to keep the indices to a minimum.
\end{remark}

Notice that $Z_I\subset \Lag^{W_n}$ for all $n\geq\max{\{i~|~i\in I \} }$ so we could look at the symplectic reduction of $Z_I$.   We record the obvious:

\begin{lemma}\label{stronger}
 For $~n\geq\max{\{i~|~i\in I\} }~$the $~$symplectic $~$reduction \linebreak $\mathscr{R}:\Lag^{W_n}\rightarrow \Lag(H_W)$ takes $Z_I\subset \Lag^{W_n}$ to $Z_{I}(n)$. The stronger $\mathscr{R}^{-1}(Z_{I}(n))=Z_I$ is also true.
\end{lemma}
\noindent
\textbf{Proof:} For $n\geq\max{\{i~|~i\in I\} }$ we have that $F_I\subset W_n^\perp$ and so $\mathcal{A}_{H_I^+}\subset Lag^{W_n}$ (see Lemma \ref{flag charts}). The reduction of the Arnold chart around $H_I^+$ is the Arnold chart around $H_I^+(n)$.  The reduction in the Arnold chart being just the projection, the lemma easily follows.  \qed

\begin{definition} For every $k$-tuple $I=\{i_1,i_2,\ldots ,i_k\}$ the \textit{\textbf{Schubert variety}}\index{Schubert variety}  is the closure of $Z_I$ in $\Lag^-$, denoted $\overline{Z}_I$.

\end{definition}
\begin{lemma}
The Schubert variety $\overline{Z}_{I}$ can be described by the following incidence relations:
\begin{center}
$\overline{Z}_I=\{L\in \Lag^-~|~$dim$~L\cap W_j\geq k-p,~\forall~0\leq p\leq k, ~\forall~ j$ such that $i_p\leq j< i_{p+1}$  where $ i_{0}=0, ~i_{k+1}=\infty$ and $i_p\in I,~\forall~1\leq p\leq k\}$
\end{center}
\end{lemma}
\noindent
\textbf{Proof:} The fact that the closure is included in the right hand side is a consequence of the upper semi-continuity of the functions:
\begin{center}
$L\rightarrow \dim{L\cap H^-},~~~~L\rightarrow\dim{L\cap W_{i_1}},~~\ldots~~L\rightarrow\dim{L\cap W_{i_k}}$
\end{center}

Conversely, let us notice that for $n$ big enough we have the following obvious equalities \linebreak $\Lag^{W_n}\cap\overline{Z}_I= \mathbf{cl}_n(Z_I)=\mathscr{R}^{-1}(\overline{Z_I(n)})$ where $\mathbf{cl}_n(Z_I)$ is the closure of $Z_I$ in $\Lag^{W_n}$. 
 Now, a lagrangian that satisfies the incidence relations in the lemma is in some $\Lag^{W_n}$ and its reduced space will satisfy the same incidence relations with respect to the flag $\tilde{W}\supset\tilde{W}_1\supset\ldots \supset \tilde{W}_n$. But this means it is in $\overline{Z_I(n)}$, since the finite version of the lemma is true by a result from \cite {N1}. This concludes the proof. \qed \\\par

\noindent
\begin{remark} The Schubert variety is not included in any of the clean sets $\Lag^{W_n}$. Nevertheless, the intersection has a very simple description: $\overline{Z}_I\cap Lag^{W_n}=\mathscr{R}^{-1}(\overline{Z}_I(n))$.
\end{remark}

We can now describe the strata in the Schubert variety $\overline{Z}_I$. Notice first that if $Z_J\subset \overline{Z}_I$ then $| J|\geq |I|$ since $| J|=\dim{L\cap H^-}$ for every $L\in  Z_J$. Say $J=\{j_1<j_2<\ldots <j_l\}$ and $I=\{i_1<i_2<\ldots<i_k\}$ with $l\geq k$. We deduce that $i_1\leq j_{l-k+1}$ since $j_{l-k+1}$ records the node where the dimension of the intersection of $L\in Z_J$ with the flag  drops to $k-1$ and similarly $i_s\leq j_{l-k+s}$ for all $1\leq s\le k$. We record this:
\begin{lemma}
\begin {itemize}
\item[a)] If $Z_J\subset \overline{Z}_I$ then $ |J|=l\geq k=| I|$ and $i_s\leq j_{l-k+s}$ for all $1\leq s\le k$.
\item[b)] If $Z_J\subset\overline{Z}_I$ has codimension $N_I+1$ in $\Lag^-$ where $N_I=\sum_{i\in I }(2i-1)$ is the codimension of $Z_I$ in $\Lag^-$ then $ |J|=k+1$, $j_1=1$ and $j_{s+1}=i_s$ for all $1\leq s\le k$.
\end{itemize}
\end{lemma}

In the proof of the Theorem \ref{weak homotopy equivalence} we used the following:
\begin{corollary} \label{fund Sch}
The fundamental Schubert variety $\overline{Z}_n$ can be described by the simple incidence relation:
\[ \overline{Z}_n=\{L~|~\dim{(L\cap W_{n-1})}\geq 1\}
\]
\end{corollary}
\noindent
\textbf{Proof:} Let $J=\{j_1,\ldots ,j_l\}$, $L\in Z_J$ and $Z_J\subset\overline{Z}_n$. The previous lemma tells us that $j_l\geq n$ and so the node where the dimension of the intersection of $L$ with the flag drops to 0 is bigger than $n-1$. This proves the $''\subset''$ inclusion. The other inclusion is obvious.
\qed\\

\bigskip

\subsection{The cohomology ring and geometrical representatives}\label{cohomology ring and geometrical representatives}

Our plan is to define geometrical representatives for the most important cohomology classes of $\Lag^-$, i.e., for certain canonical generators of its cohomology ring. The candidates are of course the stratified spaces $\overline{Z}_I$.     We identify the class whose underlying space is $\overline{Z}_I$ with a class expressed in terms of these generators.  When pulled-back to oriented, closed manifolds, via suitable maps, each of these cohomology classes has  a Poincare dual that is nothing else but the homology class determined by the preimage of the set $\overline{Z}_I$ together with an induced orientation. 

We first introduce the canonical generators for $H^*(\Lag^-)$. 

Following \cite{N1}, the groups $U(n)$ have  canonically defined cohomology classes $x_i\in H^{2i-1}(U(n),\mathbb{Z})$.  On the product $S^1\times U(n)$ there is a canonically defined complex vector bundle $E_n$ of rank $n$ called the universal bundle. \index{universal bundle}
  The  bundle is obtained by modding out  the $\mathbb{Z}$ action on the $\mathbb{Z}$-equivariant bundle:
  \[ \mathbb{R}\times U(n)\times\mathbb{C}^n\rightarrow \mathbb{R}\times U(n)
  \]
  The action on the total space is given by 
  \[ k(t,U,v):=(t+k, U, U^kv),\quad \quad \forall (t,U,v)\in \mathbb{R}\times U(n)\times\mathbb{C}^n, \quad  k\in\mathbb{Z}
  \]
 whereas on the base space, $\mathbb{Z}$ acts in the obvious way on the $\mathbb{R}$ component.
  The classes $x_i$ are transgressions of the Chern classes of $E_n$, i.e.,
  \[ x_{i}(n):=\int_{S^1} c_i(E_n)
  \]
  The classes $x_i(n)\in H^{2i-1}(U(n))$, $1\leq i\leq n$, generate the cohomology ring of $U(n)$, i.e.,
  \[ H^*(U(n),\mathbb{Z}) \simeq \Lambda(x_1,\ldots,x_n)\]
  It is not hard to see that via the canonical inclusion
  \[ S^1\times U(n)\hookrightarrow S^1\times U(n+1)
  \]
  the bundle $E_{n+1}$ pulls-back to give a bundle isomorphic to $E_n\oplus\underline{\mathbb{C}}$ hence the class $x_i(n+1)$ pulls back to $x_i(n)$. Since the canonical inclusion
  \[ i_n:U(n)\rightarrow U(\infty)
  \]
  induces an isomorphism in cohomology in degree $i<1/2(n+1)$ we conclude that there exists a \textit {unique} class
    \[ x_i\in H^{2i-1}(U(\infty)), \quad \mbox{ such that } \quad i_n^*(x_i)=x_i(n),\quad \forall n\]
    \begin{remark} The classes $x_i$ are in a suitable sense, which we will not care to make precise, the transgressions of the Chern classes for the infinite dimensional vector bundle associated to the principal $U(\infty)$ over  $S^1\wedge U(\infty)$ whose classifying map is given by the Bott periodicity map
    \[ S^1\wedge U(\infty)\rightarrow BU(\infty)
    \]  
    \end{remark}
    \begin{remark}\label{unitary equivalence} Let us note that if $V$ and $W$ are two  complex vector spaces with hermitian metrics then any unitary isomorphism between them $U:V\rightarrow W$ induces the same isomorphism on the cohomology rings of the unitary groups $U(V)$ and $U(W)$ simply because every unitary map on a vector space is homotopic with the identity. This is why we take a less invariant point of view.
\end{remark}
      
We now  fix a  complete, decreasing flag on $H^-$
\[  H^-\supset W_1\supset W_2\supset\ldots
\] 

A flag defines a natural subset $\Lag(\infty)\subset\Lag^-$, namely
\[ \Lag(\infty):=\lim_{n}\Lag(H_{W_n})
\]
where $H_{W_n}=JW_n^\perp\oplus W_n$.  Each $\Lag(H_{W_n})$ can be identified in a canonical way via Arnold's theorem with $U(JW_n^\perp)$ (inside $H_{W_n}$ it is $JW_n^\perp$ that plays the role of the horizontal lagrangian). On the other hand, a choice of  basis $\langle e_1,e_2, \ldots e_n\rangle$ of $JW_n^\perp$ identifies the unitary group $U(JW_n^\perp)$ with $U(n)$.   Remark \ref{unitary equivalence}  guarantees that we get canonically defined  cohomology classes  in $U(JW_n^\perp)$, $\Lag(H_{W_n})$ and $\Lag(\infty)$.

  \noindent
   \textbf{Notation:} We will use the same notations $x_i(n)$ and $x_i$ for the cohomology classes one gets in $\Lag(H_{W_n})$ and $\Lag(\infty)$ via the Cayley graph map.

  Theorem \ref{weak homotopy equivalence}  justifies the following definition
  \begin{definition} The fundamental \textit{\textbf{transgression classes}}\index{transgression classes} on $\Lag^-$ are the unique cohomology classes $z_i\in H^{2i-1}(\Lag^-)$ that   pull-back  to the classes $x_i\in H^{2i-1}(\Lag(\infty))$ via the  weak homotopy equivalence:
  \[ \xymatrix {  \Lag(\infty) \ar@{^{(}->}^{i}[r] & \Lag^-} , \quad \quad \quad x_i=i^*(z_i)
  \]
  where $i$ is the natural inclusion map. 
  
  For every set of positive integers $I=\{i_1,\ldots i_k\}$, define the product class  $z_I\in H^{N_I}(\Lag^-)$ to be the cup product of fundamental transgression classes:
  \[z_I=z_{i_1}\wedge z_{i_2}\wedge \ldots \wedge z_{i_k}
  \]
  \end{definition}
  \begin{remark} The classes $z_I$ are uniquely characterized by the property that
   \[ i_n^*(z_I)=x_I(n)
   \]
   for all natural inclusions $i_n:\Lag_{H_{W_n}}\rightarrow\Lag^-$,  for all $n$ big enough.
  \end{remark}

We now turn to our Schubert varieties $\overline{Z}_I$. In the Appendix \ref{appendix B} we describe in detail how can one build cohomology classes out of  a cooriented quasi-manifold. We summarize the main definitions and procedures:

 \begin{definition} \label{def: quasisubman} Let  $X$ be a Banach manifold. A \textit{\textbf{quasi-submanifold}}\index{quasi-submanifold} of  $X$ of codimension  $c$ is a closed subset $W\subset X$ together with a decreasing filtration by closed subsets
 \[
 W=\eF^0\supset\eF^1\supset \eF^2\supset \eF^3\subset \cdots
 \]
 such that  the following hold.
 
 \begin{itemize}
 
 \item $\eF^1=\eF^2$.
 
 \item The strata $\eS^k=\eF^k\setminus \eF^{k+1}$,  are submanifolds of $X$ of codimension $k+c$.
 
 \end{itemize}
 
 The    quasi-submanifold is called \textit{\textbf{coorientable}}\index{quasi-submanifold ! coorientable} if $\eS^0$  is coorientable. A coorientation of  a quasi-submanifold is then a  coorientation of its top stratum. 
 \end{definition}

  The main ingredients to define a cohomology class from a coorientable quasi-manifold  are:
\begin{itemize} 
\item A Thom isomorphism of the top dimensional stratum $\eS^0$, which is  a submanifold  and closed subset  of $X\setminus \eF^2$. This depends on the choice of a coorientation.
\item An extension isomorphism  in cohomology, over the singular stratum $\eF^2$, which exists because $\eF^2$ has codimension at least two bigger than $\eS^0$.
\end{itemize}

 A choice of a coorientation $\bom$ of the top stratum  defines a Thom isomorphism (of the closed submanifold $\eS^0\subset X\setminus \eF^2$):
\[  H^0(\eS^0)\simeq H^c(X\setminus \eF^2,X\setminus \eF^0)
\]
On the other hand by Proposition \ref{vanish codim} the quasi-submanifold $\eF^2$ has homological codimension at least $c+2$ and so the restriction map
\[ H^c(X)\rightarrow H^c(X\setminus\eF^2)
\]
is an isomorphism.  The cohomology class determined by the pair $[W,\bom]$ is the image of $1\in H^0(\eS^0)$ via the composition:
\[ H^0(\eS^0)\simeq H^c(X\setminus \eF^2,X\setminus \eF^0)\rightarrow H^c(X\setminus \eF^2)\simeq H^c(X)
\]
\begin{remark} \label{diff strat} A legitimate question is what role does the filtration play in the definition of the cohomology class? In the appendix \ref{appendix B} we show that if a  quasi-submanifold $W$ comes with two different filtrations $(W,\eF)$ and $(W,\eG)$  which have common refinement $(W,\eH)$, where by refinement we understand that $\eH^2\subset \eF^2\cup \eG^2$ and the coorientation on $W\setminus\eH^2$ restricts to the coorientations of $W\setminus\eF^2$ and $W\setminus\eG^2$ then they define the same cohomology class. It is possible that any two filtrations of  a quasi-manifold have a common refinement. However we could not prove that. 
\end{remark}
\vskip0.2in

To state the next result  we introduce a bit of notation and terminology.

\noindent
\textbf{Notation:} 
\[ Z_I^\circ:=Z_I\cup Z_{I\cup{1}}, \quad \quad\quad\quad \partial Z_I:= \overline{Z}_I\setminus  Z_I^\circ
\]

\begin{definition} We call the \textit{\textbf{standard filtration}}\index{standard filtration} on $\overline{Z}_I$ the following

\[ \eF^0:=\overline{Z}_I~; \quad \quad \eF^1=\eF^2:=\partial Z_I~;\quad \quad \eF^k:=\bigcup_{\substack{Z_J\subset \overline{Z}_I\\ N_J\geq N_I+k}} Z_J~.
\]

\end{definition}

\begin{theorem}\label{Schubert and canonical} The standard stratification on the Schubert variety $\overline{Z}_I$ turns it into a coorientable quasi-submanifold of $\Lag^-$ of codimension $N_I$. There exists a canonical choice of a coorientation $\omega_I$ on the top stratum such that  the following equality of cohomology classes holds
\[ [\overline{Z}_I,\omega_I]=z_I\]
\end{theorem}

Before we go into the proof, a short digression on coorientation is necessary. 

In \cite{N1}, Nicolaescu  uses the theory of currents to build out of the finite dimensional Schubert variety, $\overline{Z}_I(n)\subset \Lag(n)$, endowed with an orientation, a homology class which he shows is Poincare dual to the class $x_i(n)$, which is seen as a class in $\Lag(n)$ via the Arnold isomorphism.   
We summarize the main results:

\begin{prop} The sets $Z_I(n)^\circ$ are orientable smooth, subanalytic manifolds of codimension $N_I$ in $\Lag(n)$ .
\end{prop}
\noindent
\textbf{Proof:} See \cite{N1}, Lemma 5.7.
\qed\\

\noindent
\textbf{Orientation conventions I:} \index{Orientation conventions I} The space of unitary operators $U(n)$ is naturally oriented as follows. The Lie algebra $\mathfrak{u}(n)$ is the set of $n\times n$ skew-symmetric matrices. We identify it with the set of self-adjoint matrices $\Sym(n)$ by the map:
\[ \mathfrak{u}(n)\rightarrow \Sym(n),\quad \quad \quad B\mapsto -iB
\]
 Let 
 \[ \theta_1,\ldots \theta_n, (\alpha_{ij})_{1\leq i<j\leq n} ,(\beta_{ij})_{1\leq i<j\leq n} : \Sym(n)\rightarrow \mathbb{R}
 \] be the  linear functionals on $\Sym(n)$ defined as follows:
 \[ \theta_i(A)=\langle A e_i,e_i\rangle, \quad \quad \alpha_{ij}=\Real \langle Ae_i, e_j \rangle,\quad \quad \beta_{ij}= \Imag \langle Ae_i, e_j \rangle
 \]
 The vector space $\Sym(n)$ is oriented  by the following element of $\bigwedge^{n^2}(\Sym(n))^*$.
  \[ \theta_1\wedge \ldots \wedge \theta_n\wedge\bigwedge_{1\leq i<j\leq n}( \alpha_{ij}\wedge\beta_{ij})
 \] 
 In \cite{N1} (see Example 5.5) L. Nicolaescu introduced coorientations for all the Schubert cells $Z_I(n)$ of $U(n)$ using the basis $\{e_1,\ldots ,e_n\}$.  Roughly the idea is the following.  First, one identifies  $T_{\id}U(n)$ with $T_{U_I}U(n)$ via left multiplication by $U_I$ where
 \[ U_I:= \left\{\begin{array} {ccc} \id & \mbox{ on } &\langle \{e_j~|~j\in I^c\}\rangle\\  -\id &\mbox{ on } &\langle \{e_i~|~i\in I\}\rangle  \end{array}\right.
 \]
 Second, one uses the fact that in the Arnold chart centered at $U_I$ the equations for the Schubert cell $Z_I(n)$ are linear. The equations describing $Z_I(n)$ in this chart are exactly the ones given by Lemma \ref{Schubert cell description}.  The coorientation is induced by the linear exterior form  of rank $N_I$ on $T_{\id}U(n)$ (thought as the space of self-adjoint matrices).
 \[ \bigwedge_{i\in I} \theta_{i} \wedge \bigwedge_{k<i, i\in I}( \alpha_{ki}\wedge\beta_{ki})
 \]
 "transported" via left multiplication by $U_I$ to $T_{U_I}U(n)$. In other words the differential at $0$ of the composition
 \[ \xymatrix {\Sym(n) \ar^{i\cdot}[r] & \mathfrak{u}(n)\ar^{\exp}[r] & U(n)\ar^{U_I\cdot}[r] & U(n)} 
 \]
 takes the previous form to a coorientation form of $Z_I(n)$ at $U_I$.
 
 The coorientation of the corresponding Schubert cells $Z_I(n)$ in the Lagrangian Grassmannian $U(n)$ are induced via the Cayley graph diffeomorphism.
 
In this section we show that the coorientations on $Z_I(n)$ induce coorientations on  the Schubert cells $Z_I\in\Lag^-$. On the other hand Section \ref{Gen Red} we show that the Schubert cells $Z_{k}$ have a natural coorientation coming from the an explicit description of the normal bundle. The two coorientations on $Z_{k}$ are in fact one and the same. \qed\\

\medskip

\begin{prop}\label{finite cycles} The closed set $\overline{Z}_I(n)$ with the  canonical orientation $\omega_I$  is an analytical cycle and so it defines a homology class in $H_{n^2-N_I}(U(n))$ which is Poincare dual to $x_I(n)$.
\end{prop}
\noindent
\textbf{Proof:} See \cite{N1},  Theorem 6.1 \qed\\

It is clear that an orientation on $\mathscr{Z}_I(n)^\circ$ induces a coorientation on the same space by using the normal bundle first convention.

In what follows we will make precise the idea that the \textit{coorientation} of $Z_I(n)$ does not depend on $n$. 

 Let $N>n>\max{\{i~|~i\in I\}}$ and consider the natural inclusion of complex vector spaces
  $\mathbb{C}^n\oplus\mathbb{C}^n\hookrightarrow\mathbb{C}^{N}\oplus\mathbb{C}^N$. This induces an inclusion of lagrangians
  \[ \Lag(n)\hookrightarrow\Lag(N)  , \quad \quad L \rightarrow L+ \mathbb{C}^{N-n}\oplus 0
\] 
Notice that the following diagram where the vertical maps are given  by the Cayley graph diffeomorphisms commutes:
\[ \xymatrix { U(n)\ar@{^{(}->}^{\oplus 1}[r] \ar^{\mathscr{C}}[d] & U(N)  \ar^{\mathscr{C}}[d]  \\
 \Lag(n) \ar@{^{(}->}[r]  & \Lag(N)
}
\]
 
 We will denote by $V:=0\oplus\mathbb{C}^{N-n}\subset 0\oplus \mathbb{C}^N$  the  vector subspace of the vertical component  generated by the last $N-n$ elements of the canonical basis. The varieties $Z_I(N)^\circ$ are defined with respect to the canonical flag of $0\oplus\mathbb{C}^N$. Notice that every lagrangian in $Z_I(N)^\circ$ is clean with $V$ and we can look at the symplectic reduction with $V$.  
 
 The following is obvious
 \begin{lemma} The normal bundle of $Z_I(N)^\circ$ is canonically isomorphic with the pull-back of  the normal bundle of $Z_I(n)^\circ$ via the reduction map. 
 \[\mathscr{R}:=\mathscr{R}_{(n,N)}:\Lag^{V}(N)\rightarrow\Lag(n)
 \]
 \end{lemma}
 \textbf{Proof:}  By Lemma  \ref{stronger},  the restriction of $\mathscr{R}$ to $Z_I(N)^\circ$ is a vector bundle over  $Z_I(n)^\circ$.
 
It is clear that the pull-back of any vector bundle with an orientation is a vector bundle with an orientation. 
\begin{lemma}\label{orientation preserving} The canonical isomorphism between the normal bundle of $Z_I(N)^\circ$ and the pull-back of the normal bundle of $Z_I(n)^\circ$ is orientation preserving.
\end{lemma}
\noindent
\textbf{Proof:} Let us notice first that $\mathscr{R}^*(x_I(n))=x_I(N)$. Indeed due to the fact that $\mathscr{R}$ is a vector bundle over $\Lag(n)$ with zero section $i_n:\Lag(n)\rightarrow\Lag^V(N)$ we have that
\[ i_n^*\circ\mathscr{R}^*(x_I(n))=x_I(n)=i_n^*(x_I(N))
\]
Now the coorientation $\omega_N$ on $Z_I(N)^\circ$ is chosen so that
\[ [\overline{Z}_I(N),\omega_N]= x_I(N)=\mathscr{R}^*(x_I(n))
\]
On the other hand the reduction $\mathscr{R}$ is transversal to $Z_I(n)^\circ$ and so by Proposition \ref{transversality Banach} the preimage of $Z_I(n)^\circ$, which is $Z_I(N)^\circ$, together with the induced coorientation define a cohomology class equal to $\mathscr{R}^*([\overline{Z}_I(n),\omega_I(n)])=\mathscr{R}^*(x_I(n))$.
\qed\\

We are now ready for\\
\noindent
\textbf{Proof of Th. \ref{Schubert and canonical}:} Notice that for $n\geq \max\{i\in I\}$ we have $Z^\circ_I\subset\Lag^{W_n}$.   By Lemma \ref{stronger} we have that $\mathscr{R}_n^{-1}(Z_I^\circ(n))=Z_I^\circ$. So the normal bundle to $Z_I^\circ$ is canonically isomorphic with the pull-back of the normal bundle of $Z_I(n)^\circ$. Hence it induces an orientation. By Lemma \ref{orientation preserving} different  reductions induce the same orientation $\omega_I$ on the normal bundle of $Z_I^\circ$. 

This way $\overline{Z}_I$ gets the structure of a cooriented quasi-submanifold of codimension $N_I$ hence it defines a cohomology class $[\overline{Z}_I,\omega_I]\in H^{N_I}(\Lag^-)$.
Notice that for $n$ big enough the quasi-submanifold $\overline{Z}_{n+1}=\Lag^-\setminus\Lag^{W_n}$  has codimension bigger than $N_I+1$ and so (see Appendix \ref{appendix B})
\[ H^{N_I}(Lag^-)\simeq H^{N_I}(Lag^{W_n})\]  
Moreover the class defined by the cooriented quasi-submanifold $\overline{Z}_I$ in $\Lag^-$ pulls back to give the class defined by the cooriented quasi-submanifold $\overline{Z}_I\setminus \overline{Z}_{n+1}$ in $\Lag^{W_n}$. We denote this later class in $\Lag^{W_n}$ also by $[\overline{Z}_I,\omega_I]$.

Since $\mathscr{R}_n$ is a vector bundle an application of Proposition \ref{transversality Banach} says that 
\[ [\overline{Z}_I,\omega_I]= \mathscr{R}_n^*([\overline{Z}_I(n),\omega_I(n)])= \mathscr{R}_n^*(x_I(n))
\]
On the other hand, $i_n^*\mathscr{R}_n^*(x_I(n))=x_I(n)$ where $i_n:\Lag(H_{w_n})\rightarrow\Lag^{W_n}$ is the natural inclusion and of course the same property stays true if we consider the inclusion $\Lag(H_{w_n})\hookrightarrow\Lag^-$, provided we replace  $\mathscr{R}_n^*(x_I(n))$ by its extension. 
But this is saying that the class $[\overline{Z}_I,\omega_I]$ pulls back via the natural inclusions to the  classes $x_I(n)$. This is exactly the property that characterizes $z_I$. 
 \qed\\

\begin{definition} The  triple composed of the Schubert variety $\overline{Z}_I$  with the standard filtration and the coorientation $\omega_I$ is called the Schubert cocycle  or the  geometric representative of $z_I$.  The cohomology class it represents is denoted by
\[ [\overline{Z}_I,\omega_I]
\]
\end{definition}

We consider now families of vertical, Fredholm lagrangians. By that we simply mean smooth maps
\[ F:M\rightarrow \Lag^-
\]
\begin{definition} A map $F:M\rightarrow\Lag^-$ is said to be (standard) transversal to $\overline{Z}_I$ if it is transversal to every stratum in the standard stratification.
\end{definition}

\begin{lemma} Any smooth family $F:M\rightarrow \Lag^-$ can be deformed by a smooth homotopy to a family  transversal  to $\overline{Z}_{k}$. 
\end{lemma}
\noindent
\textbf{Proof:} Since $M$ is compact, transversality with $\overline{Z}_I$ means actually transversality of the reduced family with $\overline{Z}_I$ for $n$ big enough. Transversality with Whitney stratified spaces is an open, dense condition in the space of all smooth maps $G:M\rightarrow \Lag (n)$.
\qed\\
\begin{prop}\label{pull back and preimage} Let  $M$ be a closed oriented manifold  and let $F:M\rightarrow\Lag^-$ be a family transversal to $Z_I$. Then $F^{-1}(Z_I)$ is quasi-submanifold of $M$ with a naturally induced coorientation $F^*\omega_I$ and 
\[ [F^{-1}(Z_I), F^*\omega_I]= F^*[\overline{Z}_I,\omega_I]
\]
\end{prop}
\noindent
\textbf{Proof:} The pull-back of the normal bundle to $Z_I^\circ$ is naturally isomorphic with the normal bundle to $F^{-1}(Z_I^\circ)$ and the coorientation $F^*\omega_I$ is the one induced via this isomorphism. For the rest, see Proposition \ref{transversality Banach}.\qed\\

In the infinite dimensional context, Poincar\'e Duality does not make sense. Instead we aim for an expression of Poincar\'e duality for families of lagrangians parametrized by a closed, oriented manifold $M$. One way to build homology classes out of stratified spaces is via   the theory of analytic cycles, which we already mentioned, used by Nicolaescu to prove the duality $Z_I(n)=x_I(n)$.

Another way, which is more appropriate to the point of view we take in this paper is via Borel-Moore homology. In the Appendix \ref{appendix B}  we describe the relevant aspects. 

 Inside  an oriented manifold $M$ of dimension $n$, any oriented  quasi-submanifold $\eF$ of dimension $d$   defines a  Borel-Moore \index{Borel-Moore homology} homology class as follows. Every smooth, oriented manifold $S$ of dimension $d$ has an orientation class $[S]\in H^{BM}_d(S)$. In the case  of an oriented  quasi-submanifold $\eF$ of dimension $d$,  $S:=\eF\setminus\eF^2$ represents the top stratum. 
 This class can be extended to a class in $H^{BM}_{d}(\eF)$ because the absence of  singularities in codimension one implies that we have an isomorphism $H^{BM}_{d}(\eF)\simeq H^{BM}_d(S)$. Finally, this class can be pushed-forward to a class in the ambient space $M$. It turns out that when $M$ is closed and oriented this class is Poincar\'e dual to the cohomology class determined by $\eF$ with the coorientation induced in the obvious way.

In the case when $\eF$ is compact, e.g. when $M$ is compact,   then the Borel-Moore homology group of $\eF$ coincides with the singular homology group. 

We summarize our discussion:

\begin{theorem}\label{intersection} Let $F:M\rightarrow \Lag^-$ be a smooth map from an oriented, closed manifold $M$ of dimension $n$ to $\Lag^-$. Suppose $F$ is transversal to $\overline{Z}_I$. Then the preimage $F^{-1}(\overline{Z}_I)$   has a naturally  induced orientation and so it defines a homology class $ [F^{-1}(Z_I)]_M\in H_{n-N_I}(M)$ which is Poincar\'e dual to the class $F^*[Z_I]$. 
\end{theorem}

\begin{remark} The fundamental class of an oriented quasi-submanifold can be defined without  appeal to Borel-Moore homology, provided something stronger is true. Suppose  $\eF=\eF^0\supset \eF^1=\eF^2\supset\ldots \supset \eF^k$ is an oriented quasi-submanifold of dimension $d$ such that every pair $(\eF^i, \eF^{i+1})$ is a good pair, i.e. there is an open neighborhood $U\subset \eF^{i}$ of $\eF^{i+1}$ that retracts to $\eF^{i}$. For example, if  the stratification satisfies the Whitney condition then Goresky has shown that $\eF$ can be triangulated in such a way that  the triangulation respects the filtration (see \cite {Gor1}, Prop.5). Let $U$ be a neighborhood of $\eF^1$ in $\eF$ that retracts to $\eF^1$. Then the  Poincar\'e duality for a manifold with boundary says that
\[ H^i(\eF\setminus U)\simeq H_{d-i}(\eF \setminus U, \partial (\eF\setminus U))
\] 

This implies that $H^{i}(\eF\setminus \eF^1)\simeq H_{d-i}(\eF ,\eF^1)$ (both are isomorphic with the  $(d-i)$-th Borel-Moore homology group of $\eF\setminus \eF^1$). Therefore the canonical class $1$ from $H^0(\eF\setminus \eF^1)$  gives a class $a\in H_d(\eF,\eF^1)$. The map
\[ H_d(\eF)\rightarrow H_d(\eF,\eF^1)
\]
is an isomorphism because there are no singularities in codimension $1$, i.e. $\eF^1=\eF^2$. For the proof  one uses the fact that $(\eF^i, \eF^{i+1})$ is a good pair. Then one can "extend" the  class $a$ to a class in $H_d(\eF)$ which is the fundamental class of the quasi-submanifold. 

The following considerations justify the fact that our quasi-submanifolds fit into the picture just described.
For a map $F:M\rightarrow \Lag^-$ transversal to the quasi-manifold $\overline{Z}_I$ the preimage $F^{-1} (\overline{Z}_I)$ is always a Whitney stratified space. To see why it is Whitney stratified notice that since $M$ is compact then $F(M)\subset\Lag^{W_n}$ for $n$ big enough and it is easy to see that $F$ is transversal to $\overline{Z}_I$ implies $F$ transversal to $\overline{Z}_I(n)$ and $F^{-1} (\overline{Z}_I)=F^{-1}(\overline{Z_I(n)})$. Now $Z_I(n)$ is a Whitney stratified space and Whitney property is preserved under transversal pull-backs. \qed
\end{remark}

\bigskip

\subsection{Generalized Reduction}\label {Gen Red}

In this section we take the first steps towards doing intersection theory. Recall that the standard stratification of $\overline{Z}_I$ has as its top stratum the set $Z_I^\circ$. It turns out that at least in the case when $\sharp I=1$ there is a better stratification of $\overline{Z}_k$, where $k\in\mathbb{N}^*$, which comes with a natural coorientation and is more suitable for intersection theory. Here we define this stratification and describe the normal bundle of the maximum stratum. This top stratum of the new stratification contains $Z_I^\circ$ and it also defines a quasi-submanifold structure for $Z_I$.  By Remark \ref{diff strat} each of the stratifications defines the same cohomology class.

 In section $\ref {Symp red}$ we described the process of symplectic reduction with cofinite isotropic space as a differentiable map going from the set of clean lagrangians to a finite dimensional Lagrangian Grassmannian. The symplectic reduction is a well defined process on the entire $\Lag^-$, the trouble being that it is not a continuous map everywhere. However it is continuous and in fact differentiable on certain submanifolds of $\Lag^-$.  First a definition. 
\begin{definition}
For a fixed, finite codimensional of codimension $p$ subspace $W\subset H^-$  let 
\[ \Lag^W(k):=\{L\in \Lag^- ~|~ \dim{L\cap W}=k\}
\]
be the space of lagrangians that intersect $W$ along a space having fixed dimension $k$. We call these lagrangians \textit{\textbf{$k$-clean}} \index{$k$-clean} or just $\textit{\textbf{clean}} $ when $k=0$.
\end{definition}

\noindent
\begin{remark} For a complete, decreasing flag
\[ H^-:=W_0\supset W_1\supset W_2\supset \] 
we have
 \[
  Z_{k}^\circ\subset \Lag^{W_{k-1}}(1)\subset \overline{Z_{k}}
 \]
 Compare with Corollary \ref{fund Sch}.
\end{remark}
In the rest of this section we will prove that $ \Lag^W(k)$ is a coorientable submanifold of $\Lag^-$ and  we will identify the normal bundle of $ \Lag^W(k)$ with a certain tautologically defined bundle over $\Lag^W(k)$. 

\begin{lemma} For every $k$-clean lagrangian $L$ let $V=L\cap W$ and $V^\perp$ be its orthogonal in $W$. Then $L$ is clean with the isotropic space $W_L:=JV\oplus V^\perp$. 

Moreover, the symplectic reduction with $W$ coincides with the symplectic reduction with $W_L$, that is if  
\[
\mathscr{R}_{W}(L)=\Range{ P_{H_W}|_{L\cap W^\omega}}
\]
then
\[ \mathscr{R}_{W}(L)=\mathscr{R}_{W_L}(L)
\]
Therefore it is a well-defined map $\mathscr{R}:\Lag^W(k)\rightarrow\Lag {(H_W)}$.
\end{lemma}
\noindent
\textbf{Proof:} The first claim is obvious. For the second notice that $H_{W_L}=H_W$.

What we have to compare are the projections of 
    \[ L\cap (J(L\cap W)\oplus V^\perp\oplus H_W) \mbox{ and } L\cap ((L\cap W)\oplus V ^\perp\oplus H_W)
    \] 
    onto $H_W$.  Let us notice that $L\cap (J(L\cap W)\oplus V^\perp\oplus H_W) =L\cap(V^\perp\oplus H_W)$. Indeed if one writes 
    \[ 
    \begin{array}{ccccccc}
    x &=& a &+&b&+& c\\
    L &&J(L\cap W) &&V^\perp &&H_W 
    \end{array}
    \]
    notice that $a\in JL\perp L$ and $a\perp b,c$ so $a=0$. 
    
    Now $P_{H_W}(L\cap(V^\perp\oplus H_W))\subset P_{H_W}(L\cap ((L\cap W)\oplus V ^\perp\oplus H_W))$ obviously  and the other inclusion follows from noticing that 
    if 
    \[\begin{array}{ccccccc}
    x &=& a &+&b&+& c\\
    L &&(L\cap W) &&V^\perp &&H_W 
    \end{array}
    \]
    then $x-a\in L\cap (V^\perp\oplus H_W) $ and $P_{H_W}(x)=P_{H_W}(x-a)=c$. 
    \qed \\
    
  \begin{corollary}  \label {dimension}  For every lagrangian $L\in\Lag^W(k)$ the intersection $L\cap W^\omega$ has dimension equal to $k+p=\dim{L\cap W} +1/2\dim{H_W}$. Moreover, using the same notations as in the lemma,  $L\cap W^\omega$ decomposes orthogonally as
 \[
 L\cap W^\omega=L\cap W\oplus L\cap (V^\perp\oplus H_W)
 \]

  \end{corollary}
  \noindent
  \textbf{Proof:} The  image of the projection $P_{H_W}:L\cap W^\omega\rightarrow H_W$ has dimension equal to $1/2\dim{H_W}$ and the kernel is just $L\cap W$. 
  
  Clearly the two spaces that appear in the sum are orthogonal and they are both subsets of $L\cap W^\omega$. So it is enough to prove that $L\cap V^\perp\oplus H_W$ has dimension $p$. But we saw in the proof of the lemma that
  \[ L\cap (V^\perp\oplus H_W)=L\cap J(L\cap W)\oplus V^\perp\oplus H_W=:L\cap W_L^\omega
  \]
  
  Moreover the projection $P_{H_W}:L\cap W_L^\omega\rightarrow H_W$ is injective and its image is a lagrangian in $H_W$ which has dimension $p$.
  \qed\\
   
   The notations we used in the previous lemma will be used throughout this section. 
\begin{definition} For every lagrangian $L\in\Lag^W(k)$ let $V:=L\cap W$, $V^\perp$ be the orthogonal complement of $V$ in $W$ and let $\ell$ be the symplectic reduction of $L$ with $W$. The space $L_W:=\ell\oplus V\oplus JV^\perp$ is called the \textit{\textbf{associated lagrangian}}\index{associated lagrangian} or simply the \textit{\textbf{associate}}.
\end{definition}   
   
 \begin{lemma} \label {associate}
 \begin{itemize} \item[a)] The associated lagrangian is in $\Lag^W(k)$ and every lagrangian is in the Arnold chart of its associate. Moreover $L$ is given by the graph of $JX$ where $X\in\Sym(L_W)$ has the block decomposition:
 \[
\left(\begin{array}{ccc}
0 & 0 & X_2^*\\
0 & 0  & 0 \\
X_2 &0 &  X_4
\end{array}
\right) 
\]

\item[b)]  Let $W=V\oplus V^\perp$ be an orthogonal decomposition of $W$ such that $V$ is k-dimensional and let $\ell\subset  H_W$ be a lagrangian. Then $\ell\oplus V\oplus JV^\perp$ is in $\Lag^W(k)$ and the $\Lag^W(k)$  is described in the Arnold chart  $\mathcal{A}_{\ell\oplus V\oplus JV^\perp}$ by linear equations. More precisely given $S\in\Sym(\ell\oplus V\oplus JV^\perp)$ then $\Gamma_{JS}\in\Lag^W(k)$ if and only if its $V\times V$ and $V\times \ell$ blocks are zero.

\item[c)] The space $\Lag^W(k)$ is a submanifold of $\Lag^-$ of codimension $k^2+2pk$ and the symplectic reduction map:
\[
\mathscr{R}:\Lag^W(k)\rightarrow \Lag(H_W), \quad \quad L\rightarrow~\Range{ P_{H_W}|_{L\cap W^\omega}}\]
is differentiable.
\end{itemize}
\end{lemma} 
\noindent
\textbf{Proof:} a) The fact that the associated lagrangian is indeed a lagrangian is a simple check. Now 
$L_W\cap W=V$, hence clearly $L_W$ is in $\Lag^W(k)$.

For the second claim, notice that $(L, V^\perp)$ is a Fredholm pair and $V^\perp$ and $L_W^\perp$ ar commensurable, so $(L,L_W^\perp)$ is a Fredholm pair.  Moreover the intersection $L\cap L_W^\perp$ is trivial. Indeed let
\[ \begin {array} {ccccccc}
x&= &a&+&b&+&c\\
L & &J\ell & &JV&&V^\perp
\end{array}
\]
Then  $b\in L^\perp$ and so $b=0$. From $x=a+c$ it follows that $x\in L\cap W^\omega$ and $a=P_{H_W}(x)\in \ell$ so $a=0$. This implies $x=c\in L\cap V^\perp=\{0\}$.

For the last part notice that if $L$ is the graph of an operator $JS:L_W\rightarrow L_W^\perp$ then $JS|_{L\cap L_W}=0$. Simply because $JS=P_{L_{W^\perp}}\circ (P_{L_W}|_L)^{-1}$. On the other hand $V\subset L\cap L_W$. This and the self-adjointness implies that the middle row and column of $S$ are zero. The vanishing of the top, left block follows from the following considerations. The symplectic reduction of any lagrangian in the Arnold chart of $L_W$ with $W_L:=JV\oplus V^\perp$ is just the graph of the $\ell\times \ell$ block of the self-adjoint operator on $L_W$ that gives $L$. But the only operator $\ell\rightarrow \ell$ for which $\ell:=\mathscr{R}_{W_L}(L)=\mathscr{R}_{W}(L)$ is the graph of, is the zero operator.\\

b)  The first claim, that $\ell\oplus V\oplus JV^\perp\in \Lag^W(k)$ is trivial.

 Now, every lagrangian in the Arnold chart $\mathcal{A}_{\ell\oplus V \oplus JV^\perp}$ is just the graph of an operator $JS$ where   $S\in \Sym{(\ell\oplus V \oplus JV^\perp)}$. So $S$ has a block decomposition
\[
  S=\left(\begin{array}{ccc}
  S_{\ell,\ell} & S_{V,\ell} &S_{JV^\perp,\ell} \\
  S_{\ell, V} & S_{V, V} & S_{JV^\perp ,V} \\
  S_{\ell, JV^\perp} & S_{V, JV^\perp} & S_{JV^\perp, JV^\perp}
  \end{array}
  \right)
  \]

The condition $v+JSv\in W$ where $v=(v_1,v_2,v_3)\in \ell\oplus V\oplus JV^\perp$ translates to the sum
\[
  \begin{array}{ccc}
   v_1\quad +\quad v_2\quad+\quad v_3&+& (JS_{\ell,\ell}v_1+JS_{V,\ell}v_2+JS_{JV^\perp,\ell}v_3)\quad+ \\
  \quad \ell ~~\quad\quad \quad V \quad \quad \quad JV^\perp& &J\ell \\
   (JS_{\ell, V}v_1+JS_{V,V} v_2+JS_{JV^\perp, V}v_3) &+&(JS_{\ell, JV^\perp}v_1+J S_{V, JV^\perp} v_2+JS_{JV^\perp, JV^\perp}v_3) \\
 JV& & V^\perp
   \end{array}
  \]
being in $V\oplus V^\perp$. Since $\ell\oplus J\ell\oplus V\oplus JV^\perp\perp V\oplus V^\perp $ we get 
\[v_1=v_3= (JS_{\ell,\ell}v_1+JS_{V,\ell}v_2+JS_{JV^\perp,\ell}v_3)=(JS_{\ell, V}v_1+JS_{V,V} v_2+JS_{JV^\perp, V}v_3) =0 \mbox{ and }
\]
\[v_2+(JS_{\ell, JV^\perp}v_1+J S_{V, JV^\perp} v_2+JS_{JV^\perp, JV^\perp}v_3)\in V\oplus V^\perp 
\]

We conclude that in order for $v+JSv$ to be in $W$ one must have $v_2\in \Ker{T}$ where \linebreak $T:=(S_{V,\ell},S_{V, V} ):V \rightarrow \ell\oplus V$. Also $\Gamma_{JS}\cap W$ is the graph of the restriction ${(JS_{V, JV^\perp} |_{\Ker{T}})}$. 
The only way the graph of  $JS_{V, JV^\perp} |_{\Ker{T}}$ can have dimension equal to the dimension of $V$ is if $\Ker{T}=V$, that is
   \[
   (S_{V,\ell},S_{V, V} )=0
   \]

Hence the intersection $\mathcal{A}_{\ell\oplus V\oplus JV^\perp}\cap\Lag^W(k)$ consists of graphs of operators whose $V\times V$, $V\times \ell$ and $\ell\times V$ blocks are zero.\\

c) Every  lagrangian $L\in \Lag^W(k)$ is in the Arnold chart of its associate which is of the type required by part $b)$. In these charts  $\Lag^W(k)$ is described by linear equations and one can very fast see that the codimension is the one indicated.

In the Arnold chart of $\ell\oplus V\oplus JV^\perp$ the symplectic reduction of any lagrangian $L_0\in\Lag^W(k)$ with $W$ is the graph of the projection onto the $\ell\times\ell$ block and the differentiability follows.
\qed\\

  We would like to say something about the diffeomorphism type of $\Lag^W(k)$. For that end let us notice that beside the symplectic reduction, $\mathscr{R}_W$, there is another natural map one can define on $\Lag^W(k)$ namely:
    
    \[ \mathscr{D}: \Lag^W(k)\rightarrow \Gr(k,W), \quad \quad	\mathscr{D}(L)=L\cap W
    \]
    where $\Gr(k,W)$ is the grassmannian of $k$-dimensional subspaces of $W$. Before we proceed to study this map let us recall a few well-known facts about the infinite grassmannian.
    
     For every Hilbert space $H$ the set $\Gr(k,H) $ gets the structure of a metric space by considering each subspace being represented by the corresponding projection and considering the norm topology on the set of all these projections. It gets the structure of (complex) Banach manifold by the following simple lemma:
\begin{lemma} For every $k$-dimensional subspace $V_0\subset H^-$ the map
\[ \Hom(V_0, V_0^{\perp})\rightarrow \Gr(k,H),  \quad \quad T\rightarrow\Gamma_T
\]  sets a homeomorphism between $\Hom(V_0, V_0^{\perp})$ and  $\{V\in\Gr(k,H)~|~V\cap V_0^\perp=\{0\}\}$ which is an open subset of $\Gr(k,H) $. \end{lemma}
\noindent
\textbf{Proof:} When $V\cap V_0^\perp=\{0\}$ the orthogonal projection $P_{V_0}|_V:V\rightarrow V_0$ is a linear isomorphism so $V$ is just the graph of a linear map from $V_0$ to $V_0^\perp$. Now the condition $V\cap V_0^\perp=\{0\}$ implies $V_0\cap V^\perp=\{0\}$. Otherwise, due to the dimension constraints $V_0+V^\perp$ would be a proper subset of $H$ which, by taking orthogonal complements would say that $V\cap V_0^\perp\neq\{0\}$. But $V\cap V_0^\perp\oplus V_0\cap V^\perp=\{0\}$ is equivalent with $P_{V}-P_{V_0^\perp }$ is invertible and that proves that $\{V\in\Gr(k,H)~|~V\cap V_0^\perp=\{0\}\}$ is an open set. The continuity of the maps is immediate since the $P_{\Gamma_T }$ can be computed explicitly in terms of $T$ and $T^*$.\
\qed \\

\begin{corollary} The tangent space of $\Gr(k,H)$ is naturally isomorphic with the homomorphism bundle associated with the tautological bundle and its dual: $T\Gr(k,H^-)\simeq \Hom(\tau, \tau^\perp)$ where $\tau:=\{(V,v)~|~v\in V\}\subset \Gr(k,H)\times H$.
\end{corollary}
We will denote the open sets $\{V\in\Gr(k,H)~|~V\cap V_0^\perp=\{0\}\}$ by $\mathcal{A}_{V_0}^{gr}$.

\begin{definition} The \textit{\textbf{generalized reduction}}\index{generalized reduction} is the map:
\[
\mathscr{R}:\Lag^W(k)\rightarrow \Lag(H_W)\times \Gr(k,W) 
\]
\[L\rightarrow (\Range P_{H_W}(L\cap W^\omega), L\cap W) 
\]
\end{definition}
\noindent
\begin{remark} Notice that in the case $k=0$ we get what we called symplectic reduction since the second component is just a point. This is why we prefer to keep the notation $\mathscr{R}$. In fact to eliminate any possibility of confusion we will denote symplectic reduction from now on by $\mathscr{R}^1$ since that is the first component in our generalized reduction.
\end{remark}

The generalized reduction behaves very much like the symplectic reduction meaning it inherits the structure of a vector bundle whose fiber we will identify in a moment.
First let us see that $\mathscr{R}$ comes with a natural section namely
\[ \mathscr{S}: \Lag(H_W)\times \Gr(k,W) \rightarrow \Lag^W(k),\quad \quad  (\ell,V)\rightarrow \ell\oplus V\oplus JV^\perp
\]

Every associate lagrangian lies on this section.

\begin{theorem} \label{tangent gen red} \begin{itemize} 
\item[a)] The tangent space of $ \Lag^W(k)$ along the section defined by $\mathscr{S}$ can  be naturally identified with the vector subbundle of $T\Lag^-$ whose fiber at $\ell\oplus V\oplus JV^\perp$ consists of self-adjoint operators $S\in \Sym(\ell\oplus V\oplus JV^\perp)$ which have the following block decomposition:
\[S= \left( \begin{array}{ccc}
S_1 &0 &S_2^*\\
0 &0 &S_3 ^* \\
S_2 & S_3 &S_4
\end{array}
\right)
\]

\item [b)] The generalized symplectic reduction is differentiable  and $(\Ker~d\mathscr{R})|_{\Lag(H_W)\times \Gr(k,W)}$ can be identified with the vector subbundle of $T\Lag^W(k)|_{ \Lag(H_W)\times \Gr(k,W) }$ whose fiber at $\ell\oplus V\oplus JV^\perp$ consists of self-adjoint operators $S\in \Sym(\ell\oplus V\oplus JV^\perp)$ which have the following block decomposition:
\[
S= \left( \begin{array}{ccc}
0 &0 &S_2^*\\
0 &0 & 0 \\
S_2 & 0 &S_4
\end{array}
\right)
\]
\item[c)] The natural map
 \[
  \mathscr{N}:(\Ker~d\mathscr{R})|_{\Lag(H_W)\times \Gr(k,W)}\rightarrow \Lag^W(k), \quad \quad
 \mathscr{N} (\ell\oplus V\oplus JV^\perp, S)=\Gamma_{JS}
 \]
 is a diffeomorphism that makes the diagram commutative:
 
 \[ \xymatrix{ \ar @{} [dr] |{}
(\Ker{d\mathscr{R}})\bigr |_{\Lag(H_W)\times \Gr(k,W)} \ar[rr]^{\quad\mathscr{N}} \ar[dr] && \Lag^W(k) \ar[ld]_{\mathscr{R}} \\
 & \Lag(H_W) \times \Gr(k,W) &                    }
\]
 
\item [d)] The space $\Lag^W(k)$ is diffeomorphic with $\Lag(H_W)\times (\tau^\perp)^p\oplus \Sym(\tau^\perp)$ where $\tau^\perp$ is the orthogonal complement of the tautological line bundle over $\Gr(k,W)$.
\end{itemize}
\end{theorem}
\noindent
\textbf{Proof:} a) This is  obvious since as we saw in the proof of the Lemma $\ref {associate}$,  in the charts centered at $L=\ell\oplus V\oplus JV^\perp$ the manifold $\Lag^W(k)$ can be described exactly as the set of those self-adjoint operators with the claimed block decomposition. \\

b) In what concerns the differentiability, we only have to prove that the second component, $\mathscr{R}^2$, is differentiable. For that we again send to the proof of Lemma $\ref {associate}$ where we saw that in the Arnold chart $\mathcal{A}_{\ell\oplus V \oplus JV^\perp}$ the intersection $\Gamma_{JS}\cap W$ is just $\Gamma_{JS_{V,JV^\perp}}$ for every $S\in \Lag^W(k)$.  

The second claim is also obvious when one works in the Arnold charts centered at $L=\ell\oplus V \oplus JV^\perp$ since then $d_{L}\mathscr{R}$ is just the projection on the $\ell\times\ell$ and $V\times JV^\perp$ blocks.\\

c) We construct an inverse for $\mathscr{N}$. To every $L\in \Lag^W(k)$ we associate the lagrangian $\ell\oplus V\oplus JV^\perp$ where $V=L\cap W$ and $\ell$ is the symplectic reduction with $W$. In the Arnold chart centered at $\ell\oplus V\oplus JV^\perp$ the lagrangian $L$ is a graph $\Gamma_{JS}$ where $S$ has to be of the type:
\[
S= \left( \begin{array}{ccc}
0 &0 &S_2^*\\
0 &0 & 0 \\
S_2 & 0 &S_4
\end{array}
\right)
\]
This is because in the lemma $\ref {associate}$, $V=L\cap W=\Gamma_{JS_{V,JV^\perp}}$ and so $JS_{V,JV^\perp}:V\rightarrow V^\perp$ has to be zero.  So the inverse to $\mathscr{N}$ associates to $L$ the lagrangian $\ell\oplus V\oplus JV^\perp$ and the two operators, $S_2$ and $S_4$, which are the projections onto the $(\ell,JV^\perp)$ and $(JV^\perp,JV^\perp)$.

d) By b) and c) we have identified the fiber of $\Lag^W(k)$ over $\Lag_{H_W}\times \Gr(k, W)$ at $(\ell, V)$ to be the vector space $\Hom(\ell,JV^\perp)\oplus \Sym{(JV^\perp)}$. We know that the tautological bundle over $\Lag^W(k)$ is  naturally trivializable so the bundle with fiber $\Hom(\ell,JV^\perp)$ over  $\Lag_{H_W}\times \Gr(k, W)$  is naturally isomorphic with the bundle $\Hom (W^\perp, JV^\perp)$ where the lagrangian $W\in \Lag(H_W)$ is just the orthogonal complement of $W$ in $H^-$.
 A choice of a basis on $W^\perp$ proves that $\Hom (W^\perp, JV^\perp)$ is in fact $(\tau^\perp)^p$. \qed\\
 
 \noindent
\begin{example}
  We will describe the spaces $\Lag^{W}(k)$ when $\dim{H}=2$. The spaces $\Lag^W(0)$ are open subsets of $\Lag^-$ and were described in section $\ref {Symp red}$. 
  
    Let $H=\mathbb{C}^2$ and $W_1=\langle e_2\rangle$. Then
  \[ \Lag^{H^-}(1)=Z_{\{1\}}\cup Z_{\{2\}} 
  \]
 is a $3$-dimensional sphere minus a point.  The point is 
  \[\Lag^{H^-}(2)=Z_{\{1,2\}}
  \] 
  In terms of unitary operators these spaces correspond to sets of operators for which $\Ker {1+U}$ has dimension $1$ or $2$.
  
  Another non-trivial space which has the diffeomorphism type of the circle is 
  \[ \Lag^{W_1}(1)=Z_{\{2\}}\cup Z_{\{1,2\}}
  \]
  \end{example}

 One of our general goals is to describe how one can do intersection theory in $\Lag^-$. We would therefore need a description of the normal bundle of $\Lag^W(k)$ in $\Lag^-$. It will be enough to describe a splitting of the differential of the inclusion $\Lag^W(k)\hookrightarrow\Lag^-$. We have a canonical choice for this splitting in the charts along the zero section as the Theorem \ref{tangent gen red} shows.   In order to find a global characterization we will use transition charts. The next two lemmata  are very important. 
 \begin{lemma} Let $L,L_0\in\Lag^-$ be two lagrangians such that $L\in\mathcal{A}_{L_0}$. Then the differential at $L$ of the transition map between the Arnold chart  centered at $L_0$ and the Arnold chart centered at $L$  is the map:
 \[ d_L:\Sym{(L_0)}\rightarrow \Sym{(L)} \quad \quad d_L(\dot{S})=P_{L}|_{L_0}\circ \dot{S}\circ P_{L_0}|_{L}
 \]
 \end{lemma}
 \noindent
\textbf{Proof:} Let $L_1$ be in $\mathcal{A}_{L_0}\cap\mathcal{A}_{L}$. This means that $L_1$ can be described both as $\Gamma_{JX}$ where $X\in\Sym{(L_0)}$ and $\Gamma_{JS}$ where $S \in\Sym{(L)} $ It is not hard to see what $S$ should be.
\[ JS=P_{L^\perp}\circ (I,JX)\circ [P_{L}\circ (I,JX)]^{-1}
\]
The image of the map $(I,JX):L_0\rightarrow \hat{H}$ gives the lagrangian $L_1$ and the inverse of $P_{L}\circ (I,JX)$ is a well-defined operator $L\rightarrow L_1$ since $L_1$ is in $\mathcal{A}_{L}$. We consider the function:
\[ F:\Sym{(L_0)}\rightarrow\Sym{(L)}, \quad \quad F(X)=-JP_{L^\perp}\circ (I,JX)\circ [P_{L}\circ (I,JX)]^{-1}
\]
Notice that for $X_0=-JP_{L_0^\perp}\circ (P_{L_0}|_L)^{-1}$ we have $F(X_0)=0$ since $X_0\in\Sym{(L_0)}$ is the self-adjoint operator such that $L=\Gamma_{JX_0}$. We compute the differential of $F$ at $X_0$:
\[ d_{X_0}F(\dot{S})=-JP_{L^\perp}\circ (0,J\dot{S})\circ  [P_{L}\circ (I,JX_0)]^{-1}-JP_{L^\perp}\circ(I,JX_0)\circ[\ldots]=~~~~
\]
\[ =-JP_{L^\perp}\circ (0,J\dot{S})\circ  [P_{L}\circ (I,JX_0)]^{-1}~~~~~~~~~~~~~~~~~~~~~~~~~~~
\]
The reason for the cancellation of the second term is that the image of $(I,JX_0)$ is in $L$.  It is easy to see that $[P_{L}\circ (I,JX_0)]^{-1}=P_{L_0}|_{L}$ the restriction to $L$ of the projection onto $L_0$. Also since $P_{JL}(Jv)=JP_L(v)$ for any lagrangian $L$ and for any $v\in \hat{H}$ we get that $-JP_{L^\perp}\circ (0,J\dot{S})=P_L\circ \dot{S}$. So
\[ d_{X_0}F(\dot{S})=P_{L}|_{L_0}\circ \dot{S}\circ P_{L_0}|_{L}
\]
and this is our $d_L$.
\qed\\

It is convenient to have another description of the differential of the transition map. For that end let us recall that Arnold's theorem provides a canonical unitary isomorphism:
\[ \tilde{U}:L_0\rightarrow L, \quad \quad \tilde{U}(v):=\frac{1}{2}[(1+U)v+iJ(1-U)v] \quad \forall ~v\in L_0
\]
where $U\in \mathscr{U}(L_0)$ is the Cayley transform of the self-adjoint operator $X_0\in\Sym(L_0)$ that gives $L$ as a graph of  $JX_0:L_0\rightarrow L_0^\perp$. Notice first that the projection $P_{L}|_{L_0}$ has a description in terms of the same self-adjoint operator $X_0$. The orthogonal $L^\perp$ is the switched graph of $-(JX_0)^*=X_0J$. So in order to find the projection $P_{L}|_{L_0}$ in terms of $X_0$ one needs to solve the system
\[\left\{\begin{array}{ccccc}
a&=&v &+ &X_0J w \\
0&= &JX_0 v &+& w
\end{array}
\right.
\]
where $a,v\in L_0$ and $w\in L_0^\perp$. That is easy to do and one gets 
\[ v=(1+X_0^2)^{-1}(a)
\]
which delivers the expression for the projection:
\[ P_{L}|_{L_0}(a)=(1+X_0^2)^{-1}(a)+JX_0(1+X_0^2)^{-1}(a)
\]
We now plug in 
\[X_0=i(1+U)^{-1}(1-U)
\] 
to conclude that 
\[ P_{L}|_{L_0}(a)=\frac{1}{2}\tilde{U}((1+U^*)(a))
\] 
Since $P_{L_0}|_{L}=(P_{L}|_{L_0})^*$ we have just proved the following result:
\begin{lemma} \label {transition}
The map $d_L$ in the previous lemma can be written as:
\[ 
d_L(\dot{S})=\frac{1}{4}\tilde{U}(1+U^*)\dot{S}(1+U)\tilde{U}^*
\]
\end{lemma}

\begin{definition} Let $j:E\rightarrow F$ be an injective morphism of vector bundles over a smooth Banach manifold $\mathscr{X}$. An \textit{\textbf{algebraic complement}} $G$ of $E$ is a vector bundle over that splits $j$. This means that there exists an injective morphism $k: G\rightarrow F$ such that
\[ F=E\oplus G
\]
\end{definition}
\noindent
\textbf{Notation:} Let $F:\mathscr{X}_1\rightarrow \mathscr{X}_2$ be a smooth immersion of Banach manifolds.  An algebraic complement of the tangent bundle $T\mathscr{X}_1$  is denoted by $N\mathscr{X}_1$.

\begin{lemma}  Let  $F:\mathscr{X}_1\rightarrow \mathscr{X}_2$ be a smooth immersion of Banach manifolds. Then every algebraic complement of $T\mathscr{X}_1$ is naturally isomorphic with the normal bundle $\nu\mathscr{X}_1$. 
\end{lemma}
\noindent
\textbf{Proof:}  The natural projection $N\mathscr{X}_1\rightarrow \nu\mathscr{X}_1$ is an isomorphism.\qed \\

We have all we need for proving the following
\begin{prop} \label {algebraic complement}
\begin {itemize}
\item [a)] Every lagrangian $L\in\Lag^W(k)$ has an orthogonal decomposition $L=\ell\oplus L\cap W\oplus \Lambda$ where $\ell$ is the orthogonal complement of $L\cap W$ in $L\cap W^\omega$ and $\Lambda$ is the orthogonal complement of $L\cap W^\omega$ in $L$. Then the space of operators $S\in \Sym{(L)}$ with block decomposition:
\[
S= \left( \begin{array}{ccc}
0 & S_1^* & 0\\
 S_1& S_2& 0 \\
 0 & 0 & 0
\end{array}
\right)
\] is an algebraic complement of $T_L\Lag^W(k)$. 
\item[b)] The algebraic complement  of $T\Lag^-|_{\Lag^W(k)}$ described above is an finite dimensional,  orientable bundle. If $k=1$, it has a natural orientation.
\end{itemize}
\end{prop}
\noindent
\textbf{Proof:} a) The claim is clearly true for any associate lagrangian $L_W$ by Lemma $\ref {associate}$. We want to use the transition maps between two different Arnold charts at $L$, namely the one given by $L_W$ and the one centered at $L$ to show that the claim is true in general.

 In order to avoid any confusion we will let $\R_W{L}=:\ell_0\subset L_W$. So $L_W=\ell_0\oplus L\cap W\oplus JV^\perp$. 
 
 By definition $\ell_0=P_{H_W}(\ell)$. We are looking for a relation between $\ell$ and $\ell_0$ in terms of the unitary isomorphism  $U$. Here $U\in\mathscr{U}(L_W)$ is the Cayley transform of the self-adjoint operator $X$ whose graph is $L$. That is:
 \[ X=i\frac {1-U}{1+U}, \quad \quad U=\frac{i-X}{i+X} \quad \mbox { and } \quad L=\Gamma_{JX}
 \]
 
 It is not hard to see from what we just said that $\ell$ is the graph of the restriction $JX|_{\ell_0}$. Now 
 $\tilde{U}:L_W\rightarrow L$ has the following expression.
 \[ \tilde{U}v=\frac{1+U}{2}v+iJ\frac{1-U}{2}v
 \]
 
 In other words 
 \[ 2\tilde{U}(1+U)^{-1}w=w+JSw
 \]
 
We conclude that 
\[2\tilde{U}(1+U)^{-1}\ell_0=\ell.\]

Let $\tilde{\ell}:=\tilde{U}^*\ell$ and $\tilde{\Lambda}:=\tilde{U}^*\Lambda$. The previous identity says that
\begin{equation} \label {ell and ell0}
 \frac{2}{1+U}\ell_0=\tilde{\ell} \quad \mbox { or } \quad -i(i+X)\ell_0=\tilde{\ell}
\end{equation}

 Since $\tilde{U}_{L\cap W}=id$ we get that $\tilde{U}^*$ takes the decomposition $L=\ell\oplus L\cap W\oplus \Lambda$ to an orthogonal decomposition $L_W=\tilde{\ell}\oplus L\cap W\oplus \tilde{\Lambda}$. The operators $S\in\Sym{(L)}$ with the given block decomposition go via conjugation by $\tilde{U}$   to operators $\tilde{S}\in\Sym{(L_W)}$ with the same type of block decomposition relative $L_W=\tilde{\ell}\oplus L\cap W\oplus \tilde{\Lambda}$.

We realize, looking at Lemma $\ref {transition}$, that due to dimension constraints the only thing one needs to prove is that the equation:
\begin {equation} \label {alg}
\frac {1+U^*}{2}B \frac{1+U}{2}= S
\end {equation}
has only the trivial solution $B_i=0$, $S_i=0$, where
\[ B= \left( \begin{array}{ccc}
B_1 &  0 & B_2^*\\
  0 & 0 &  B_3^*\\
 B_2 & B_3 & B_4
\end{array}
\right)\quad \mbox { and } \quad
S = \left( \begin{array}{ccc}
0 & S_1^* & 0\\
 S_1& S_2& 0 \\
 0 & 0 & 0
\end{array}
\right)
\]
The main point here is that the block decomposition of $B$ is relative $L_W=\ell_0\oplus L\cap W\oplus JV^\perp$ and decomposition of $S$ is relative $L_W=\tilde{\ell}\oplus L\cap W\oplus \tilde{\Lambda}$.  Notice that $(\ref {alg})$ can be written as 
\[B=-(i-X) S (i+X)
\]
This is the same thing as 
\[ \langle Bv,w \rangle =-\langle (i-X)S(i+X)v,w\rangle=\langle S(i+X)v,(i+X)w \rangle, \quad \forall v,w\in L_W
\]
We take first $v\in \ell_0$ and $w\in L\cap W$. Relation $(\ref {ell and ell0})$ and  $X=0$ on $L\cap W$ imply
\[ 0=-i\langle S(i+X)v, w\rangle =-i\langle S_1(i+X)v,w\rangle
\]
We conclude that $S_1\equiv 0$.  Similarly taking $v,w\in L\cap W$ we get $S_2\equiv 0$ which finishes the proof.\\

 b) Let us notice that we have two tautological bundles over $\Lag^W(k)$ namely 
 \[
 \begin {array}{ccc}
  \vartheta &\hookrightarrow & \Lag^W(k)\times W \\
  \downarrow & &\downarrow \\
  \vartheta^\omega &\hookrightarrow &\Lag^W(k) \times W^\omega
  \end{array}
 \]  
whose fiber at $L$ consists of $L\cap W$ and $L\cap W^\omega$ respectively. We have of course that $\vartheta$ is a subbundle of $\vartheta^\omega$ and if we let  $\theta$ be the orthogonal complement of $\vartheta$ in $\vartheta^\omega$, then the bundle described in the statement is none other than $\Sym{(\vartheta)}\oplus \Hom{(\theta,\vartheta)}$. Hence it is the direct sum of a complex bundle, always naturally oriented and the bundle of self-adjoint endomorphisms associated to a complex bundle. But this last one is up to isomorphism the bundle associated to the principal bundle of unitary frames with the action of conjugation of the unitary group on its Lie algebra (the gauge transformations bundle). So it is clearly orientable.

In the case $k=1$, $\vartheta$ is a line bundle and $\Sym{(\vartheta)}$ is oriented by the identity.
\qed\\

\medskip

\noindent
\textbf{Orientation conventions II:}\index{Orientation conventions II} Notice that we have the following inclusion $Z_{\{k\}}\subset\Lag^{W_{k-1}}(1)$ of manifolds of codimension $2k-1$. In Section \ref{cohomology ring and geometrical representatives} we showed how the finite dimensional Schubert cells $Z_{\{k\}}(n)$ induce a coorientation on $Z_{\{k\}}$. We would like to say that the coorientation described there is the same as the natural orientation of the algebraic complement of $T\Lag^{W_{k-1}}(1)$ explained above. 

The connection between the two is Proposition \ref{useful for orientation} .  It is enough to consider the finite dimensional case, i.e. we will work with $U(n)$, $\Lag(n)$ and $Z_{\{k\}}(n)$.

Let $U_k\in Z_{\{k\}}(n)$ be the orthogonal reflection with the $-1$ eigenspace given by $\langle e_k\rangle$.
 We  consider that the coorientation at $U_{k}$ of $Z_{\{k\}}(n)\subset U(n)$  is induced by the  linear form on $T_{\id}U(n)$ of rank $2k-1$ 
\[  \tag{*} \theta_{k} \wedge \bigwedge_{j<k}( \alpha_{jk}\wedge\beta_{jk}).
\]
This form can be transported to a linear form on $T_{U_k}U(n)$ by the differential at $0$ of the map
\[ \xymatrix {  \Sym(n) \ar[r] \ar @/_1pc/ [rr]_{\Cay_k} & \mathscr{U}(n) \ar^{U_{k}\cdot}[r] & \mathscr{U}(n),}
 \quad \quad A\mapsto \frac{i-A}{i+A}=:U \mapsto U_k U
\]

The reason for which we can use the Cayley transformation, in place of the exponential map as we did in Section \ref{cohomology ring and geometrical representatives} is because the differential at $0$ of the Cayley transform
\[ A\mapsto \frac{i-A}{i+A}
\]
is the identity, after we identify $\mathfrak{u}(n)$ with $\Sym(n)$ via multiplication by $-i$.

We consider the coorientation on $Z_{\{k\}}(n)\subset \Lag(n)$ at $H_{k}^+$ to be the one obtained from the natural orientation on $\Sym(\langle f_k\rangle)\oplus \Hom(\langle f_k\rangle,\langle e_1,\ldots ,e_{k-1}\rangle)$.
 By Proposition  \ref{useful for orientation}, in order for the Cayley graph map $\mathscr{O}:U(n)\rightarrow\Lag(n)$ to be coorientation preserving, i.e.,  to take the coorientation at  $U_{k}$ of $Z_{\{k\}}(n)\subset U(n)$ to the coorientation at $H_{k}^+$ of $Z_{\{k\}}(n)\subset \Lag(n)$ we need to check that the (differential of the) map
 \[\Sym(n)\rightarrow\Sym(H_{k}^+),\quad \quad \quad A\mapsto \tilde{U}_kA\tilde{U}_k^{-1}
 \]
 is coorientation preserving, i.e. that it takes the form $(*)$ to a positive multiple of the orientation form of 
 \[ \Sym(\langle f_k\rangle)\oplus \Hom(\langle f_k\rangle,\langle e_1,\ldots ,e_{k-1}\rangle)
 \]
 This is straightforward.\qed\\

 \bigskip

\subsection{Local intersection numbers} \label{IntFor}

We are now ready to do intersection theory on $\Lag^-$. 

 In this section $M$ will be a closed, oriented manifold of fixed dimension $2k-1$, unless otherwise stated. This is the codimension  of the Schubert variety $\overline{Z}_{k}$. Let $F:M\rightarrow \Lag^-$ be a smooth map.  We will call such a map a (smooth, compact) family of lagrangians. 
   
 In Section \ref{cohomology ring and geometrical representatives} we defined the  transversality of $F$ to $\overline{Z}_{k}$ to be transversality on every stratum in the \textit {standard stratification}. In the case when $M$ has complementary dimension this implies that $F$ can only meet the top stratum,  $Z_k^\circ$.  In this section we weaken this condition of transversality by defining a new stratification on $\overline{Z}_k$ whose top stratum contains $Z_k^\circ$.
  
We saw in the previous section that if
\[ H^-:=W_0\supset W_1\supset W_2\supset \]
is a complete, decreasing flag of $H^-$, we have
 \[
  Z_{k}^\circ\subset \Lag^{W_{k-1}}(1)\subset \overline{Z_{k}}
 \]
where $\Lag^{W_{k-1}}(1)$ is a smooth submanifold.  

 \begin{definition} The \textit{\textbf{non-standard stratification}}\index{non-standard stratification}  on $\overline{Z_k}$:
 \[  \overline{Z_k}:=S_0\supset S_2\supset S_3 \supset 
 \]
 has at its highest stratum the manifold $S_0\setminus S_2:=\Lag^{W_{k-1}}(1)$, while the other strata, $S_{i}\setminus S_{i+1}$ are unions of $Z_J\subset \overline {Z_{k}}$ each of which  has  codimension $(2k-1)+i$ in $\Lag^-$.
 
 A function $F:M\rightarrow\Lag^-$ is (non-standard) transversal to $\overline{Z}_k$ if it is transversal to every stratum in the non-standard stratification.
 \end{definition}
 
 \begin{remark} In the rest of this paper the stratification on $\overline{Z}_k$ is the non-standard and the notion of transversality we use is the one adapted to this stratification.
\end{remark}
\begin{remark} The two stratifications on the Schubert variety $\overline{Z}_k$ define the same cohomology class in $\Lag^-$. See Remark \ref{diff strat} and appendix \ref{appendix B} for details.
\end{remark}
 
  By Proposition \ref{algebraic complement},  $\Lag^{W_{k-1}}(1)\hookrightarrow \Lag^-$ has an algebraic complement $N\Lag^{W_{k-1}}(1)$ which is naturally oriented as follows. 
  
   Let $L\in \Lag^{W_{k-1}}(1)$, $V:=L\cap W_{k-1}$ and $\ell$ be the orthogonal complement of $L\cap W_{k-1}$ in $L\cap W_{k-1}^\omega$.  The algebraic complement to $T_L\Lag^{W_{k-1}}(1)$ is the vector subspace  of $\Sym{(L)}$ of operators coming from 
 \[ \Sym(V)\oplus \Hom(\ell, V)
 \] 
  The space $V$ is one dimensional and so $\Sym(V)$ is a one dimensional \textit {real} vector space, naturally oriented by the identity map. A non-zero operator $A\in\Sym(V)$ is positively oriented if the following number is positive:
 \[ \langle Av,v \rangle \mbox{ for any  } v\in L\cap W_{k-1}
 \]
 The canonical orientation on $\Hom(\ell, V)$ is given by the following data. Let  $v$  be a unit vector in $V$ and $\{g_1,g_2,\ldots ,g_{k-1}\}$ be a complex orthonormal basis for $\ell$. We say that a basis $T_1,\ldots T_{2k-2}$ is positively oriented for $\Hom(\ell, V)$  if the following determinant is positive:
 
 \[ \left | \begin{array}{ccccc}
 \Real \langle T_1g_1,v \rangle & \Imag \langle T_1g_1,v \rangle & \ldots & \Real \langle T_1g_{k-1},v \rangle & \Imag \langle T_1g_{k-1},v \rangle\\
 \Real \langle T_2g_1,v \rangle & \Imag \langle T_2g_1,v \rangle & \ldots & \Real \langle T_2g_{k-1},v \rangle & \Imag \langle T_2g_{k-1},v \rangle\\
 \ldots & \ldots & \ldots &\ldots &\ldots \\
 \Real \langle T_{2k-2}g_1,v \rangle & \Imag \langle T_{2k-2}g_1,v \rangle & \ldots & \Real \langle T_{2k-2}g_{k-1},v \rangle & \Imag \langle T_{2k-2}g_{k-1},v \rangle
 
 \end{array}
 \right |
 \]
 One can check that the orientation does not depend on the choice of $v$ or of the basis $\{g_1,g_2,\ldots ,g_{k-1}\}$.
 
 The following is straightforward:
 \begin{lemma} For a self-adjoint operator $S\in \Sym(L)$ the $\Sym(V)\oplus \Hom(\ell, V)$ block is described in the orthonormal basis $\{v,g_1,g_2,\ldots ,g_{k-1}\}$ as the operator:
 \[ \begin{array} {ccc}
  v &\rightarrow &\langle Tv,v\rangle v \\
 g_1&\rightarrow &\langle Tg_1,v\rangle v \\
  & \ldots & \\
  g_{k-1} &\rightarrow &\langle Tg_{k-1},v\rangle v
 \end{array}
 \]
  \end{lemma}

This lemma and the previous observations prompts the following definition:

\begin{definition}\label{determinant} 
Let $F:M\rightarrow \Lag^-$ be an oriented family of lagrangians of dimension $2k-1$ transversal to $\overline{Z_k}$ and let $p\in F^{-1}(\overline{Z_k})=F^{-1}(\Lag^{W_{k-1}}(1))$ be a point in $M$.

 Let $\{\epsilon_1, \ldots \epsilon_{2k-1}\}$ be an oriented basis for $M$ at p, $v$ be a unit vector in $F(p)\cap W_{k-1}$ and $\{g_1,g_2,\ldots ,g_{k-1}\}$ be a complex orthonormal basis of $\ell(p)$, the orthogonal complement of $F(p)\cap W_{k-1}$ in $F(p)\cap W_{k-1}^\omega$. 

The \textit{\textbf{intersection number}}\index{intersection number} at $p$ of $F$ and $\overline{Z_k}$, denoted $\sharp (M\cap \overline{Z_k})_p$ is the sign of the determinant
\[ 
\left | \begin{array}{ccccc}
\langle  d_pF(\epsilon_1)v,v\rangle & \Real \langle d_pF(\epsilon_1)g_1,v \rangle & \ldots & \Imag \langle d_pF(\epsilon_1)g_{k-1},v \rangle \\
\langle  d_pF(\epsilon_2)v,v\rangle & \Real \langle d_pF(\epsilon_2)g_1,v \rangle & \ldots & \Imag \langle d_pF(\epsilon_2)g_{k-1},v \rangle \\
\ldots & \ldots &\ldots &\ldots \\
\langle  d_pF(\epsilon_{2k-1})v,v\rangle & \Real \langle d_pF(\epsilon_{2k-1})g_1,v \rangle & \ldots & \Imag \langle d_pF(\epsilon_{2k-1})g_{k-1},v \rangle 
\end{array}
\right | 
\]
\end{definition}
 
 Theorem $\ref{intersection}$ implies the following:

\begin{prop} Let $F:M\rightarrow \Lag^-$ be a smooth family of lagrangians transversal to $Z_k$. The following equality holds:
\[ 
\int_M F^*([ {Z_k}, \omega_{k}])=\sum_{p\in  F^{-1}(\overline{Z_k})}\sharp (M\cap \overline{Z_k})_p
\]
where the integral represents the evaluation of a cohomology class on the fundamental class of $M$.
\end{prop}

It is useful to have a formula for the intersection number in terms of projections.
\begin {lemma} \label{inters formula} 
Let $F:M\rightarrow\Lag^-$ be a smooth family of lagrangians transversal to $Z_k$ and let $P:\Lag^-\rightarrow \mathscr{B}(\hat{H})$ be the smooth map that takes a lagrangian to its orthogonal projection. Denote by $P_F:M\rightarrow\mathscr{B}(\hat{H})$ the composition $-JP\circ F$. The intersection number $\sharp (M\cap \overline{Z_k})_p$ is equal to the sign of the determinant:
\[ 
\left | \begin{array}{ccccc}
\langle  d_pP_F(\epsilon_1)v,v\rangle & \Real \langle d_pP_F(\epsilon_1)g_1,v \rangle & \ldots & \Imag \langle d_pP_F(\epsilon_1)g_{k-1},v \rangle \\
\langle  d_pP_F(\epsilon_2)v,v\rangle & \Real \langle d_pP_F(\epsilon_2)g_1,v \rangle & \ldots & \Imag \langle d_pP_F(\epsilon_2)g_{k-1},v \rangle \\
\ldots & \ldots &\ldots &\ldots \\
\langle  d_pP_F(\epsilon_{2k-1})v,v\rangle & \Real \langle d_pP_F(\epsilon_{2k-1})g_1,v \rangle & \ldots & \Imag \langle d_pP_F(\epsilon_{2k-1})g_{k-1},v \rangle 
\end{array}
\right | 
\]
\end{lemma}
\noindent
\textbf{Proof:} By Lemma \ref{projection derivative} we have 
\[ d_P(\dot{S})=\left(\begin{array}{cc}
0 & \dot{S}J_L^{-1} \\
J_L\dot{S} &  0
\end{array}
\right)
\]  and  so $(-Jd_LP)|_{\Sym(L)}=id_{\Sym(L)} $.  Since all the vectors $v,g_1,\ldots ,g_{k-1}$ belong to $L:=F(p)$ the proposition follows.
\qed\\

The intersection numbers when $k=1$ have received a particular attention.\\
\noindent
\textbf{Notation:}  Let $\Mas:=\Lag^{H^-}(1)=\{L~|~\dim L\cap H^-=1\}$ be the top stratum in the non-standard stratification of $\overline{Z}_1$.

The notation is justified by the following definition.
\begin{definition} For every family $F:S^1\rightarrow \Lag^-$, transversal to $\overline{\Mas}$ the intersection number
\[ \sum_{p\in F^{-1}(\Mas)} \sharp (M\cap \Mas)_p
\]
is called the \textit{\textbf{Maslov index}}.
\end{definition}

\begin{prop}
 The Maslov index is a homotopy invariant that provides an isomorphism:
 \[ \pi_1(\Lag^-)\simeq \mathbb{Z}
 \]
 \end{prop}
 \noindent
 \textbf{Proof:} This is obvious in the light of the fact that the Maslov index is the evaluation over $S^1$ of the pull-back of the cohomology class determined by $\Mas$,  namely $[Z_{1}, \omega_1]$.
 \qed

\begin{lemma} The Maslov index of a family $F:S^1\rightarrow \Lag^-$ can be computed by the formula:
\[ \sum_{p\in F^{-1}(\Mas)} sgn \langle -J\dot {P}_pv_p,v_p\rangle
\]
where $v_p$ is a non-zero vector in $F(p)\cap H^-$ and $\dot{P}_p$ is the derivative at $p$ of the family of the associated projections. 
\end{lemma}
\noindent
\textbf{Proof:} This is just a particular case of $\ref {inters formula}$.
\qed\\
\begin{remark} Although we defined the Maslov index for families of lagrangians parametrized by the circle,  one can use the same definition for families parametrized by the interval $[0,1]$. The Maslov index is then a homotopy invariant of maps with the end-points fixed.
\end{remark}

The following observations lead to an interesting formula. Let $k\geq 1$ and denote by $Z_{\geq k}$ the union of Schubert cells:
\[ Z_{\geq k}:= \bigcup_{i\geq k} Z_i
\]
\begin{lemma} The set $Z_{\geq k}$   is a closed a subspace and a smooth submanifold of codimension $2k-2$ in  $\Mas$. Moreover the following set equalities hold:
\begin{itemize}
\item[(a)] $Z_{\geq k}=\mathscr{R}^{-1}(\mathbb{P}(W_{k-1}))$ where $\mathscr{R}$ is the generalized reduction
\[ \mathscr{R}:\Mas\rightarrow \mathbb{P}(H^-)
\]
\item[(b)] $Z_{\geq k}=\Mas\cap \Lag^{W_{k-1}}(1)$
\item[(c)] $\overline{ Z}_{\geq k}=\overline{ Z}_{ k}$
\end{itemize}
\end{lemma}
\noindent
\textbf{Proof:} The first set equality proves that  $Z_{\geq k}$ is a smooth manifold of codimension $2k-2$ 
in  $\Mas$ since the generalized reduction is a vector bundle and $\mathbb{P}(W_{k-1})$ is a closed submanifold of $\mathbb{P}(H^-)$. \qed\\
\begin{remark}  A new stratification of $\overline{Z}_k$ with $Z_{\geq k}$ being the top stratum does not turn $\overline{Z}_k$ into a quasi-submanifold because $Z_{\geq k}$ does not contain the codimension $1$ stratum $Z_{1,k}$.\qed
\end{remark}

\begin{definition}\label{strong transversal} A smooth $2k-1$ dimensional family $F:M\rightarrow\Lag^-$ is \textit{\textbf{strongly transversal}}\index{strongly transversal} to $\overline{Z}_k$ if the following  conditions hold
\begin{itemize}
\item $F$ is  transversal to $\overline{Z}_k$ 
\item $F$ is transversal to $\overline{Z}_1$ 
\item $F^{-1}(\overline{Z}_1)=F^{-1}(\Mas)$
\end{itemize}
\end{definition}

\begin{remark} The first and the third  conditions of strong transversality imply that $F^{-1}(\overline{Z}_k)=F^{-1}(Z_{\geq k})$. Indeed, the first condition implies that $F^{-1}(\overline{Z}_k)=F^{-1}(\Lag^{W_{k-1}}(1))$, whereas the third implies that $F^{-1}(\overline{Z}_k)\subset F^{-1}(\Mas)$.
\end{remark}

\begin{remark} Every smooth family can be deformed to a family that  satisfies the first two transversality conditions.  However the third condition of strong transversality is not amenable to perturbations, since there are topological obstructions to achieving that. An example is a family for which cohomology class $F^*[\overline{Z}_{1,2},\omega_{1,2}]$ is non- trivial.
\end{remark}

Things are good when  $k=2$ since then there are no topological obstructions in that case.
\begin{lemma} \label{three strong}  Let $\dim{M}=3$. Any family $F:M^3\rightarrow\Lag^-$  can be deformed to a strongly transversal family to  $\overline{Z}_2$. 
\end{lemma}
\noindent
\textbf{Proof:} First deform the family to a map transversal to $\overline{Z}_2$ and then move it off   $ \bigcup_{k\geq 2} Z_{1,k}$ which has codimension 4 and has the property that $\bigcup_{k\geq 2} Z_{1,k}=\Lag^{W_1}(1)\setminus\Mas$. \qed 

\begin{prop}\label{strong tautology}
Let $M$ be an oriented, closed manifold and let $F:M\rightarrow\Lag^-$ be a  family, strongly transversal to $\overline{Z}_k$. Then $M^1:=\{m\in M~|~\dim {F(m)\cap H^-}=1 \}$ is a closed, cooriented submanifold of $M$ of dimension $2k$. Let  $\gamma\subset  M^1\times \mathbb{P}(H^-)$ be  the tautological bundle over $M^1$ with fiber  $\gamma_m:=L\cap H^-$. Then 
\[ \int_{M} F^*[\overline{Z}_{k},\omega_{k}]= \int_{M^1} c_1(\gamma^*)^{k-1}
\]
\end{prop}
\noindent
\textbf{Proof:} Notice that $M^1=F^{-1}(\Mas)$. The fact that $M^1$ is a cooriented submanifold of $M$  of codimension $1$ follows from the second condition of strong transversality and the fact that $\Mas$ is a cooriented submanifold of $\Lag^-$ of codimension $1$.  It is closed because it is the preimage of $\overline{Z}_1$ by the third condition of strong transversality.

We have the following commutative diagram
\[\xymatrix{ 
F^{-1}(\overline{Z}_k) \ar@{^{(}->}[d] \ar^{\quad F}[r]     & Z_{\geq k} \ar^{\mathscr{R}\quad}[r]       \ar@{^{(}->}[d]                       & \mathbb{P}(W_{k-1})\ar@{^{(}->}[d]\\ 
M^1\ar@{^{(}->}[d]\ar^{F |_{M^1}\quad}[r]                    & \Mas \ar^{\mathscr{R}}[r] \ar@{^{(}->}[d] & \mathbb{P}(H^-) \\
M \ar^{F}[r]                                              & \Lag^-                                         &
}
\]
The local intersection number of $M$ and $Z_{\geq k}$ in $\Lag^-$ at a point $m\in M$ is the local intersection number of $M^1$ and $\Mas$ at $m\in M$, which is the local intersection number of $\mathscr{R}\circ F|_{M^1}$ with $\mathbb{P}(W_{k-1})$ at $m\in M$. 

Let $\tau^*$ be the dual to the tautological bundle of $\mathbb{P}(H^-)$. Then the Poincar\'e dual of $\mathbb{P}(W_{k-1})$ in $\mathbb{P}(H^-)$ is $c_1(\tau^*)^{k-1}$ where $c_1$ is the first Chern class of $\tau^*$. The total intersection number of $M^1$ and $\mathbb{P}(W_{k-1})$ is the evaluation of the pull-back of the Poincar\'e dual to $\mathbb{P}(W_{k-1})$ on $M^1$. The following observation finishes the proof 
\[  \gamma=(\mathscr{R}\circ F|_{M^1})^*(\tau)
\]
  \qed

\section{Applications} 
In this chapter we aim to relate the theory developed so far to index theory. We introduce a criterion for deciding when a general family of operators is continuous/differentiable. We describe   how the  Atiyah-Singer classifying space for $K^{-1}$, that is a certain connected component of the set of bounded, self-ajoint operators relates to $\Lag^-$.   The section on the odd Chern character is based on standard results and is designed to make the connection between the fundamental cohomology classes of $\Lag^-$ and index theory. In the last two sections we give concrete local intersection formulae for different families of self-adjoint operators.

\bigskip

\subsection{Differentiable families}\label{differentiable families}

In practice, families of lagrangians come from closed, self-adjoint operators. In order to be able to do differential topology  one needs an easy criterion to decide when these families are differentiable, especially when one has in mind to work with operators whose domain varies, such as elliptic boundary problems. In this section we give such a criterion  and some examples.

We start by recalling Kato's definition of differentiability.  Let $B$ be a smooth manifold.

\begin{definition} [\textit{Kato}, \cite {K}, Ch.VII-1.2] Let $(T_b)_{b\in B}$ be a family of closed,densely defined, self-adjoint  operators with domains $D(T_b)_{b\in B}$ . The family is continuous/differentiable if there exist continuous/differentiable families of \textit{bounded} operators $S_b,R_b:H\rightarrow H$ such that $\Range(S_b)=D(T_b)$ and $T_bS_b=R_b$ for all $b$.
\end{definition}   

Our operators will always be Fredholm so we concentrate on them.  

\noindent
\textbf{Notation:} Let $\SFred$ \index{$\SFred$}be the set of all closed, densely-defined, self-adjoint, Fredholm operators.

To each closed, self-adjoint, Fredholm operator one can associate its switched graph. More precisely one has a map:
\[ \tilde{\Gamma}:\SFred\rightarrow \Lag^-,\quad \quad\quad T\mapsto  \tilde{\Gamma}_T:=\{(Tv,v)~|~v\in D(T)\}
\]
\begin{lemma} Let $T:B\rightarrow \SFred$ be a family of operators. Suppose $\tilde{\Gamma}\circ T:B\rightarrow \Lag^-$ is continuous/differentiable. Then $T$ is continuous/differentiable in the sense of Kato. 
\end{lemma}
\noindent
\textbf{Proof:} If $\tilde{\Gamma}\circ T$ is continuous/differentiable then the family of  corresponding unitary operators $U_b=\mathscr{C}^{-1}(\tilde{\Gamma}( T_b))$ is continuous/differentiable and we take 
\[S_x=-i(1-U_x),\quad\quad\quad \quad  R_x=1+U_x\]
in Kato's definition.\qed\\

In order to state the main result of this section we introduce some terminology. We keep the notations from the previous sections.

\begin{definition}
A symplectic operator on $\hat{H}$ is a bounded, invertible operator  that satisfies the relation
\[ X^*JX=J
\]
Let $\eSy(\hat{H})$ be the group of all symplectic operators on $\hat{H}$.
\end{definition}

\begin{lemma} 
\begin{itemize} 
\item [(a)] Let $L\subset \hat{H}$ be any closed subset of $\hat{H}$ and $X:\hat{H}\rightarrow \hat{H}$ be any bounded, invertible operator. Then
\[ (XL)^\perp=(X^*)^{-1}L^\perp
\]
\item[(b)]The group $\eSy(\hat{H})$ acts on $\Lag$.
\item[(c)] Let $L\in \Lag$ and $X\in\eSy(\hat{H})$ an let $\psi_L:\hat{H}\rightarrow \hat{H}$ be defined by
\[\psi_L(X)=X\circ P_L+(X^*)^{-1}\circ P_{L^\perp}.
\]
Then  the reflection in $XL$ is given by $R_{XL}=\psi_L(X)\circ R_L\circ \psi_L(X)^{-1}$.
\end{itemize}
\end{lemma}
\noindent
\textbf{Proof:} $(a)$ First $X^*(XL)^\perp\perp L$ because if $x\in L$ and $w\in (XL)^\perp$ then
\[ \langle X^*w,v\rangle=\langle w, Xv\rangle= 0
\]
This is the same thing as $X^*(XL)^\perp\subset L^\perp$. Replacing $L$ by $L^\perp$ one also gets that $X^*(XL^\perp)^\perp\subset L$. We claim that 
\[ X^*(XL)^\perp + X^*(XL^\perp)^\perp=L+L^\perp=\hat{H}
\]
and therefore the previous two inclusions are in fact equalities.  Indeed, since all the spaces involved are closed
 \[ ((XL)^\perp+(XL^\perp)^\perp)^\perp=XL\cap XL^\perp=\{0\}
 \] and so $(XL)^\perp+(XL^\perp)^\perp=\hat{H}$ and the claim follows.\\
 
 $(b)$ We have $JXL=(X^*)^{-1}JL=(X^*)^{-1}L^\perp=(XL)^\perp$.\\
 
 $(c)$ One checks  immediately that for $v\in L$ and $w\in L^\perp$
 \[ R_{XL}\psi_L(X)(v)=\psi_L(X)R_L(v)=Xv, \quad R_{XL}\psi_L(X)(w)=\psi_L(X)R_L(w)=-(X^*)^{-1}w.
 \]
 
\begin{prop}  The group $\eSy(\hat{H})$ is a Banach-Lie group modelled on the infinite Lie algebra of operators $Y$ that satisfy:
\[ Y^*J+JY=0
\]
The natural action $\eSy(\hat{H})\times \Lag\rightarrow \Lag$ is differentiable. 
\end{prop}

\noindent
\textbf{Proof:} The first part is standard. We prove the second part in two steps. First, we show that if one fixes a lagrangian, say $H^+$, the map
\[\eSy(\hat{H})\rightarrow \Lag,\quad \quad \quad X\rightarrow XH^+
\]
is differentiable. By Corollary \ref{reflection diff} this is equivalent to proving that the map of associated reflections is differentiable. By part $(c)$ of the previous lemma this is obvious.

Second, we notice that if $X\in\eSy(\hat{H})$ and $L\in \Lag$ then
\[ XL= X\mathscr{O}(H^+, (\mathscr{C})^{-1}(L))H^+
\]
where $(\Gamma^{\small{c}})^{-1}$ is the inverse of the Cayley graph map.  The map 
\[ \eSy(\hat{H})\times \Lag\rightarrow \eSy(\hat{H}), \quad \quad (X,L)\rightarrow X\mathscr{O}(H^+, (\mathscr{C})^{-1}(L))
\] is differentiable since $\mathscr{O}$ and $(\mathscr{C})^{-1}$ are.
\qed\\

The following theorem gives a useful, general criterion for differentiability. Let $B$ be the open, unit ball in $\mathbb{R}^n$
\begin{theorem} \label{differentiability} 
 Let $(T_x)_{x\in B}:D(T_x)\subset H\rightarrow H$  be a family indexed by $B$ of  densely-defined, closed, self-adjoint, Fredholm operators. Let $H_0:=D(T_0)$ and suppose $H_0$ comes equipped with an inner product such that:
\begin{itemize}
\item[(1)] the inclusion $H_0\rightarrow H$ is bounded and
\item[(2)] the operator $T_0:H_0\rightarrow H$ is bounded.
\end{itemize} 
Suppose there exists a differentiable family of bounded, invertible operators $U:B\rightarrow \mathscr{G}L(H)$ such that 
\begin{itemize}
\item[(a)] $U_x^*(H_0)=D(T_x)$;
\item[(b)] the new family of operators $\tilde{T}_x:=U_xT_xU_x^*$  is a differentiable  family of bounded  operators in $\mathscr{B}(H_0,H)$.
\end{itemize}
 Then the family of switched graphs associated to $(T_x)_{x\in B}$ is differentiable at zero in $\Lag^-$. 
\end{theorem}
\noindent
\textbf{Proof:}  Let us notice that the family of operators on $\hat{H}$
\[ \hat{U}_x= \left(\begin{array}{cc} U_x &0\\ 0 &(U_x^*)^{-1}
\end{array}\right)
\]
is a differentiable family of symplectic operators such that $\hat{U}_x\tilde{\Gamma}_{T_x}=\tilde{\Gamma}_{\tilde{T}_x}$. Hence it is enough to prove the differentiability of the family $\tilde{T}_x$. 
We will suppose from now on that all operators $T_x$ are defined on the same domain such that  $(a)$ and $(b)$ are satisfied for $U_x=\id$.

  Suppose now that $\Ker {T}_0=\{0\}$. Then $\tilde{\Gamma}_{T_0}$ is in the Arnold chart $\mathcal{A}_{H^+}$ and for $\|x\|$ small enough the switched graphs $\tilde{\Gamma}_{T_x} $ are in the same open set since $\tilde{\Gamma}_{T_x}\cap H^+=\Ker{T_x}$ and $L\rightarrow L\cap H^+$ is an upper semi-continuous function. Each $T_x$ in this smaller set has a self-adjoint inverse $
S_x:H\rightarrow H_0\subset H$. The differentiability of the family $T_x$ in $\mathscr{B}(H_0,H)$ is equivalent with the differentiability of the family $S_x\in \mathscr{B}(H,H_0)$. This implies the same property for $S_x$ seen as a family in $\mathscr{B}(H,H)$ since $i:H_0\rightarrow H$ is differentiable. 
Since the switched graph of $T_x$ is the graph of $S_x$ it follows that $\tilde{\Gamma}_{T_x}$ is differentiable. \par

  The way to put us in the situation $\Ker {T}_0=\{0\}$ is by adding a real constant $\lambda$ to the family $T_x$ constant for which $\Ker{(T_0+\lambda)} =\{0\}$. We can do this because $T_0$ is Fredholm. In order to justify that this does not change anything we look back at Arnold's isomorphism. We claim that if the unitary operator $U_x$ corresponds to $T_x$, i.e. $\mathscr{C}(U_x)=\tilde{\Gamma}_{T_x}$, then
\[
U_x(\lambda):=(2i+\lambda(1-U_x))^{-1}(2iU_x+\lambda (1-U_x))
\]
corresponds to $T_x+\lambda$. Indeed "the denominator", $2i+\lambda(1-U_x)$ is invertible for every real constant $\lambda$ since $\frac{2i+\lambda }{\lambda }$ has modulus one if and only if $\lambda\in -i+\mathbb{R}$. One can easily verify that $U_x(\lambda)U_x(\lambda)^*=I$. It is also a matter of routine to check that 
\[
\left(\begin{array}{c}
(1+U_x(\lambda))v\\
-i(1-U_x(\lambda))v
\end{array}
\right)=
\left(\begin{array}{c}
((1+U_x)-\lambda i(1-U_x)) [2i(2i+\lambda(1-U_x))^{-1}v]\\
(-i(1-U_x))[2i (2i+\lambda(1-U_x))^{-1}v]
\end{array}
\right)
\]
and this proves the claim. 
\qed \\

 The next result gives a practical way to decide when condition $(b)$ in the previous theorem is satisfied for a family of operators defined on the same domain. 
 \begin{lemma}\label{norm diff} Let $H_0\subset H$ be a dense subset in $H$,  are endowed with  an inner products such that the inclusions $H_0\subset H$ is continuous.  Let $T, g_2: B\rightarrow \mathscr{B}(H_0,H)$, $ g_1:B\rightarrow\mathscr{B}(H),$   be three families of operators such that
 \begin{itemize}
 \item[(a)] $g_1(0)=\id$ and  $g_2(0)=0$;
 \item[(b)] $T(x)=g_1(x)T(0)+g_2(x),\quad \forall x\in B$.
 \end{itemize}
 Suppose $g_1$ and $g_2$ are differentiable at zero. Then $T$ is differentiable at zero. 
\end{lemma} 
\noindent
\textbf{Proof:}  Let $dg_i$, $i\in\{1,2\}$ be the differentials at $0$ of $g_i$. Then 
\[ \|g_1(x)T_0+g_2(x)-T(0)-dg_1(x)T_0-dg_2(x)\|\leq \|(g_1(x)-g_1(0)-dg_1(x))\|_{H} \|T_0\|_{H_0,H}+\]
\[+\|(g_2(x)-g_2(0)-dg_2(x))\|_{H_0,H}
\]
Dividing by $\|x\|$ and taking the limit  $x\rightarrow 0$ finishes the proof.
\qed\\

\begin{corollary} If in the previous proposition $g_2$ is a continuous map of bounded operators $g_2:B\rightarrow\mathscr{B}(H,H)$ then the claim stays true.
\end{corollary}
\noindent
\textbf{Proof:} This is based on the fact that $\|T\|_{H_0,H}\leq \|T\|_{H,H}$.\qed\\

\begin{definition}\label{affine family} A family of  operators $(T_x)_{x\in B}:H_0\rightarrow H$ is called affine if 
\[ T_x-T_0\in\mathscr{B}(H),\quad\quad\quad \forall x\]
\end{definition}
\begin{corollary} An affine family of operators $(T_x)_{x\in B}$ is differentiable at $0$  if the associated family of bounded operators $T_x-T_0$ is differentiable.
\end{corollary}

\begin{corollary} \label{univ family} The universal family $T: U(N)\rightarrow \Lag^-(L^2[0,1])$ (see Example \ref{ex univ family}) is differentiable.
\end{corollary}
\noindent
\textbf{Proof:}
We use the previous criterion to prove differentiability at $1$. Take a chart at $1$ in $U(N)$. For example one can take $\phi:\Sym(N)\rightarrow U(N)$ be the Cayley transform or $\phi(A)=e^{iA}$. For every $A\in \Sym(N)$ let $U_A:C^{\infty }_{\phi(A)}([0,1])\rightarrow C^{\infty }_{1}([0,1])$ be the operator defined by:
\[U_A(f)(t)=\phi(tA)^{-1}f(t)=\phi(-tA)f(t)\]
One checks easily that these operators extend to a differentiable family of unitary operators $U: \Sym(N)\times L^2[0,1]\rightarrow L^2[0,1]$. The resulting family is:
$$U_AT_{\phi(A)}U_A^{-1}(f)=-i\frac{df}{dt}-i\phi(-tA)\left(\frac{d}{dt}\phi(tA)\right)f$$
This is a differentiable family of bounded operators $L^{1,2}([0,1])\rightarrow L^2([0,1])$ by the previous corollary. \qed

\bigskip

\subsection{Index Theory}\label {index}

This section is inspired by \ref{N3}.

The classifying space for odd $K$-theory, $\Lag^-$ is not the usual space one uses in index theory. In their work \cite{AS2}, Atiyah and Singer looked at the homotopy type of the space $\BFred\subset\Sym$ of bounded, self-adjoint, Fredholm  operators endowed with the norm topology. They proved it has three connected components $\BFred_+$, $\BFred_{-}$ and $\BFred_*$ characterized by :
\begin{itemize}
\item $T\in\BFred_{\pm} \Leftrightarrow$   T has only positive/negative essential spectrum 
 \item $T\in \BFred_* \Leftrightarrow $  T has both positive and negative essential spectrum 
 \end{itemize}
The first two components, $\BFred_{\pm}$ are contractible and  $\BFred_*$ is  classifying for $K^{-1}$. 

Let $B$ be  compact space. It follows from Atiyah-Singer results that every element in $K^{-1}(B)$ can be represented by the homotopy class of a continuous map 
\begin{equation} \label {analytic index} B\rightarrow \BFred_* \quad \quad \quad\quad \quad \quad b\mapsto  T_b
\end{equation}
The homotopy class of such a family of operators is the analytic index of the family. Notice that 
the map
\[ [0,1]\times B\rightarrow  \BFred_*\quad \quad \quad  b\mapsto  \frac{T_b}{\sqrt{1+tT_b^2}}
\]  
 provides a homotopy between the initial map and the associated family of "zeroth order" operators
\[B\rightarrow  \BFred_* \quad \quad \quad \quad \quad b\mapsto  \frac{T_b}{\sqrt{1+T_b^2}}
\]
 This leads to  the standard trick  that allows one to define an  analytic index  for a family of \textit{unbounded} self-adjoint, Fredholm operators with spectrum stretching to both $\pm\infty$. To make the ideas more precise we introduce the following function
 \[ \Ri:\mathbb{R}\rightarrow\mathbb{R}, \quad \quad\quad \quad \quad x\mapsto \frac{x}{\sqrt {1+x^2}}
 \]
 Recall that $\SFred $ is the set of all closed, densely defined, self-adjoint,  Fredholm  operators  on the Hilbert space $H$. The map
 \[ \Ri :\SFred\rightarrow \BFred \quad \quad \quad \quad \quad T\mapsto \Ri(T)
 \]
 is an injection.
 
 \begin{definition}   The \index{Riesz topology} \textit{\textbf{Riesz topology}} on $\SFred$ is the topology induced by the metric 
 \[ d(T_1,T_2)=\|\Ri(T_1)-\Ri(T_2)\|
 \]
 A function $f:B\rightarrow \mathscr{C}\mathscr{S}$ is called Riesz continuous if it is continuous with respect to the Riesz topology.
 
 The \textit{\textbf{Atiyah-Singer index}}  of a Riesz continuous, family of   operators $T:B\rightarrow\SFred$  for which $\Ri (T(b))\in \BFred_*$ is the homotopy class of  the map
\[ \Ri\circ T: B\rightarrow \BFred_*
\]
\end{definition}
\begin{remark} In order  to define an analytic index in the unbounded case, using the Atiyah-Singer classifying space $\BFred_*$ one \textit {needs}  the family to be Riesz continuous. At the other extreme, if all the operators involved are bounded, then Riesz continuity is equivalent with the norm continuity and  $\Ri\circ T$ is homotopic  with $T$.\qed
\end{remark}

The vertical, Lagrangian Grassmannian, $\Lag^-$ suggests a different approach. A self-adjoint, Fredholm operator, $T:D(T)\subset H\rightarrow H$ bounded or unbounded gives rise to a vertical, Fredholm lagrangian, namely its switched graph.
\[ T\rightarrow \tilde {\Gamma}_{T}:=\{(Tv,v)~|~v\in D(T)\}
\]
\begin{definition} A family $T:B\rightarrow \SFred$ of Fredholm operators is \textit{\textbf {gap continuous}} if the map
\[ \tilde{\Gamma} \circ T: B\rightarrow\Lag^-,\quad \quad \quad b\mapsto \tilde{\Gamma}_{T_b}
\]
is continuous. 

The \textit{\textbf{graph  index}} of a gap continuous family of Fredholm operators $T:B\rightarrow\SFred$ is the homotopy class of the map
\[ \tilde{\Gamma} \circ T: B\rightarrow\Lag^-.
\] 
\end{definition}

\begin{lemma} \label{boundedgap} The map 
\[ \tilde{\Gamma}:\BFred\rightarrow\Lag^- \quad \quad \quad T\rightarrow \tilde{\Gamma}_{T} \]
is continuous hence every  map $T:B\rightarrow \BFred$ which is continuous in the norm is also gap continuous.
\end{lemma}  
\noindent
\textbf{Proof:} The reflection in the switched graph can be computed explicitly in terms of the operator and  Corollary \ref{reflection diff} finishes the proof.\qed\\

\begin{lemma} Riesz continuity is invariant under conjugation, i.e. if $T:B\rightarrow\SFred$ is Riesz continuous and $U:B\rightarrow \mathscr{U}(H)$ is a continuous family of unitary operators then the family $\tilde{T}:B\rightarrow \SFred$ such that $\tilde{T}_b:=U_bT_bU_b^*$  for all $b\in B$ is continuous.
\end{lemma}
\noindent
\textbf{Proof:} This is straightforward in light of
\[ \Ri(\tilde{T}_b)=U_b\Ri{(T_b)}U_b^*
\]

Let $H_0$ be a dense subspace in $H$ such that there exists an inner product $\langle\cdot,\cdot\rangle_0$ that makes the inclusion $H_0\rightarrow H$ continuous.

\begin{remark} If $T:H_0\rightarrow H$ is a bounded operator, than the topology defined by the graph norm of $T$ on $H_0$, (see \ref{graph norm}) is weaker than the topology of the norm $\|\cdot\|_0$. In other words the identity map:
\[ (H_0,\|\cdot\|_0)\rightarrow (H_0,\|\cdot\|_g)
\]
is continuous as one can easily see. On the other hand if $T$ is 
Fredholm  then the norms are equivalent. Indeed  there exists a constant  $C_1>0$ such that
\[ \|v\|_g=\|v\|_H\geq C_1\|v\|_{H_0},\quad\quad \forall v\in \Ker{T}
\]
simply because $\Ker{T}$ is finite dimensional. Moreover there exists $C_2>0$ such that 
\[ \|v\|_g\geq \|Tv\|_H\geq C_2\|v\|_{H_0}, \quad \quad \forall v\in\Ker{T}^\perp 
\]
because $T\bigr|_{\Ker{T}^\perp}:\Ker{T}^\perp\rightarrow \Ran{T}$ is invertible. \qed
\end{remark}

\begin{prop} Let $T:B\rightarrow \SFred(H)$ be a family of self-adjoint, Fredholm operators such that $T_b:H_0\rightarrow H$ is bounded and the family is continuous at $b_0$ as a map $T:B\rightarrow\mathscr{B}(H_0,H)$. Then $T$ is Riesz continuous at $b_0$.
\end{prop}
\noindent
\textbf{Proof:} 
See  Proposition 1.7 in \cite{N3} and   Theorem VI.5.12 in \cite{K}.\qed

\begin{definition} Let $T:B\rightarrow \SFred$ be a family of self-adjoint, Fredholm operators, let $b_0\in B$ and  let $H_0:=D(T_{B_0})$ be endowed with the graph norm of $T_{b_0}$. The family  is called \textit{\textbf{nice}} at $b_0\in B$ if  there exist a \textit {continuous} family of unitary operators $U:B\rightarrow\mathscr{U}(H)$ such that 
\begin{itemize}
\item[(a)] $H_0=U_bD(T_b)$;
\item[(b)] The new family $(\tilde{T})_{b\in B}$, $\tilde{T}_b:=U_bT_bU_b^*$ is continuous as a family of bounded operators $T:B\rightarrow\mathscr{B}(H_0,H)$.
\end{itemize}

A family is called nice if it is nice at every point.\qed
\end{definition}

The following result is a consequence of what we said above.
\begin{prop} Every nice family of operators is Riesz continuous.
\end{prop}

\begin{example} Every continuous, affine family of operators is nice.
\end{example}

\begin{example} The universal family is nice.
\end{example}

\begin{prop} Every Riesz continuous family of operators is gap continuous.
\end{prop}
\noindent
\textbf{Proof:} See Lemma 1.2 in \cite{N3}. \qed\\

\begin{theorem} The graph map $\tilde{\Gamma}:\BFred_*\rightarrow\Lag^-$ is a weak homotopy equivalence and for every Riesz continuous family $T:B\rightarrow \SFred$ the Atiyah-Singer index coincides with the graph index.
\end{theorem}
\noindent
\textbf{Proof:} See Proposition 3.1 and Theorem 3.3 in \cite{N3}.\qed\\

 The previous Theorem says that the graph index is the appropriate notion of analytic index one has to look at.
 
 \begin{definition} Let $F:M\rightarrow\SFred$ be a family of  self-adjoint, Fredholm operators parametrized by a compact topological space $M$. Then $F$ is said to be continuous if $\tilde{\Gamma}\circ F$ is continuous. The analytic index of a continuous family $F$, denoted $[F]$ is the homotopy class   of the map
\[\tilde{\Gamma}\circ F:M\rightarrow \Lag^-
\] 
\end{definition}
\begin{remark} All nice families of operators are continuous. 
\end{remark}

\bigskip

\subsection{The Chern Character}\label {Chern}

 Let $M$ be a finite, CW-complex, hence compact.  The Chern character is a ring homomorphism:
\[ \ch: K^0(M)\rightarrow H^{\even}(M,\mathbb{Q})
\]
 The suspension isomorphism, which is actually taken to be the definition of  $K^{-1}$,  helps us extend the Chern character to the odd case:
 \[ \xymatrix{
  \tilde{K}^{-1}(M)  \ar[d]_{\ch}  \ar @{=} [r]^{\Sigma} & \tilde{K}^{0}(\Sigma M) \ar[d]^{\ch}  \\
 \tilde{H}^{\odd}(M,\mathbb{Q})  \ar  [r]^{\Sigma \quad}        &    \tilde{H}^{\even}(\Sigma M,\mathbb{Q})         }
 \]
 It is well-known that $U(\infty)$ is a classifying space for $K^{-1}$. Hence every element in $\tilde{K}^{-1}(M)$ can be represented by the homotopy class of a (pointed) map $f:M\rightarrow U(\infty)$.
 Let $[f]\in K^{-1}(M)$ be the element this map represents. Then $\Sigma f:\Sigma M\rightarrow \Sigma U(\infty)$ represents an element in $\tilde{K}^0(\Sigma M)$ which corresponds to $f$ via the suspension isomorphism. The previous commutative diagram can be written as
 \begin{equation} \label{chern character1}
 \Sigma\circ \ch {[f]}=\ch{([\Sigma f])}
 \end{equation}
 
 A short digression is necessary at this point. The space $\Sigma U(\infty)$ comes with a principal $U(\infty)$ bundle $\breve{U}$, namely the bundle obtained with the clutching map given by the identity. More precisely one starts with the trivial $U(\infty)$ bundle over $[0,1]\times U(\infty)$ and identifies  $(0,U, g)$ with $(1,U, Ug)$  for all $(U,g)\in U(\infty)\times U(\infty)$. 
 
 This is an old acquaintance of ours.  Indeed the pull-back of this bundle to $\Sigma U(n)$ is nothing else but the bundle frame bundle associated to the vector bundle $E_n$ which we considered in Section \ref{cohomology ring and geometrical representatives}. 
 
  Another way of looking at these bundles is via the periodicity map (see \cite{Q}, page 224-225)
 \[ \Sigma U(n)\rightarrow \Gr (n,2n)\hookrightarrow \Gr(n,\infty)\simeq BU(n)
 \] where the first map is given explicitly as follows 
 \[ [0,\pi]\times U(n) \rightarrow \Gr(n,2n),\quad \quad \quad (t,U)\rightarrow  
 \cos{t}  \left(\begin{array}{cc} 1  & 0 \\  0 & -1 \end{array} \right)-\sin{t} \left(\begin{array}{cc} 0 & U^{-1}\\  U &0 \end{array} \right)
 \]
 The right hand side is an involution of $\mathbb{C}^n\oplus \mathbb{C}^n$.
 The bundle $E_n$ is the pull-back of the universal  $U(n)$vector bundle  $EU(n)$. In the same way $\breve{U}$ comes from the universal $U(\infty)$-bundle over $BU(\infty)$. 
 
 Now every continuous map $f:M\rightarrow U(\infty)$ defined on a compact set $M$ is homotopy equivalent with a map (which we denote by the same letter) $f:M\rightarrow U(n)$. The class $[f]\in K^{-1}(M)$ or the class $[\Sigma f]\in \tilde{K}^{0}(\Sigma M) $ can be represented by the bundle $(\Sigma f)^* E_n$  (which determines  a stable isomorphism class). Using equation (\ref{chern character1}) we get that 
 \begin{equation} \label{chern sigma} \ch{[f]} =\Sigma^{-1} \ch((\Sigma f)^* E_n))=\Sigma^{-1}((\Sigma f)^*\ch E_n)
 \end{equation}
 
 The inverse of the suspension isomorphism $\Sigma$ is easy to describe. It is the composition
\[ \xymatrix {
 \tilde{H}^{\even}(\Sigma M, \mathbb{Z})\ar [r]^{\pi^*\quad} & \tilde{H}^{\even}(S^1\times M,\mathbb{Z})\ar [r]^{~~/{dt}} & H^{\even-1} (M,\mathbb{Z})
 }
\]
where $\pi:S^1\times M\rightarrow \Sigma M$ stands for the projection and $/{dt}$ stands for the slant product with the orientation class of $S^1$. So 
\begin{equation} \label{chern sigma2}\Sigma^{-1}((\Sigma f)^*\ch {E_n}) =(\pi^*(\Sigma f)^*\ch {E_n})/{dt}=((\Sigma f\circ p)^*\ch {E_n})/{dt}=\end{equation}
\[=((id_{S^1}\times f)^*\ch{E_n})/{dt}= f^*(\ch{E_n}/dt)
\]
The class $\ch{E_n}/dt\in H^{odd}(U(n),\mathbb{Q})$ is called the \textit{transgression} class of the Chern character. Of course, one can do slant product componentwise and get, for each positive integer $k$ a class:
\[ \ch_{2k-1}^\tau:= \ch_{2k}{(E_n)}/{dt}\in H^{2k-1}(U(n),\mathbb{Q})
\]
There is nothing special about the Chern character. The same transgression process can be applied to any characteristic class of $E_n$, in particular to the Chern classes and we have already done this in Section \ref{cohomology ring and geometrical representatives} where we denoted those classes by $x_i$. We use a different notation now which is more appropriate to this context.
\[ c_{2k-1}^\tau:=c_k(E_n)/{dt}\in H^{2k-1}(U(n),\mathbb{Z})
\] 
There is a very simple relation between $\ch_{2k-1}^\tau$\index{$\ch_{2k-1}^\tau$} and $c_{2k-1}^\tau$:

\begin{lemma} 
\[ \ch_{2k-1}^\tau= \frac{(-1)^{k-1}}{(k-1)!} c_{2k-1}^\tau
\]
\end{lemma}
\noindent
\textbf{Proof:} First of all, $\ch_{2k}{(E_n)}\in \tilde{H}^{2k}(S^1\times U(n),\mathbb{Q})$ is a polynomial in the variables \linebreak $c_1(E_n),\ldots ,c_k(E_n)$ and the coefficient of $c_k(E_n)$ is $(-1)^{k-1}/(k-1)!$. On the other hand, every element in $H^{2k}(S^1\times U(n),\mathbb{Z})$ is a sum:
\[ z=x + Dt\wedge y
\]
where $x\in H^{2k}(U(n),\mathbb{Z})$, $y\in H^{2k}(U(n),\mathbb{Z})$ and $Dt\in H^1(S^1,\mathbb{Z})$ satisfies $Dt(dt)=1$. We claim that for every characteristic class of $E_n$ its $H^{2k}(U(n),\mathbb{Z})$ component vanishes. Indeed the class $x$ is the pull-back of $z$ via the inclusion $\{1\}\times U(n)\rightarrow S^1\times U(n)$ and the claim follows by noticing that the pull-back of the bundle $E$ is trivial over $U(n)$.

We conclude that the cup product of any two characteristic classes of $E_n$ is zero and so we have 
\[ \ch_{2k}{(E_n)}=\frac{(-1)^{k-1}}{(k-1)!} c_{k}(E_n)
\]
which after taking the slant product gives the identity we were after.
\qed\\

Suppose now that $M$ is a closed, oriented manifold and $f:M\rightarrow\Lag^-$. Theorem \ref{Schubert  and canonical} says that the pull-back $f^*x_{k}=f^*[Z_{k},\omega_{k}]$. On the other hand 
 by the previous lemma, relations (\ref{chern sigma}), (\ref{chern sigma2}) and Proposition \ref{pull back and preimage} we have the following result
 
 \begin{prop} \label{chern character}  Let $M$ be a closed manifold and let $f:M\rightarrow \Lag^-$ be a smooth map transversal to $Z_k$.  The  following  holds:
\[  \ch_{2k-1}([f])= \frac{(-1)^{k-1}}{(k-1)!} f^*[Z_k,\omega_k]=\frac{(-1)^{k-1}}{(k-1)!} [f^{-1} (Z_k), f^*\omega_k]
 \]
 \qed
\end{prop}

\medskip

Let us take now $M:=S^{2N-1}$. On one hand we have an isomorphism:
\[ \pi_{2N-1}(\Lag(N))\rightarrow \pi_{2N-1}(\Lag^-)
\] 
because if $W_N$ is a subspace of codimension $N$ in $H^-$ then  $\Lag^-\setminus\Lag^{W_N}$ has codimension $2N+1$ in $\Lag^-$, being equal with the Schubert variety $\overline{Z}_{N+1}$ and so every map  $S^{2N-1}\rightarrow\Lag^-$ can be homotoped to a map $S^{2N-1}\rightarrow \Lag^{W_N}$. The later space is a vector space over $\Lag(N)$. 

On the other hand we have a morphism
\[\tag{*} \pi_{2N-1}(\Lag(N))\rightarrow H^{2N-1}(S^{2N-1},\mathbb{Z}), \quad \quad [f]\rightarrow f^{*}[\overline{Z}_{N},\omega_{N}]
\]
Moreover this morphism is injective. Indeed, by what was said above we have
 \[  f^{*}[\overline{Z}_{N},\omega_{N}]=c^{\tau}_{2N-1}(\Sigma f)^*E_N=\Sigma^{-1}(c_{N}(\Sigma f)^*E_N)\]
 This means that $f^{*}[\overline{Z}_{N},\omega_{N}]=0$ if and only if the Euler class of the complex bundle $(\Sigma f)^*E_N$ over $S^{2N}$ is zero. But the Euler class is the only  obstruction to trivializing a rank $N$ complex bundle over $S^{2N}$. Hence the classifying map $\Sigma f$ has to be homotopically trivial. On the other hand, via Bott periodicity the suspension map
 \[ \pi_{2N-1}(U(N)) \rightarrow \pi_{2N}(\Sigma U(N)), \quad \quad [f]\rightarrow [\Sigma f]
 \]
 is an isomorphism. Hence if $[\Sigma f]=0$ then $[f]=0$.
 
 The morphism $(*)$ is not surjective. If we compose it with the isomorphism obtained by integrating over the fundamental class of $S^{2N-1}$, 
 \[H^{2N-1}(S^{2N-1},\mathbb{Z})\rightarrow \mathbb{Z},\quad \quad \quad \alpha\mapsto\int_{S^{2N-1}}\alpha
 \] 
 then we get a morphism $ \pi_{2N-1}(\Lag(N))\rightarrow \mathbb{Z}$. Its image is in the subgroup $(n-1)!\mathbb{Z}$ by Bott divisibility theorem which is saying that the Chern character of every rank $N$ complex vector bundle over $S^{2N}$ is an integer, because it is the index of the twisted signature operator. (see Theorem IV.1.4 in \cite{LM}) This implies that the Euler class of that bundle is divisible by $(N-1)!$. In fact the image is the whole subgroup $(n-1)!\mathbb{Z}$ since there exists an element  $a\in\tilde{K}(S^{2N})$ such that $\ch{a}=\PD{(pt)}$ (see Theorem 24.5.3 in \cite{Hir}). We can represent this by a map $\Sigma{f}:S^{2N}\rightarrow BU(N)$. The clutching map associated to the $U(N)$-principal bundle over $S^{2N}$ induced by $\Sigma{F}$ is the (desuspension) map $f:S^{2N-1}\rightarrow U(N)$ whose Chern character is $\PD{(pt)}$.
 
 The previous discussion leads to the following result.

 \medskip
 
 \begin{theorem} \label{homotopy type and cycles} The map 
 \[\Pi_N: \pi_{2N-1}(\Lag^-)\rightarrow\mathbb{Z},\quad\quad\quad \Pi_N([f:S^{2N-1}\rightarrow \Lag^-])=\int_{S^{2N-1}}f^*[\overline{Z}_{N},\omega_{N}]
 \]
 is injective and the image is the subgroup $(n-1)!\mathbb{Z}$.
  \end{theorem}
  \begin{corollary}
 The homotopy type of a map $f: S^{2N-1}\rightarrow \Lag^-$ is determined by the integer
 \[ \int_{S^{2N-1}}f^*[\overline{Z}_{N},\omega_N]
 \]
 which is always divisible by $(N-1)!$.  If $f$ is transversal to $\overline{Z}_{N}$ then this integer is the total intersection number of $f$ and $\overline{Z}_N$.
 \end{corollary}

\begin{remark} Any map $f:S^{2N-1}\rightarrow \Lag^-$ can be deformed to a map $S^{2N-1}\rightarrow \Lag(N)$. After identifying $\Lag(N)$ with $U(N)$ one gets a map $\Lag(N)\rightarrow S^{2N-1}$ coming from the fibration $p: U(N)\rightarrow S^{2N-1}$. The degree of the composition $p\circ f:S^{2N-1}\rightarrow S^{2N-1}$ is exactly the integer from the theorem. 
 \end{remark}
 
\bigskip

\subsection{Intersection for Families of Operators} \label{InterFamOp}

The intersection formulae in the  Section \ref{IntFor} were given in terms of the differential of the family of lagrangians or, what is more or less the same thing, the differential of the associated projections. In practice these lagrangians are switched graphs of self-adjoint, Fredholm operators.   In order to adapt those intersection formulae to the case of operators, the first thing to do is  to make sure that we know what we mean by the differential of a family of operators.  This is clear in the case when all the operators involved are bounded and that is what we do next. In the case when the operators are unbounded but the family is affine (see Definition \ref{affine family}) then the formulae of this section hold with minimal changes.

 The following is straightforward.
\begin{lemma} Let $T$ be a bounded self-adjoint, Fredholm operator. Then its switched graph
\[ \tilde{\Gamma}_T:=\{(Tv,v)~|~v\in H\}
\] is a vertical, Fredholm lagrangian.
\end{lemma}
Let $M$ be a  smooth manifold. Recall a definition:
\begin{definition}\label{differentiability def} Let $F:M\rightarrow\SFred$ be a family of  self-adjoint, Fredholm operators. Then $F$ is said to be smooth/ continuous if $\tilde{\Gamma}\circ F$ is smooth/continuous. The \textit{\textbf{analytic index}}\index{analytic index} of a continuous family $F$, denoted $[F]$ is the homotopy class   of the map
\[\tilde{\Gamma}\circ F:M\rightarrow \Lag^-
\] 
\end{definition}
\begin{remark} All families of operators satisfying the conditions of Theorem \ref{differentiability} are \linebreak smooth/continuous. 
\end{remark}

The  homotopy class of a continuous family $F$ determines an element in $K^{-1}(M)$ also called  the index.
\begin{definition} The \textit{\textbf{cohomological index}} \index{cohomological index} of a continuous family $F$ is  denoted by $\ch [F]\in H^{\odd}(M,\mathbb{Q})$ and represents the cohomology class obtained by applying the Chern character to the analytic index. 
\end{definition}

Let $W$ be a  codimension $k-1$ subspace  of $H^-$. We consider the associated $2k-1$ codimensional  cocycle whose underlying space is the following Schubert variety 
\[ \overline{Z}_W=\{L\in\Lag^-~|~\dim{L\cap W}\geq 1\}
\] 
\begin{definition} A smooth family $F:M\rightarrow\SFred$ is said to be in general position with respect to $W$ \index{general position}  if  $\tilde{\Gamma}\circ F$ is  transversal to the Schubert variety $\overline{Z}_W$ with the non-standard stratification.
\end{definition}
 If $M$ has complementary dimension to $\overline{Z}_W$, i.e., $\dim {M}= 2k-1$, the condition to be in general position with respect to $W$ implies that there are only a finite number of points $p\in M$ such that 
 \begin{equation} \label{dker} \dim \tilde{\Gamma}_{F(p)}\cap W=1
 \end{equation}
   This means that
  \[ \dim{\Ker{(F(p))}\cap W}=1
  \]

 \noindent
 \textbf{Notation:}  Let  $F:M\rightarrow \SFred$ be a smooth family  in general position with respect to $W$. For every $p\in M$ such that $\dim{\Ker{(F(p))}\cap W}=1$ denote by $\epsilon_p\in \{\pm 1\}$ the intersection number at $p$ of $\tilde{\Gamma}\circ F$ with $\hat{Z}_W:=\Lag^{W}(1)=\{L~|~\dim{L\cap W}=1\}$. 
 
\begin{theorem} \label{chern and poincare dual} Let $M$ be a closed oriented manifold of dimension $2k-1$,  let \linebreak $F:M\rightarrow\SFred$ be a smooth family of self-adjoint, Fredholm operators and let  $W\subset H$ be a codimension $k-1$ subspace such that $F$ is in general position  with respect to $W$. Denote by $M_W$ the set $M_W:=\{p\in M~|~\dim{\Ker{(F(p))}\cap W}=1\}$. Then 
\[ \PD \ch_{2k-1}([F])=\frac{(-1)^{k-1}}{(k-1)!} \sum_{p\in M_W} \epsilon_pp 
\] where the term on the left is the Poincar\'e dual to the $2k-1$ component of the cohomological index. 
\end{theorem}
\noindent
\textbf{Proof:} This is a restatement of Theorem \ref{intersection} using Proposition \ref {chern character}. \qed\\

Our main goal in this section is to give a formula for the intersection numbers $\epsilon_p$. This is a local problem.  We first take up the case of bounded operators. A simple but important result is

\begin{lemma} If $T\in\Sym(H)$ is a bounded self-adjoint, Fredholm operator the projection $P^-:\tilde{\Gamma}_T\rightarrow H^-$ is a Banach space isomorphism.
\end{lemma}
\noindent
\textbf{Proof:} Straighforward.\qed \\

Let $B$ be the unit ball in $\mathbb{R}^{2k-1}$.
The next result relates the operator differential to the graph differential.
\begin{lemma} \label{differential}
 Let $T:B \rightarrow \Sym (H)$ be a family of bounded, self-adjoint, Fredholm operators, differentiable at zero. Then the family of switched graphs $(\tilde{\Gamma}_{T_x})_{x\in B }$ is differentiable at zero. Moreover, for every unit vector $v\in\mathbb{R}^n$,  the following equality  holds between the graph and the operator partial derivatives of the family at $0$
\[
P_0^-\circ \frac{\partial\tilde{\Gamma}}{\partial v}\Bigr |_{0}\circ (P_0^-)^{-1}=(1+T_0^2)^{-1}\circ \frac {\partial T}{\partial v}\Bigr |_0 \quad \in\Sym(H)
\] 
Here $P_0^-$ is the projection of the switched graph of $T_0$ onto $H^-$. 
\end{lemma}
\noindent
\textbf{Proof:} For differentiability see Theorem \ref{differentiability}.

For $\|x\|$ small the switched graph of $T_x$ is in the Arnold chart of $\tilde{\Gamma}_{T_0}$. Therefore it is  the graph of an operator $JS_x:\tilde{\Gamma}_{T_0}\rightarrow J\tilde{\Gamma}_{T_0}$, where $S_x\in\Sym{(\tilde{\Gamma}_{T_0})}$. We fix such an $x$. We are looking for an expression for $P_0^-S_x(P_0^-)^{-1}$ as an operator on $H$. 

 For every vector $v\in H$ we have a decomposition:
 \[ (T_xv,v)=(T_0z,z)+J(T_0y,y)=(T_0z+y, z-T_0y)
\]
It is not hard to see that 
\[\begin{array}{ccc}
y&=& (1+T_0^2)^{-1}(T_x-T_0)v\\
v &=& (1+T_0 T_x)^{-1}(1+T_0^2) z
\end{array}
\]
The last relation makes sense, since $1+T_0T_x$ approaches  the invertible operator $1+T_0^2$. The operator $P_0^-S_x(P_0^-)^{-1}:H\rightarrow H$ is nothing else but the correspondence $z\rightarrow y$ hence the expression:
\[ P_0^-S_x(P_0^-)^{-1}= (1+T_0^2)^{-1}(T_x-T_0) (1+T_0 T_x)^{-1}(1+T_0^2)
\]
Differentiating this expression with respect to $x$ finishes the proof.
\qed\\

In order not to repeat ourselves we give the following
\begin{definition} A  smooth family of bounded, self-adjoint, Fredholm operators $F:B\rightarrow \Sym(H)$   is called  \textit{\textbf{localized}}\index{localized family} at $0$ with respect to $W$ if the following two conditions hold
\begin{itemize}
\item  $F$ is in general position with respect to $W$;
\item $(\tilde{\Gamma}\circ F)^{-1}(\hat{Z}_W)=\{0\}$
\end{itemize} 
\end{definition}

The fact that switched graph of $F(0)$ is in $Z_W$  implies that $1\leq\dim{\Ker{F(0)}}\leq k$ by Corollary \ref{dimension}.

We treat first a particular non-generic case. 

\begin{prop} \label {fam oper} Let  $F: B\rightarrow \BFred$ be a  family of  self-adjoint, Fredholm operators localized at $0$ with respect to $W$. Suppose that $\dim{\Ker {F(0)}}=k$. Let  $\phi\in\Ker{F(0)}\cap W$ be a unit vector, let $\phi^\perp$ be the orthogonal complement of $\langle v\rangle$ in $\Ker {F(0)}$ and let  $\{\psi_1,\ldots \psi_{k-1}\}$ be an orthonormal basis of $\phi^\perp$.
 
 The intersection number, $\epsilon_0$, is given by the sign of the determinant:
\[
\left | \begin{array}{ccccc}
\langle    {\partial_1 F}\phi,\phi\rangle & \Real \langle \partial_1 F \psi_1,\phi \rangle & \ldots & \Imag \langle  \partial_1 F \psi_{k-1},\phi \rangle \\
\langle  \partial_2 F\phi,\phi\rangle & \Real \langle  \partial_2 F\psi_1,\phi \rangle & \ldots & \Imag \langle \partial_2 F\psi_{k-1},\phi \rangle \\
\ldots & \ldots &\ldots &\ldots \\
\langle  \partial_{2k-1} F\phi,\phi\rangle & \Real \langle \partial_{2k-1} F\psi_1,\phi \rangle & \ldots & \Imag \langle \partial_{2k-1} F\psi_{k-1}, \phi \rangle 
\end{array}
\right | 
\]
where $\partial_{i}F$ is the partial derivative of $F$ at zero in the  $i$-th coordinate direction of $\mathbb{R}^{2k-1}$. 
\end{prop}
\noindent
\textbf{Proof:} Let $\tilde{F}:=\tilde{\Gamma}\circ F$.

Since $\dim\Ker {F(0)}=k$ we get that  
\[ \tilde{\Gamma}_{JF(0)}\cap W^\omega=\tilde{\Gamma}_{JF(0)}\cap H^-=\Ker{F(0)}\] 
 and so the vectors $g_1,\ldots , g_{k-1}$ in the definition of the intersection number \ref{determinant} are all in the domain of $F(0)=H^-$ and we can take them all in $\Ker {F(0)}$.
We want to replace the partial derivatives of $\tilde{F}$ in that intersection formula with the partial derivatives of $F$. 

The claim that proves the lemma is: 
\[ \langle d_0\tilde{F}(x)g, \phi \rangle = \langle \frac{\partial F}{\partial x}(0)g , \phi\rangle
\]
for every unit vector $x$ and every $g\in\langle\{\phi, \psi_1,\ldots ,\psi_{k-1}\}\rangle$. In order to prove the claim let  $P_0^-$ be the projection of the switched graph of $ F(0)$ onto $H^-$ and let $w:=(1+F_0^2)^{-1}\circ \frac {\partial F}{\partial x}(0)g$. Then 
\[(P_0^-)^{-1}\circ (1+F_0^2)^{-1}\circ \frac {\partial F}{\partial x}\Bigr |_0 \circ P_0^- (0,g)=(F_0w, w)
\]
 Therefore,  by using Lemma $\ref {differential}$ we get 
\[ \langle d_0\tilde{F}(x)g, \phi\rangle =\langle (F_0w, w), (0,\phi) \rangle= \langle w,\phi \rangle 
\]
Then 
\[ \langle w,\phi \rangle =\langle \frac {\partial F}{\partial x}(0)g, (1+F_0^2)^{-1}\phi \rangle=\langle \frac {\partial F}{\partial x}(0)g, \phi \rangle
\]
The last equality holds because $\phi\in \Ker {F_0}$.
\qed\\

   In the case $k=2$, the intersection numbers still have a quite simple description.  Suppose for now that $B$ is the three dimensional ball.
      
 \begin{prop} \label{part case}
  Let $T: B\rightarrow\BFred $ be a family of bounded, self-adjoint, Fredholm operators. Let $e\in H$ be a vector and suppose that $T$ is localized at $0$ with respect to $\langle e \rangle^\perp$.   Let $0\neq\phi$ be a generator of $ \Ker{T_0}\cap \langle e \rangle^\perp$. Then only one of the two situations is possible
  \begin{itemize}
  
  \item[I] $\dim{\Ker T_0}=1$, in which case let $\psi$ be a non-zero vector satisfying the following two relations
   \begin{equation} \label{IntEq}
  \left\{ \begin{array}{ccc}
 \langle \phi, \psi \rangle &= &0 \\
 T_0\psi &= &e 
\end{array}
 \right.
 \end{equation}
 \item[II] $\dim{\Ker T_0}=2$, in which case let $\psi\in\Ker{T_0}$ be a  non-zero vector such that $\psi\perp \phi$. 
 \end{itemize}
 Then the intersection number, $\epsilon_0$ of $T$ with $\overline{Z}_{e^\perp}$ is given by the determinant
 
\[ \left | \begin{array} {ccc}
 \langle \partial_1T\phi,\phi\rangle & \Real  \langle  \partial_1T\psi,\phi\rangle & \Imag  \langle \partial_1T\psi,\phi\rangle \\
 \langle \partial_2T\phi,\phi\rangle & \Real  \langle \partial_2T\psi,\phi\rangle & \Imag  \langle \partial_2T\psi,\phi\rangle \\
 \langle \partial_3T\phi,\phi\rangle & \Real  \langle \partial_3T\psi,\phi\rangle & \Imag  \langle \partial_3T\psi,\phi\rangle
 \end{array}
 \right |
 \]
 where $\partial_iT$ is the directional derivative of $F$ in the $i$-th coordinate direction of $\mathbb{R}^3$.
\end{prop}

 \noindent
 \textbf{Proof:} Let $W=\langle e\rangle^\perp$. The intersection of the switched graph of $T_0$ with $W^\omega$ is two dimensional. Hence the kernel of $T$ is either one or two-dimensional.  
  One vector in the intersection $\tilde{\Gamma}_T\cap W^\omega$ is $(0,\phi)$.  If the kernel of $T$ is two dimensional, then $\tilde{\Gamma}_T\cap W^\omega=\Ker{T}$ and so the second vector in the intersection formulae is a generator of the orthogonal complement in $\Ker {T}$ of $\langle \phi \rangle$. This is the condition imposed on $\psi$ in this situation.
  
  In the latter case the equation $T\psi=ae_1$ has solutions if and only if $a=0$.

 If $\dim{\Ker T_0}=1$, the condition $(T_0\alpha,\alpha)\in W^\omega$ imposes that $T_0\alpha= ae$ for some constant $a$. Of course we are looking for a solution when $a\neq 0$, since otherwise $\alpha$ is a multiple of $\phi$.   At any rate the projection $\psi$ of $\alpha/a$ to $\Ker {T_0}^\perp$ is an element of $W^\perp\subset W^\omega$ and a generator of the orthogonal complement of $\Ker{T_0}$ in $\tilde{\Gamma}_T\cap W^\omega$. It satisfies the two conditions we imposed on $\psi$.
 
 The fact that one can replace the partial derivatives of the switched graphs in the $g_1=(T\psi,\psi)$ direction by the partial derivatives of $T$ in the $\psi$ direction is a computation exactly as in $\ref {fam oper}$ where we used Lemma $\ref {differential}$.  \qed\\

We state now the general case. And $B$ is again the $2k-1$ dimensional ball.
 \begin{prop} Let  $W\subset H$ be a $k-1$ codimensional subspace and let $T:B\rightarrow\BFred$ be a family of bounded,  self-adjoint, Fredholm operators localized at $0$ with respect to $W$.  Suppose that $\dim{\Ker{T_0}}=p\leq k$. 
 Let $\phi$ be a generator of $\Ker {T_0}\cap W$ and let $\phi_1,\ldots ,\phi_{p-1}$ be a basis of the orthogonal complement of $\phi$ in $\Ker{T_0}$. 
 
 The space $W_T:=W\cap \Ran{T_0}$ has dimension $k-p$. Let  $\psi_1,\ldots ,\psi_{k-p}$ be an orthonormal basis of $ P\Bigr |_{(\Ker{T})^\perp}{T_0}^{-1}(W_T)$.

  Then the intersection number, $\epsilon_0$, of $T$ with $\overline{Z}_W$ is the sign of the determinant
  \[
\left | \begin{array}{ccccc}
\langle    {\partial_1 T}\phi,\phi \rangle &  \langle  \partial_2 T\phi,\phi\rangle  & \ldots &  \langle  \partial_{2k-1} T\phi,\phi\rangle  \\
\Real \langle \partial_1 T \phi_1, \phi\rangle  & \Real \langle  \partial_2 T\phi_1,\phi \rangle & \ldots & \Real \langle \partial_{2k-1} T\phi_1,\phi \rangle  \\
 \ldots & \ldots &\ldots &\ldots \\
 \Imag \langle  \partial_1 T  \phi_{p-1},\phi \rangle & \Imag \langle \partial_2 T\phi_{p-1},\phi \rangle & \ldots & \Imag \langle \partial_{2k-1} T\phi_{p-1}, \phi \rangle \\

\Real \langle \partial_1 T \psi_1, \phi\rangle  & \Real \langle  \partial_2 T\psi_1,\phi \rangle & \ldots & \Real \langle \partial_{2k-1} T\psi_1,\phi \rangle  \\
 \ldots & \ldots &\ldots &\ldots \\
\Imag \langle  \partial_1 T  \psi_{k-p},\phi \rangle & \Imag \langle \partial_2 T\psi_{k-p},\phi \rangle & \ldots & \Imag \langle \partial_{2k-1} T\psi_{k-p}, \phi \rangle \\
\end{array}
\right | 
\] 
 \end{prop}
 \noindent
 \textbf{Proof:} One only  needs to make sense of what the orthogonal complement of $\Ker{T_0}\cap W$ in  $\tilde{ \Gamma}_{T}\cap W^\omega$ is.
 \qed\\
 
 If the reader thinks, as we do, that these expressions for the local intersection numbers do not have a great deal of aesthetic appeal, we will try to make it up by a different global formula. This formula is not always available but we think it is worth writing it down.
 
 We make the following definition based on \ref{strong transversal}
 \begin{definition} A smooth family $T:M\rightarrow \BFred$ of bounded self-adjoint, Fredholm operators is called strongly transversal to $\overline {Z}_k$ if $\tilde{\Gamma}\circ T$ is strongly transversal to $\overline {Z}_k$
 \end{definition}  
  
 \begin{lemma} Let $\dim{M}=3$. Any smooth family $T:M\rightarrow \BFred$ can be deformed to a strongly transversal family to $\overline{Z}_2$.
 \end{lemma}
 \noindent
 \textbf{Proof:} This is just proof of Lemma \ref{three strong} with the addition that one has to make sure that in the course of deformation one stays inside $\BFred$. This is true because the map $\tilde{\Gamma}:\BFred\rightarrow\Lag^-$ is open.\qed\\

 \begin{prop}\label{strong transversality} Let  $M$ be a closed, oriented manifold of dimension $2k-1$ and let  $T:M\rightarrow \BFred$ be a strongly transversal family to $\overline{Z}_k$. Then $M^1:=\{m\in M~|~\dim{\Ker{T_m}}=1\}$ is a closed, cooriented manifold. Let $\gamma\subset M^1\times \mathbb{P}(H)$ be the tautological line bundle over $M^1$ with fiber $\gamma_m=\Ker{T_m}$. Then
 \[ \int_{M} F^*[\overline{Z}_{k},\omega_{k}]= \int_{M^1} c_1(\gamma^*)^k
\]
 \end{prop}
 \noindent
 \textbf{Proof:} This is just Proposition \ref{strong tautology} formulated in terms of operators.
 
We want to describe the coorientation of $M^1$ in concrete terms. Let $m\in M^1$ and $v\in T_pM\setminus T_pM^1 $ be a vector. The vector $v$ is said to be positively oriented if given a curve $\alpha:(-\epsilon,\epsilon)\rightarrow M$ such that $\alpha\cap M^1=m=\alpha(0)$ and $\alpha'(0)=v$, the curve of operators $T\circ\alpha$ has local spectral flow equal to $+1$. This means that the eigencurve determined by alpha has a $0$ eigenvalue at $0$ and the derivative is positive. \qed\\

\bigskip

\subsection{Intersection for Families of Operators II}

The motivating example for this paper was the \textit {spectral flow}.  The idea behind the spectral flow is very simple although to put it in a general differentiable topological framework turns out to be a difficult task. Classically, one starts with  a family of  self-adjoint,  elliptic  operators $A_t$  parameterized by the circle, or by the unit interval. The elliptic operators have discret eigenvalues and  if the family is continuous the eigenvalues of the family also vary in continuous families $(\lambda_j(t))_{j\in \mathbb{Z}}$ called eigencurves. Since the family is compact  only a finite number  of eigencurves  will become zero at some moment in time. The spectral flow  is  the difference between the number of eigencurves that start with a negative and end up with a positive sign and those which start with a positive and end up  with a negative sign. In other words, the spectral flow is a count with sign of the $0$-eigenvalues and to such a  $0$-eigenvalue one associates   the sign of the derivative $\dot {\lambda (t_0)}$ of the eigencurve that contains that $0$. 

\begin{figure}[h]
\centerline{\includegraphics[angle=0, height=1.8in, width=2.5in]{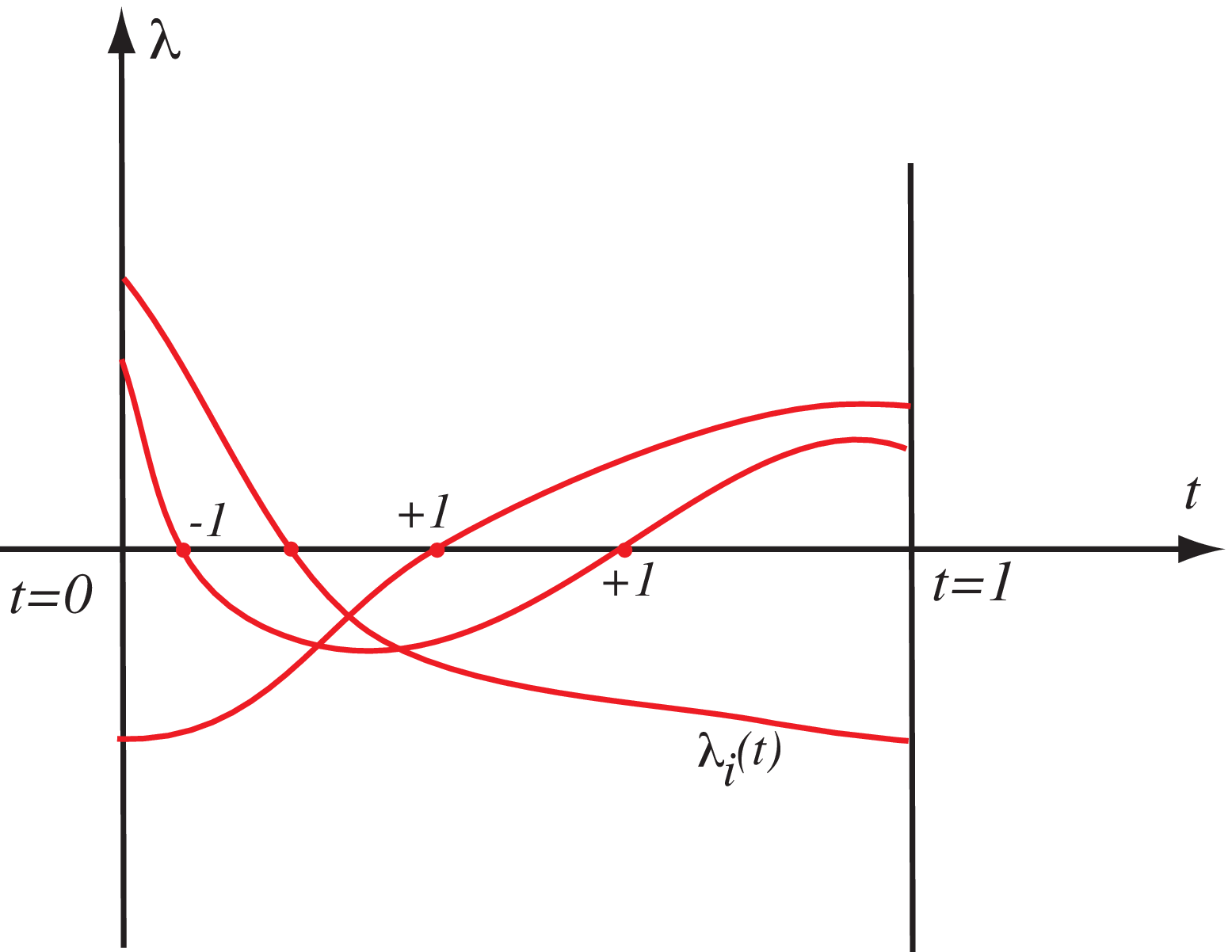}}
\caption{\sl $SF=0$}
\label{fig: 1}
\end{figure}

In order to make the picture rigorous one has to answer certain questions. What does it mean for the family of operators to be continuous? There are several possible answers: in the Riesz topology, in the gap topology or, if one  is dealing with an affine family of Dirac operators for example then one can use the topology of the underlying space.  Another question is when are the eigencurves differentiable around $0$.  One can argue that in order to count those eigencurves that become zero at some moment in time one does not need differentiability but only continuity (one fixes a finite family of eigencurves that  contains all those that go thru zero  and then the spectral flow is $[\sum^+(1)-\sum^-(1)]-[\sum^+(0)-\sum^-(0)]$ where $\sum^\pm(p)$ is the number of positive/negative end points of the eigencurves at $p\in\{0,1\}$).
We are interested in the situations when one can localize the spectral flow, arguably a more useful method of computation.

 If  the family $A_t$  satisfies all the  good conditions one wishes  then the sign of $\dot {\lambda} (t_0)$ coincides with the sign of $\langle \dot A_{t_0} v,v\rangle$ for any vector $v\in \Ker {A_{t_0}}$.  Indeed if $A_t$ and $\lambda_t$ are differentiable families of operators then so is $B_t :=A_t-\lambda_t$. The correspondence $B\rightarrow \Ker{B}$ is a differentiable map  when $\dim{\Ker{B}}$ is constant and so we can, locally around $t_0$, choose a smooth family of unit eigenvectors $v_t$ for $A_t$ and eigenvalue $\lambda_t$. That is  we have the relation:
\[ A_tv_t=\lambda_tv_t
\]
To simplify notations suppose $t_0=0$. We can differentiate at $0$ to get:
\[ \dot{A_{0}}v_{0}+A_{0}\dot{v_{0}}=\dot{\lambda_{0}}v_{0}
\]
We notice of course that $\langle A_{0}\dot{v_{0}}, v_{0} \rangle =0$ since $A_0$ is self-adjoint. So

\[ \langle \dot{A_{0}}v_0,v_0 \rangle=\langle \dot{\lambda_{0}}v_{0}, v_0\rangle=\dot{\lambda_{0}}
\]

The sign of $\langle \dot{A_{0}}v_0,v_0 \rangle$ does not depend on the choice of the vector $v_0\in\Ker{A_{t_0}}$. 

Among the good conditions one wishes of $A$ is the fact that the derivative $\dot{A}_t$ makes sense, which is the case if the operators are bounded, or more generally if the family is affine (see Definition \ref{affine family}). Another useful observation is the fact that, for a smooth family  of  bounded or affine operators $(A_t)_{t\in[0,1]}$, if $A_{t_0}$ is an isolated operator with one-dimensional kernel, meaning that $\Ker{A_t}=0$ for all $t\neq t_0$ close enough, then $\langle \dot{A_{t_0}}v_0,v_0 \rangle\neq 0$ for some vector $v_0\in\Ker{A_{t_0}}$.

These  observations lead us to the following definition. 

\begin{definition} \label{spectral flow} Let $A:[0,1]\rightarrow \SFred$ be a smooth (see Definition \ref{differentiability def}) family of bounded or affine, self-adjoint, Fredholm operators. The family is said to be in general position if the following two conditions hold
\begin{itemize}
\item There are only a finite number of points $t\in[0,1]$ such that $\Ker{A_t}\neq 0$.
\item  For every operator $A_t$ with non-trivial kernel $\dim\Ker{A_t}=1$.
\end{itemize}

Let $\mathscr{Z}:=\{t\in[0,1]~|~\Ker{A_t}\neq 0\}$. The \textit{\textbf{ local spectral flow}} at $t_0\in\mathscr{Z}$  of $A$ is the sign of the number $\langle \dot{A_{t_0}}v,v \rangle$ for some vector $v\in A_{t_0}$.

  The \textit{\textbf{spectral flow}}\index{spectral flow} of the family is the sum of the local flows, i.e.,
\[ \sum_{t_0\in\mathscr{Z}} sign \langle \dot{A_{t_0}}v,v \rangle
\]
\end{definition}

\begin{prop} The spectral flow of a family is the Maslov index of the associated family of switched graphs $\tilde{A}:[0,1]\rightarrow \Lag^-$.
\end{prop}
\noindent
\textbf{Proof:} In the bounded case this is just $\ref {fam oper}$ for $k=1$. The affine case follows the same line of argument.
\qed\\

 We want to show now that the local intersection numbers when one deals with affine families of operators is not more complicated than the bounded case. The specific case we consider are families of Dirac operators defined on the same domain where we let the connection vary.  For more on Dirac operators see \cite {LM}  or  \cite {Roe}. 
  
  Let $Cl(M)\rightarrow M$ be the bundle of Clifford algebras associated to a compact, oriented, Riemannian manifold. Let  $S:\rightarrow M$ be a Cliffford module over $M$. For us, this means that $S$ is a hermitian bundle that  comes with a bundle endomorphism $c:Cl(M)\rightarrow \End{(S)}$, called Clifford multiplication, which is a unitary Clifford algebra representation in each fiber.

 For each hermitian connection $\nabla$ on $S$, compatible with the Levi-Civita connection on $Cl(M)$ one gets a Dirac operator, $\mathscr{D}_{\nabla}$:
 \[   \xymatrix 
 { 
 C^{\infty}(S)\ar^{\nabla\quad\quad}[r] & C^{\infty}(T^*M\otimes S) \ar^{\quad \quad c}[r] & C^{\infty} (S)\\
  } 
 \]
 \[ \mathscr{D}_{\nabla}:=c\circ \nabla
 \]
 It is a known fact that  this operator can be extended to a bounded operator 
 \[
  \mathscr{D}_{\nabla}:L^{1,2}(S)\rightarrow L^2(S)
  \] which is self-adjoint, elliptic. The $L^{1,2}$ inner product on $S$ is defined with the help of a connection but it does not depend on it.
 
 The space of hermitian connections on $S$, $\mathscr{A}(S)$ is an affine space modelled by  $\Omega^1(\Ad{U(S)})$, which is the set of one forms with values in the bundle associated to the principal bundle of orthonormal frames of $S$, via the adjoint representation. On $\Ad{U(S)}$ one has a hermitian metric so one can speak of $L^2$ sections of $T^*M\otimes \Ad{U(S)}$ and this is exactly what   $\Omega^1(\Ad{U(S)})$ will represent  for us, $L^2$ sections rather than just smooth ones.  
 
 Therefore once one fixes a hermitian connection, the space $\mathscr{A}(S)$ becomes a Hilbert space. 
 
\begin{lemma}
The map  \[ 
\mathscr{D}: \mathscr{A}(S)\rightarrow \mathscr{B}(L^{1,2}(S), L^2(S)), \]
\[
\nabla \rightarrow \mathscr{D}_{\nabla}
 \]
 is differentiable.
\end{lemma} 
\noindent
\textbf{Proof:} If one fixes a connection $\nabla_0$ then one gets an induced map:
\[  \Omega^1(\Ad{U(S)})\rightarrow \mathscr{B}(L^2(S), L^2(S))
\]
\[ \nabla-\nabla_0\rightarrow \mathscr{D}_\nabla-\mathscr{D}_{\nabla_0}=c(\nabla-\nabla_0)
\]
This family is clearly differentiable. The rest follows from Lemma \ref{norm diff}.
\qed\\

By Theorem $\ref {differentiability}$ the associated family of switched graphs $\tilde{\Gamma}\circ \mathscr{D}:\mathscr{A}(S)\rightarrow\Lag^- $ is differentiable. We can actually compute this differential explicitly.  In order to do that let us first remember a classical fact about first-order, elliptic operators, namely the elliptic estimates. In our case, for example, there exists a constant $C$, depending only on $\nabla$ such that for every $\phi\in D(\mathscr{D}_{\nabla})$ one has the inequality:
  \[ \|\phi\|_{1,2} \leq C (\|\mathscr{D}_{\nabla}(\phi)\|_2+\|\phi\|_2)
  \]
  This inequality can be rephrased as saying that the map:
  \[
  \gamma:  L^{1,2}(S)\rightarrow \tilde {\Gamma}_{\mathscr{D}_{\nabla}}\subset L^2(S)\oplus L^2(S),\quad \quad \quad \quad
    \phi\rightarrow (\mathscr{D}_{\nabla}\phi,\phi)   
  \]
  is continuous and therefore an isomorphism of Banach spaces.  \footnote {a useful consequence is that  one can change the inequality sign in the elliptic estimates and still get a true sentence}  
  
  We compute the differential of $\tilde{\Gamma}\circ\mathscr{D}$ at a fixed connection $\nabla_0$. The result is a map 
  \[ \Omega^1(\Ad{U(S)}) \rightarrow \Sym{(\tilde{\Gamma}_{\mathscr{D}_0})}
  \]
  The target space of this map are self-adjoint operators on the switched graph of $\mathscr{D}_{\nabla_0}$. We can identify this space with $\Sym{(L^{1,2}(S))}$ via the map $\gamma$.
  \[ 
  \Sym{(\tilde{\Gamma}_{\mathscr{D}_{\nabla}})}\rightarrow \Sym{(L^{1,2}(S))}, \quad \quad \quad\quad
   T\rightarrow \gamma^{-1}T\gamma
  \]
 \begin{definition} If  $\nabla_0$ is a fixed connection then  the following map 
 \[ d\mathscr{D}\bigr |_0:\Omega^1(\Ad{U(S)})\rightarrow \Sym{(L^{1,2}(S))},\quad  \quad  d\mathscr{D}\bigr |_0(\omega):=\gamma^{-1} d(\tilde{\Gamma}\circ\mathscr{D})\bigr|_{\nabla_0}(\omega) \gamma
 \]
 is called the \textit{\textbf{(projected) graph differential}} of $\mathscr{D}$ at $\nabla_0$.
 \end{definition}
  
  \begin{lemma} The following relation holds between the graph differential and the differential of $\mathscr{D}$ at $\nabla_0$:
  \[ d\mathscr{D}\bigr |_0 (\omega)= (1+\mathscr{D}^2_{\nabla_0})^{-1}c(\omega) \]
  \end{lemma}
$~$ \\  
  \noindent
  \textbf{Proof:} Let $\omega\in\Omega^1(\Ad{U(S)})$ be a $1$-form with values in $(\Ad{U(S)})$ and let $\mathscr{D}_{t}:=\mathscr {D_{\nabla_0}}+tc(\omega)$ be the associated affine path of Dirac operators that starts at $\mathscr {D_{\nabla_0}}$. We will denote be $\tilde{\Gamma}_t$ the switched graph of $\mathscr{D}_{t}$.
  
  We want to express the switched graph of $\mathscr{D}_{t}$ as the graph of a bounded, self-adjoint, operator $\tilde{\Gamma}_0\rightarrow \tilde{\Gamma}_0$. So let:
   \[ (\mathscr {D}_t x,x)=(\mathscr {D}_0y,y)+J(\mathscr {D}_0z,z)=(\mathscr {D}_0y,y)+(z,-\mathscr{D}_0z)
  \]
  where $x,y$ and $z$ are in $L^{1,2}(S)$.  In order to solve the system:
\[
  \begin{array}{ccc}
  \mathscr {D}_t x&=&\mathscr {D}_0y+z\\
  x&= &y-\mathscr{D}_0z
  \end{array}
  \]
  we will  suppose first that $x,y\in L^{2,2}(S)$. Then:
  \[\label {composition} \begin{array}{ccc}
  (1+\mathscr{D}_0\mathscr{D}_t)x&=&(1+\mathscr{D}_0^2)y\\
  (\mathscr{D}_t-\mathscr{D}_0)x &= &(1+\mathscr{D}_0^2)z 
  \end{array}
  \]
 Of course, $1+\mathscr{D}_0^2$ is an invertible operator $L^{2,2}(S)\rightarrow L^2(S)$ and because the association $\nabla\rightarrow\mathscr{D}_{\nabla}$ is continuous, so is $1+\mathscr{D}_0\mathscr{D}_t$ for $t$ small enough. Therefore:
 \[\label{operator Eq}
  z=(1+\mathscr{D}_0^2)^{-1} (\mathscr{D}_t-\mathscr{D}_0) (1+\mathscr{D}_0\mathscr{D}_t)^{-1}(1+\mathscr{D}_0^2)y
 \]
 Notice that the operator on the right hand side is pseudo-differential  of order minus one  and as such it can be extended to a continuous operator $L^{1,2}(S)\rightarrow L^{1,2}(S)$. The association
 $y\rightarrow z$ is the operator $\tilde{\Gamma}_0\rightarrow \tilde{\Gamma}_0$ seen only after conjugation with the isomorphism $\gamma$.
 
 The derivative at $t=0$ of this family of operators is exactly the one that appears in the statement of the lemma.
 \qed\\
 
 \begin{remark} The previous computation is almost exactly the same as the one carried in Lemma \ref{projection derivative}. The only thing that is different is the target for the projected graph differential which is a space of self-adjoint operators on $L^{1,2}$ rather than $L^2$. This is true more generally for affine families of operators and the same relation between derivatives holds. More precisely, let $B$ be the unit ball in some $\mathbb{R}^n$ and  $F:B\rightarrow \mathscr{B}(H_0,H)$ be a smooth family of self-adjoint Fredholm operators such that $F(b)-F(0)= A_b\in\mathscr{B}{(H)}$. Then the projection map $P_0:\tilde{\Gamma}_{F(0)}\rightarrow H_0$ identifies the switched graph at $F(0)$ with $H_0$ (with the graph norm) and so the projected graph differential and the "differential" of $F$ satisfy 
 \[ P_0^{-1} d (\tilde{\Gamma}\circ F) P= (1+F(0)^2)^{-1} dA
 \]
 \end{remark}
 
 \begin{prop}  Let $F: \mathbb{R}^n\rightarrow  \mathscr{A}(S)$ be a smooth family of Dirac operators as above, or more generally a family of affine operators. Then  the local intersection formulae of section \ref{InterFamOp}  for bounded operators, still hold with the obvious modifications. For example, in the case of $\ref {part case} $ the intersection number is computed by the sign of the determinant:
  \[ \left | \begin{array} {ccc}
 \langle c(\partial_1f)\phi,\phi\rangle & \Real  \langle c( \partial_1f)\psi,\phi\rangle & \Imag  \langle c(\partial_1f)\psi,\phi\rangle \\
 \langle c(\partial_2f)\phi,\phi\rangle & \Real  \langle c(\partial_2f)\psi,\phi\rangle & \Imag  \langle c(\partial_2f)\psi,\phi\rangle \\
 \langle c(\partial_3f)\phi,\phi\rangle & \Real  \langle c(\partial_3f)\psi,\phi\rangle & \Imag  \langle c(\partial_3f)\psi,\phi\rangle
 \end{array}
 \right |
 \] 
 where $\langle\phi\rangle=\Ker \mathscr{D}_0\cap \langle e\rangle^\perp$ and $\psi$ satisfies the relations:
  \begin{equation} 
  \left\{ \begin{array}{ccc}
  \psi & \neq & 0\\
 \langle \phi, \psi \rangle &= &0 \\
 \mathscr{D}_0\psi &= &ae_1 
\end{array}
 \right.
 \end{equation}
 for some constant $a$, which is $0$ when $\Ker{D}_0$ is two dimensional.
 \end{prop}
\noindent
\textbf{Proof:} The passing from the graph derivative to the  operator derivative is given by the previous lemma. The rest goes just as in the proof of Proposition \ref{fam oper}.
\qed\\

\appendix
\section{Representatives of cohomology classes in Banach manifolds} \label{appendix B}

 We describe in this appendix how certain stratified spaces in an infinite Banach manifold define cohomology classes.  Our presentation   is  inspired from the work of B. Iversen  \cite{Ive} and G. Ruget \cite{Rug}.
 
  In the sequel  our spaces will be assumed paracompact. In fact our ambient space $X$ is assumed to be a metric Banach manifold.  For such a space, $H^k(X)$  will denote the  (\v{C}ech) cohomology of a topological space $X$ with coefficients in the constant sheaf $\bZ$.  If the space is locally contractible then the cohomology with  coefficients in the constant sheaf $\underline{\bZ}$ can be identified with singular cohomology. see \cite[Chap.III]{Bred}. 
  
   For a closed subset $\eC\subset X$ we denote by $H^\bullet_\eC(X)$ the local homology of $X$ along $\eC$ (see \cite[Sec. II.9]{Ive}).  If $S$ and $X$ are locally contractible then,
\[
H^\bullet_\eC(X)\cong H_{sing}^\bullet(X,X\setminus \eC).
\]

One important property of local cohomology is \cite[Prop. II.9.5]{Ive} :
\begin{prop}[excision exact sequence] \label{excision exact} Let $\eC_1\supset \eC_2$ be two closed subsets of the topological space $X$. Then one has the following long exact sequence:
\[ \rightarrow H^k_{\eC_2}(X)\rightarrow H^k_{\eC_1}(X)\rightarrow H^k_{\eC_1/\eC_2}(X/\eC_2)\rightarrow H^{k+1}_{\eC_2}(X)\rightarrow
\]
\end{prop}

\begin{remark} When $\eC_1$ and $\eC_2$ are locally contractible then the previous sequence corresponds to the long exact sequence in singular cohomology associated to the triple $(X, X/\eC_1, X/\eC_2)$.
\end{remark}

   For any closed subset $\eC$  we denote by  $\eH^k_{\eC}$ the sheaf on  $X$ associated to the presheaf $\widetilde{\eH}^k_{\eC}$   such that for any open set $U\subset X$ we have
\[
\Gamma(U\cap\eC,\widetilde{\eH}^k_{\eC})=H^k_{\eC\cap U}(U).
\] 
For every  $x\in X$ the stalk at $x$ of $\eH^k_{\eC}$ is denoted by $\eH^k_{\eC}(x)$ and can be described  by the inductive limit
\[
\eH^k_{\eC}(x):=\varinjlim_{U\ni x} H^k_{U\cap\eC}(U),
\]
where $U$ runs through all the open neighborhoods of $x$. Notice that this sheaf has support on $\eC$ and because of that we have:

\[ H^*(X,\eH^k_{\eC})\simeq H^*(\eC, i^{-1}\eH^k_{\eC})
\]

\begin {definition} The closed space  $\eC$ is said to  have \emph{homological codimension in $X$ at least $c$} if and only if $\eH^k_{\eC}=0$, $\forall k<c$. We write this as
\[
\codim^h_X(\eC)\geq c.
\]
\end{definition}

Observe that the above definition is local, i.e., $\codim^h_X(\eC)\geq c$ if and only if   for some open cover  $\boldsymbol{\eU}$ of  $X$  we have
\[
\codim_U^h(U\cap \eC)\geq c,\;\;\mbox{for any $U\in \boldsymbol{\eU}$}.
\]

\begin{definition}  A closed   subset $\eC\hra X$ is  said to be \emph{normally nonsingular} (or NN)\index{normally non-singular (NN)} of codimension $c$ if the following holds.

\begin{itemize}
\item For any point $w\in \eC$ there exists a neighborhood $\eN$ of $w$ in $X$ and a homeomorphism of pairs
\[
(\eN, \eN\cap\eC) \Lra \bigl(\,\bR^c\times (\eN\cap \eC),\{0\}\times \eN\cap\eC\,\bigr).
\]
\end{itemize}
\end{definition}

\begin{remark}   If $\eC\hra X$ is  NN of codimension  $c$ then $\codim^h_X\eC\geq c$.
\end{remark}

\begin{remark} If $\eC\hra X$ is   submanifold of $X$ of codimension $c$ and $\eC$ is closed as a subset of $X$ then $\eC$ is NN of codimension $c$. In particular, $\codim^h_X\eC\geq c$.
\end{remark}

If $\eC\hra X$ is NN of codimension $c$    the sheaves $\eH^k_{\eC}$ are trivial if $k\neq c$, while if $k=c$ the sheaf $\eH^c_{\eC}$ is locally isomorphic to the constant sheaf  $\underline{\bZ}$.   We say that $\eH^c_{\eC}$ is the \emph{co-orientation sheaf} of $\eC\hra X$ and we  will denote it by $\Bom_{\eC}$.

The  Grothendieck spectral  sequence  for local cohomology (see \cite[Remark 2.3.16]{Dim}) converges to $H^\bullet_\eC(X)$ and its $E_2$ term is given by
\[
E_2^{p,q}= H^{p}(X,\eH_\eC^q)=\begin{cases}
0  &q\leq c\\
H^{p}(\eC,\Omega_{\eC}) & q=c.
\end{cases}
\]

The following extension property is a consequence of has been just said
\begin{prop}[Extension property]   If  $\codim^h_X(\eS)\geq c$ then, for any  closed subset $\eC\supset \eS$ and   any $q<c-1$, the   restriction map
\[
H^q_{\eC}(X)\ra H^q_{(\eC\setminus \eS)}(\eX\setminus \eC)
\]
is an isomorphism.     We will refer to it  as the \emph{extension across $\eS$}.
\label{prop: ext}
\end{prop}

\noindent 
\textbf{Proof:}  By what has just been said $H^q_{\mathscr{S}}(X)=0$ for all $q<c$.
The rest is taken care by the excision exact sequence (Proposition \ref{excision exact}). \qed\\

\begin{corollary} If $\codim_X^h(\eS)\geq c$, then the  for any $q< c-1$ the restriction map
\[
H^q(X)\ra H^q(X\setminus \eS)
\]
is an isomorphism.
\label{cor: ext}
\end{corollary}

The Grothendieck spectral sequence gives  a \emph{Thom isomorphism}
\[
\eT_{\eC}: H^\bullet(\eC,\Omega_{\eC})\ra H^\bullet_{\eC}(X)[c]:=H^{\bullet+c}_{\eC}(X).
\]
The composition of this morphism with the natural morphism $\boldsymbol{e}_{X,\eC}:H^\bullet_{\eC}(X)\ra H^\bullet(X)$ is the \emph{Gysin map}
\[
\gamma_{\eC}:H^\bullet(\eC,\Bom_{\eC})\ra H^\bullet(X)[c].
\]

\begin{definition} Suppose $\eC\hra X$ is  NN of codimension $c$ in $X$.

\smallskip

\noindent (a) The set  $\eC$   is called \emph{coorientable} in $X$ if the co-orientation sheaf $\Bom_{\eC}$ is isomorphic to the constant sheaf  $\underline{\bZ}$ with stalk $\bZ$ at every point.
 A \emph{co-orientation} of the embedding $\eC\hra X$ is a choice of an isomorphism $\underline{\bZ}\ra \Bom_{\eC}$.  A  coorientation is uniquely determined by an element $\bom_{\eC}\in H^0(\eC,\Bom_{\eC})$ which, viewed as a section of $\Bom_{\eC}$, it  has the property   that  $\bom_{\eC}(w)$ generates the stalk $\Bom_{\eC}(w)$ for any $w\in\eC$.
 
  \noindent (b) For any coorientation $\bom_{\eC}$ of $\eC$ we define by
 \[
 \Phi_{\eC}:=\eT_{\eC}(\bom_{\eC})\in H^c_{\eC}(X),\;\;[\eC]_X: =\gamma_{X,\eC}(\bom_{\eC})\in H^c(X).
 \]

 The class $\Phi_{\eC}$ is called the \emph{Thom class} of the  (normally nonsingular) embedding $\eC\hra X$, and the \emph{element} $[\eC]^{X}$ is called the  \emph{cohomology class}  determined by the normally nonsingular co-oriented embedding $\eC\hra X$.\qed
 \end{definition}

\begin{prop} Suppose $\eC\hra X$ is a $NN$ subset of $X$ of codimension $c$ and $\eS\subset \eC$ is a closed subset  of $\eC$ such that
 \[
 \codim^h_X(\eS)\geq c+2,
 \]
such the $NN$ subset $\eC\setminus X$ is coorientable in $X\setminus \eS$. Then $\eC$ is coorientable in $X$ and any coorientation  of $\eC\setminus \eS$ in $X\setminus \eS$   extends to a coorientation of $\eC$ in $X$.
\label{prop: ext-coor}
\end{prop}
\noindent
\textbf{Proof:} Proposition \ref{prop: ext} gives an isomorphism
\[
H^c_{\eC}(X) \ra H^c_{\eC\setminus \eS}(X\setminus \eS).
\]
which fits in a commutative diagram
\[
\xymatrix{
H^0(\eC,\Bom_{\eC})  \ar^{\tau}[d] \ar[r] & H^0\bigl(\, \eC\setminus\eS,\Bom_{\eC}\mid_{\eC\setminus \eS}\,\bigr)  \ar^{\tau}[d]\\
H^c_{\eC}(X)  \ar[r]  & H^c_{(X\setminus \eS)\cap\eC}(X\setminus \eS)
}\]
where we notice that the restriction  to $\eC\setminus \eS$ of the coorientation sheaf  $\Bom_{\eC}$ is the coorientation sheaf $\Bom_{(X\setminus\eS)/(\eC\setminus \eS)}$. The vertical arrows are the Thom isomorphisms. Since the bottom   horizontal arrow  is an isomorphism, we deduce that  the same is true for the top one. In other words, the  restriction morphism
\[
\Gamma(W,\Bom_{\eC})\ra \Gamma(\eC\setminus \eS,\Bom_{\eC}\mid_{\eC\setminus \eS})
\]
is an isomorphism. 

The coorientation of $\eC\setminus\eS$ in $X\setminus \eS$ determines a section $\bom_{\eC\setminus\eS}$ of $\Bom_{\eC}$ over $\eC\setminus \eS$ such that for every $w\in \eC\setminus \eS$  the element $\bom_{\eC\setminus\eS}(w)$ is a generator of the stalk $\Bom_{\eC}(w)$.   From  the above  diagram we deduce that there exists a unique section ${\bom}_\eC$ of  $\Bom_{\eC}$ over $\eC$ that restricts to $\bom_{\eC\setminus\eS}$. We want to show that for  every $w\in \eC$, the element $\bom_{\eC}(w)$ is a generator of the stalk $\Bom_{\eC}(w)$.       We want  to check this   when $w\in \eS\subset \eC$.

Since $\codim^h_X\eS\geq c+2$ we deduce that  $\eS$ has empty  interior as a subset of $\eC$. This is because if $s$ were a point in the interior of $S$ then $\eH^c_{X/\eS}(s)\cong \eH^c_{\eC}(s)\neq 0$.   Choose  a small, connected open neighborhood $U$ of $w$ in $\eC$ such that  the restriction of $\Bom_{\eC}$ is trivial. This is possible since  the sheaf $\Bom_{\eC}$ is trivial. Note that $U\setminus \eS\neq \emptyset$.

 On the neighborhood $U$  the sections of $\Bom_{\eC}$ can be identified with locally constant functions $U\ra \bZ$. Since $U$ is connected,  any such function must be  constant. Thus $\bom_{\eC}|_U$ can be identified with a constant function  $U\ra \bZ$ whose value   at any point $w'\in U\setminus \eS$ is a generator of $\bZ$. This shows that $\bom_\eC$  is indeed a coorientation of $\eC$ in $X$.
\qed
\smallskip

  The normal nonsingularity is still a pretty strong restriction.  We want to explain how to associate a cohomology class to a closed subset   that slightly violates the normal nonsingularity condition.

  \begin{definition}    A closed subset $\eC\subset X$ is called \textit\textbf{{quasi normally nonsingular}} (for short QNN)\index{quasi normally nonsingular (QNN)} of codimension $c$ if  there exists a  closed subset $\eS\subset\eC$ such that the following hold.

  \begin{itemize}
  
  \item   $\codim^h_X(\eS)\geq c+2$, and
  
   \item $\eC\setminus \eS$  is  a NN closed subset of $X\setminus \eS$ of codimension $c$. 

  \end{itemize}

  The  set $\eS$ is called \emph{a singular locus} for the  embedding $\eC\hra X$.
\qed
\end{definition}

  \begin{remark}  The singular locus in the above definition  \emph{is not uniquely determined}  by $\eC$.  For example if $\eS$ is a singular locus and $w\in \eC\setminus \eS$ then $\eS\cup\{w\}$ is a singular locus if $\dim \eC>1$.
  \end{remark}

   Suppose $\eC\hra X$ is a  QNN subset of $X$ of codimension $c$.  Fix a   singular locus $\eS\subset \eC$.  By Proposition  \ref{prop: ext} and Corollary \ref{cor: ext} we have that $H^k_{\eS}(X)=0$,  $\forall k\leq c+1$, and an isomorphism 
\[
H^c_{\eC}(X)\ra H^c_{(\eC\setminus \eS)}(X\setminus \eS).
\]
We denote by $\eE_{\eS}$ its inverse, and we refer to it as the \emph{extension across $\eS$}.  
If the singular locus $\eS$ is such that $\eC\setminus\eS$ is coorientable in $X\setminus \eS$, then a choice of  coorientation defines   an element $\bom_{\eC,\eS}\in  H^0\bigl(\eC\setminus \eS,\Bom_{(X\setminus\eS)\mid_(\eC\setminus\eS)}\bigr)$. We denote by $\Phi_{X,\eC,\eS}$ the  element  in $H^c(X,X\setminus \eC)$ that corresponds to  $\bom_{\eC,\eS}$ via the isomorphism
 \[
 H^0(\eC\setminus \eS,\Bom_{(X\setminus \eS)/(\eC\setminus \eS)})\stackrel{\eT}{\Lra}H^c_{(\eC\setminus \eS)}(X\setminus \eS)    \stackrel{\eE_\eS}{\Lra} H^c_{\eC}(X).
 \]
\begin{prop} Suppose  $W\hra X$ is a $QNN$ of codimension $c$, such  that  for some choice of singular locus $\eS_0$ the  $NN$ set $\eC\setminus \eS_0$ is coorientable in $X\setminus \eS_0$. Then for any other choice of singular locus $\eS_1$, such that
\[
\codim^h_{X}(\eS_0\cap \eS_1)\geq c+2
\]
 the $NN$ set  $\eC\setminus \eS_1$ is  coorientable in $X\setminus \eS_1$,  and any coorientation  $\bom_0$ of $\eC\setminus \eS_0$ induces  a unique coorientation $\bom_1$ of $\eC\setminus \eS_1$ which agrees with $\bom_0$ on $\eC\setminus(\eS_0\cup\eS_1)$.  Moreover, the class $\Phi_{X,\eC,\eS_0}\in H^c_{\eC}(X)$  determined by $\bom_0$ coincides with the class $\Phi_{X,\eC,\eS_1}\in H^c_{\eC}(X)$ determined by $\bom_1$.
\label{lemma: sing-locus}
\end{prop}
\noindent
\textbf{Proof:}  We carry the proof  in three steps.

\smallskip

\noindent {\bf Step 1.} \emph{If $\eS_0\subset \eS_1$ are singular loci of $\eC$, and $\eC\setminus \eS_1$ is coorientable in $X\setminus \eS_0$, then any  coorientation  of $\eC\setminus \eS_1$ extends to a unique coorientation of $\eC\setminus \eS_0$  and we have an equality between the corresponding elements in $H^c_{\eC}(X)$}  
\[
\Phi_{X,\eC,\eS_0}=\Phi_{X,\eC,\eS_1}.
\]
  Set $X'=X\setminus \eS_0$, $\eC'=\eC\setminus \eS_0$, $\eS'=\eS_1\setminus\eS_0$. Then $\eS_1\setminus\eS_0$ is closed in $\eC\setminus \eS_0$ and 
  \[
  \codim^h_{X\setminus \eS_0}{(\eS_1\setminus \eS_0)}=\codim^h_{X}{\eS_1}\geq c+2.
  \]
   Proposition  \ref{prop: ext-coor} implies that any  coorientation of $\eC'\setminus \eS_1$ extends to a coorientation of $\eC'$ in $X'$.  Moreover, the diagram below is commutative
\[
\xymatrix {
H^c_{\eC\setminus\eS_1}(X\setminus\eS_1) \ar^{\eE_{\eS_1\setminus\eS_0}}[r] \ar^{\eE_{\eS_1}}[d] & H^c_{\eC\setminus \eS_0}(X\setminus \eS_0) \ar^{\eE_{\eS_0}}[dl]\\
H^c_{\eC}(X)
}
\]
We only need to check that $\eE_{\eS'}$ maps the Thom class of the embedding $\eC'\setminus \eS'\hra X'\setminus \eS'$  to the Thom class of the embedding $\eC'\hra X'$. This  follows from the functoriality of the Grothendieck spectral sequence which in this special case can be rephrased as saying that  the Thom isomorphism is compatible with the restriction to open sets.

\noindent {\bf Step 2.} \emph{If $\eS_0$ and $\eS_1$ are singular loci such that $\codim^h_{X}(\eS_0\cap\eS_1)\geq c+2$, then $\eS_0\cup \eS_1$ is a singular locus.} Observe  that for every $x\in X$ we have a  Mayer-Vietoris long exact sequence  (see \cite[Eq. (2.6.29)]{KS})
\[
\dotsc \ra \eH^k_{X/(\eS_0\cap\eS_1)}(x)\ra \eH^k_{X/\eS_0}(x)\oplus \eH^k_{X/\eS_1}(x)\ra \eH^k_{X/(\eS_0\cup\eS_1)}(x)\ra  \eH^{k+1}_{X/(\eS_0\cap\eS_1)}(x)\ra\dotsc,
\]
where for any closed subset $C\subset X$ we denoted by $\eH^k_{X/C}(x)$ the stalk at $x$ of the local cohomology sheaf $\eH^k_{X/C}$.

\noindent {\bf Step 3.}    Suppose  $\eS_0$ and $\eS_1$ are singular loci such that $\eC\setminus \eS_0$ is coorientable and $\codim^h_{X}(\eS_0\cap\eS_1)\geq c+2$. Fix a coorientation $\bom_0$    denote by $\Phi_{X,\eC,\eS_0}$ the element in $H^c_{\eC}(X)$ determined by $\bom_0$. 

The coorientation $\bom_0$  restricts to a coorientation $\bom_{01}$ of $\eC\setminus (\eS_0\cup\eS_1)$  that determines  an element $\Phi_{X,\eC,\eS_0\cup\eS_1}\in H^c_{\eC}(X)$. From Step 1  we deduce
\[
 \Phi_{X,\eC,\eS_0}=\Phi_{X,\eC,\eS_0\cup\eS_1}.
 \]
  By Proposition \ref{prop: ext-coor} the coorientation $\bom_{01}$ extends to a unique coorientation $\bom_1$ of $\eC\setminus \eS_1$  and  from Step 1 we conclude  that the element  $\Phi_{X,\eC,\eS_1}$ determined by $\bom_1$ coincides with the element  $\Phi_{X,\eC,(\eS_0\cup\eS_1)}$ determined by $\bom_{01}$.

\qed
 \begin{definition} (a) If $\eC$ is $QNN$  of codimension $c$ in $X$,  $\eC$ is said to be \emph{coorientable} if there exists a singular locus $\eS\hra \eC$   such that the $NN$ subset $\eC\setminus \eS\hra X\setminus \eS$ is coorientable in $X\setminus \eS$. A coorientation of $\eC$ is  defined to be a coorientation of $\eC\setminus \eS$ in $X\setminus \eS$.\qed
 
 \noindent (b)  If $\eC$ is $QNN$  of codimension $c$ in $X$, then two singular loci $\eS, \eS'$  of $\eC$ are said to be equivalent if there exists a sequence of singular loci
 \[
 \eS=\eS_0,\dotsc, \eS_n=\eS'
 \]
 such that 
 \[
 \codim^h_{X}(\eS_{i-1}\cap \eS_i)\geq c+2, ~~~\forall i=1,\dotsc, n.
 \]
 \end{definition}

 From Proposition \ref{lemma: sing-locus} we deduce that if $\eC$  is a co-oriented QNN subset of codimension $c$ in $X$,  then  the \emph{cohomology} class
\begin{equation}
[\eC]^X:= \boldsymbol{e}_{X,\eC}\bigl(\,\Phi_{X,\eC,\eS}\,\bigr)\in H^c(X)
\label{eq: coho-rep}
\end{equation}
depends only on the equivalence  class  of singular locus $\eS$ and it is called the \emph{cohomology class determined by $\eC$}  and  the (equivalence class of the) singular locus $\eS$. Above $\boldsymbol{e}_{X,\eC}$ denotes the natural extension morphism $H^c_{\eC}(X)\ra H^c(X)$.

The next results follows immediately from the above definitions.

\begin{prop}\label{transversality Banach} Suppose $X,\eY$ are smooth Banach manifolds, $E\subset X$ is a closed subset  such that $\codim^h_{X} E\geq c+ 2$,  $\eC\hra \eY$ is a cooriented QNN of codimension $c$, and $f:X\setminus E\ra \eY$ is a smooth map  with  the following properties.

 \begin{itemize} 
 \item There exists a  singular locus $\eS$ for $\eC$ such that
 \[
 \codim^h_{X}\bigl(\, f^{-1}(\eS)\cup E\,\bigr)\geq c+2.
 \]
 \item The   restriction of $f$ to $X\setminus (f^{-1}(\eS)\cup E)$ is transversal to $\eC\setminus \eS$.
 \end{itemize}

 Then the following hold.
 \begin{enumerate}
 
 \item[(a)]  The subset $f^{-1}(\eC)$ is a cooriented  $QNN$ subspace of $X\setminus E$ of codimension $c$ with singular locus  $f^{-1}(\eS)$.
 
 \item[(b)]  The subset  $f^{-1}(\eC)\cup E$ is a canonically  oriented $QNN$  subspace of codimension $c$ in $X$  with sinylar locus  $f^{-1}(\eS)\cup E$ .
 
 \item[(c)]  The canonical inclusion $i: X\setminus E\ra X$  induces an isomorphism 
 \[
 i^*: H^c(X)\ra H^c(X\setminus E)
 \]
and 
 \[
 i^*[f^{-1}(\eC)\cup E]^X=[f^{-1}(\eC)]^{X\setminus E}=f^*[\eC]^\eY.
 \]
 \end{enumerate}
 \label{prop: blow}
 \end{prop}

 The $QNN$ subspaces  of a Banach  manifold may be difficult to recognize due to the homological codimension conditions.  We recall now a definition (see \ref{def: quasisubman})
    
 \begin{definition} \label{quasi-submanifold}   A \textit{\textbf{quasi-submanifold}} of  $X$ of codimension  $c$ is a closed subset $\eF\subset X$ together with a decreasing filtration by closed subsets
 \[
 \eF=\eF^0\supset\eF^1\supset \eF^2\supset \eF^3\subset \cdots
 \]
 such that  the following hold.
 
 \begin{itemize}
 
 \item $\eF^1=\eF^2$.
 
 \item The \textit {\textbf{strata}} $\eS^k=\eF^k\setminus \eF^{k+1}$,  are submanifolds of $X$ of codimension $k+c$.
 
 \end{itemize}
 
 The    quasi-submanifold is called \textit{\textbf{coorientable}} if $\eS^0$  is coorientable. A \emph{coorientation} of  a quasi-submanifold is then a  coorientation of its top stratum. 
 
 The stratification is said to be finite if there exists an $n$ such that $\eF^{n}=\emptyset$. 
 
 \noindent  (b)    If $f: \eY\ra X$   is a smooth map,  and $\eF$ is a quasi-submanifold of $X$, then $f$ is said to be transversal to $\eF$ if it is transversal to every stratum  of $\eF$.\qed
 
\end{definition}

\begin{prop} \label{vanish codim} Any quasi-submanifold  $\eF=\eF^0\supset\eF^1=\eF^2\supset\cdots $ of codimension $c$  in a Banach manifold  $X$ is a $QNN$ subset of codimension $c$ with singular locus $\eF^2$.
\end{prop}
\noindent
\textbf{Proof:}   It suffices to prove that 
\[
\codim^h_X\eF^0\geq c\;\;\mbox{and}\;\;\codim^h_X\eF^2\geq c+2.
\]
The fact that $\eF^1=\eF^2$ plays no role in the proof of these inequalities so we prove only the first one.

It is enough to show that given $w\in\eF^0$ 
\[ H^k(U,U\setminus \eF^0)=0 ~~\forall k\leq c
\]
for all small open neighborhoods $U$ of $w$. But for $U$ open small enough $U\cap \eF^0$ is a stratified space with a finite stratification because there exists an $n$ such that $w\in \eF_n\setminus\eF_{n+1}$ and $\eF_{n+1}$ is closed. So without restriction of the generality we can suppose the that the stratification is finite. 

  We now use induction on the number of strata and the excision exact sequence for local cohomology to prove the result. Indeed say $\eF^{n+1}=\emptyset$.  Then there exists a maximal $N< n$ such that $\eF_N$ is a nonempty, closed submanifold of codimension $c+N$ in $X$. Therefore $\eF_N$ is normally non-singular and so it has homological codimension at least $c+N$. Suppose we have proved that  $\eF^1$ has $\codim^h\geq c+1$. Then in the long exact sequence:
\[
H^k(X, X\setminus \eF^1)\rightarrow H^k(X, X\setminus \eF^0)\rightarrow H^k(X\setminus \eF^1, X\setminus \eF^0)\rightarrow H^{k+1}(X,X\setminus \eF^1)
\]
the first and the last group are zero for all $k< c$. On the other hand 
\[ H^k(X\setminus \eF^1, X\setminus \eF^0)= H^k(X\setminus \eF^1, (X\setminus \eF^1)\setminus (\eF^0\setminus \eF^1))
\]
Now $\eF^0\setminus \eF^1$ is a closed submanifold of $X\setminus \eF^1$ of codimension $c$ so it has homological codimension at least $c$ so the previous group also vanishes for $k<c$ and this finishes the proof.
\qed

\smallskip

   We summarize the previous discussion.   Any cooriented, codimension $c$ quasi-submanifold $\eF\hra X$ determines 
   \begin {itemize} 
   \item [i)] a   Thom class $\Phi_{\eF}\in H^c_{\eF}(X)$;
   \item[ii)] a cohomology class $[\eF]^X\in H^c(X,\bZ)$. 
   \end{itemize}

It is clear the the preimage of a cooriented quasi-submanifold $\eF\hra X$ of codimension $c$  via a smooth map $F: \eY\ra X$ transversal to the strata of  $\eF$ is a  quasi-submanifold of $\eY$ of codimension $c$ equipped with  a natural coorientation and
\[
\bigl[\,F^{-1}(\eF)\,\bigr]^{\eY}=F^*\bigl(\,[\eF]^{X}\,\bigr).
\]

In finite dimensions the cohomology class  associated to a cooriented quasi-submanifold  class is intimately related to Poincar\'{e} duality.  For any locally compact space $X$ we denote by $H_\bullet^{BM}(X)$ the \emph{Borel-Moore homology}.\index{Borel Moore homology}  In the particular case when $X$ admits a compactification $\overline{X}$ such that the pair $(\overline{X},X)$ is a CW-pair then one can take the definition of the Borel-Moore homology to be 
\[ H_k^{BM}(X):=H_k(\overline{X},\overline{X}\setminus X)
\]
where on the left we mean singular homology. For example if $X$ is a compact differentiable  manifold the Borel-Moore homology coincides with the usual singular homology.

The "classical" Poincare-Alexander duality (see for example Th. 8.3  in \cite{Bredon}) says that for  an oriented manifold  $X$ of dimension $n$ and a compact subset $K\subset X$ the following groups are isomorphic:
\[H^p(K)\simeq H_{n-p}(X,X\setminus K)
\]
where on the left we have \v{C}ech cohomology and on the right we have singular cohomology.  Suppose now we want to switch the role of homology and cohomology in the previous isomorphism. Then  Poincar\'e duality (see Th. IX. 4.7 in \cite{Ive}) has the form:

\[ 
 H_{K}^p(X)=H^p(X, X\setminus K)\simeq H^{BM}_{n-p}(K)
\]
for any closed subset $K\subset X$.

When  $K=X$ Poincar\'e duality takes the form:
\[ H^p(X)=H_{n-p}^{BM}({X})
\]

Let  $X$ be an oriented smooth manifold of dimension $n$,  with orientation class $[X]\in H_n^{BM}(X)$, and   $\eF\hra X$  a  cooriented  quasi-submanifold of $X$ of codimension $c$.   The coorientation of $\eF$ defines an orientation  of the top stratum $\eF^\circ:=\eF\setminus \eF^2$ of $\eF$ and thus a canonical element 
\[
\mu_{\eF^\circ}\in H_{n-c}^{BM}(\eF^\circ).
\]

On the other hand we have:
\begin{prop}\label{extension Borel} Let $\eF$ is a quasi-submanifold of codimension $c$ inside a $n$-dimensional manifold $X$. Then 
\[ H^{BM}_k(\eF)=0~~~ \forall k>n-c
\]
\end{prop}
\noindent
\textbf{Proof:} This follows by induction on strata from the long exact sequence:
\[ H^{BM}_k(\eF^1)\rightarrow H_k^{BM}({\eF})\rightarrow H_k^{BM}({\eF\setminus \eF^1})\rightarrow H^{BM}_{k-1}(\eF^1)
\]
and from the fact that the Borel-Moore homology of a $p$-dimensional vanishes in dimension bigger then $p$.
\qed

\smallskip

We deduce that
\[
H_{n-c}^{BM}(\eF)\ra H_{n-c}^{BM}(\eF^\circ)
\]
is an isomorphism  and thus  there exists   an element $\mu_{\eF}\in H_{n-c}^{BM}(\eF)$ that maps to $\mu_{\eF^\circ}$. If $i$ denotes the canonical inclusion $\eF\hra X$, we obtain an element
\[
[\eF]_{X}:=i_*[\mu_{\eF}]\in H_{n-c}^{BM}(X)
\]
called the \textit{\textbf{(Borel-Moore) homology class}}  determined by the  cooriented    quasi-submanifold $\eF$. 

\begin{prop}\label{BorelMoorePoincare}
The class $\mu_{\eF}\in H_{n-c}^{BM}(\eF)$ is Poincar\'e dual to the Thom class $\Phi_{\eF}\in H^c_{\eF}(X)$ and the class $[\eF]_{X}\in H_{n-c}^{BM}(X)$ is Poincar\'e dual to $[\eF]^{X}\in H^c_{\eF}(X)$.
\end{prop}
\noindent
\textbf{Proof:} see \cite{Ive} Ch. X.4.

The way this relates with the theory of analytic cycles which was used by Nicolaescu, \cite {N1} to construct Poincar\'e duals to the generators of the cohomology ring of $U(n)$ is as follows. Let $X$ be a compact subanalytic manifold and let $\mathscr{F}$ be a quasi-submanifold which is subanalytic. We can choose a triangulation of $\mathscr{X}$ that is compatible with the stratification
\[ X\supset\mathscr{F}\supset\mathscr{F}^2\supset\ldots
\]
After some barycentric subdivisions we can assume that a simplicial neighborhood of $\mathscr{F}^2$  in $\mathscr{F}$ deformation retracts to $\mathscr{F}^2$. In this case we have
\[ H^{BM}_{n-c}(\eF^\circ)=H_{n-c}(\mathscr{F}, \mathscr{F}^2)
\]
The orientation on $\mathscr{F}^\circ$ induces orientations on the top dimensional simplices contained in $\mathscr{F}$. The codimension condition on $\mathscr{F}^\circ$ insures that fact that the sum of these top dimensional simplices with  orientations defines a relative homology class in $H_{n-c}$ as its 
boundary lies in the simplicial neighborhood and the $n-c-1$-homology of this negihborhood is zero.

Hardt has described another model of homology based on subanalytic currents. His theory satisfies the Eilenberg-Steenrod axioms and thus, for any  compact triangulated subanalytic set $X$ we have a canonical isomorphism
\[ H_*^{simplicial}(X)\rightarrow H_*^{Hardt}(X)
\]
Via the above isomorphism the class $[\eF]_{X}$ coincides with the current of integration over $\mathscr{F}^\circ$.

\newpage

\printindex
\end{document}